\documentclass[12pt]{amsart}

\usepackage{enumerate}
\usepackage{epsfig}
\usepackage{amsfonts,amsmath,amssymb,amsthm}
\usepackage{stmaryrd}
\usepackage{wasysym}
\usepackage{xcolor}
\usepackage{xspace} 
\usepackage{hyperref} 
\hypersetup{pdftex, bookmarksnumbered, a4paper, colorlinks=true, linkcolor=darkblue, citecolor=darkblue, urlcolor=blue} 
\usepackage{hypcap} 

\newtheorem{theorem}{Theorem}[]
\newtheorem{proposition}[theorem]{Proposition}
\newtheorem{lemma}[theorem]{Lemma}
\newtheorem{corollary}[theorem]{Corollary}

\newtheorem{definition}[theorem]{Definition}

\theoremstyle{definition} 
\newtheorem{example}[theorem]{Example}

\newtheorem{remark}[theorem]{Remark}

\newtheorem{observation}[theorem]{Observation}
\newtheorem{question}[theorem]{Question}

\newcommand{\cl}{\ensuremath{\prec}} 
\newcommand{\cle}{\ensuremath{\preccurlyeq}} 
\newcommand{\tle}{\ensuremath{\trianglelefteqslant}} 
\newcommand{\vartriangledown}{\ensuremath{\rotatebox[origin=c]{180}{$\triangle$}}} 
\newcommand{\set}[2]{\ensuremath{\left\{#1\,\middle|\,#2\right\}}} 
\newcommand{\eqdef}{\mbox{~\raisebox{0.2ex}{\scriptsize\ensuremath{\mathrm:}}\ensuremath{=} }} 
\newcommand{\ssm}{\ensuremath{\smallsetminus}} 
\newcommand{\dotprod}[2]{\langle#1 \,|\, #2\rangle} 
\newcommand{\R}{\ensuremath{\mathbb{R}}} 
\newcommand{\Z}{\ensuremath{\mathbb{Z}}} 
\newcommand{\N}{\ensuremath{\mathbb{N}}} 
\newcommand{\cM}{\mathcal{M}} 
\newcommand{\cS}{\mathcal{S}} 
\newcommand{\cT}{\mathcal{T}} 
\newcommand{\cC}{\mathcal{C}} 
\newcommand{\cN}{\mathcal{N}} 
\newcommand{\LL}{\ensuremath{\mathbb{L}}} 
\newcommand{\UU}{\ensuremath{\mathbb{U}}} 
\newcommand{\GG}{\ensuremath{\mathbb{G}}} 

\newcommand{\ie}{\textit{i.e.}~} 
\newcommand{\etc}{etc.} 

\newcommand{\gon}[1]{\mbox{$#1$-gon}}
\newcommand{\kcross}[1]{\mbox{$#1$-crossing}}
\newcommand{\ktri}[1]{\mbox{$#1$-trian}\-gu\-la\-tion}
\newcommand{\krel}[1]{\mbox{$#1$-rele}\-vant}
\newcommand{\kbound}[1]{\mbox{$#1$-boun}\-dary}
\newcommand{\kirrel}[1]{\mbox{$#1$-irre}\-le\-vant}
\newcommand{\kstar}[1]{\mbox{$#1$-star}}

\newcommand{\kdepth}[1]{\mbox{$#1$-depth}}
\newcommand{\kkernel}[1]{\mbox{$#1$-kernel}}
\newcommand{\kpointed}[1]{\mbox{$#1$-pointed}}
\newcommand{\kalter}[1]{\mbox{$#1$-alter}nation}

\newcommand{\pt}[1]{\mbox{$#1$-pseu}\-do\-trian\-gu\-la\-tion}
\newcommand{\mpt}{multi\-pseudo\-trian\-gu\-la\-tion\xspace}
\newcommand{\piantiperiodic}{\mbox{$\pi$-anti}\-periodic}
\newcommand{\increasing}[1]{\mbox{$#1$-increa}sing}
\newcommand{\decreasing}[1]{\mbox{$#1$-decrea}sing}
\newcommand{\greedy}[1]{\mbox{$#1$-greedy}}

\newcommand{\svs}{\vspace{.2cm}} 

\graphicspath{{figures/}}
\newcommand{\fref}[1]{Figure~\ref{#1}}
\definecolor{darkblue}{rgb}{0,0,0.7} 
\newcommand{\defn}[1]{\emph{\color{darkblue} #1}}

\begin{document}

\title[Multitr., pseudotr. and primitive sorting networks]{\mbox{Multitriangulations,~pseudotriangulations} \\ and primitive sorting networks}

\author{Vincent Pilaud}
\address{CNRS \& LIX (UMR 7161), \'Ecole Polytechnique}
\email{vincent.pilaud@lix.polytechnique.fr}
\urladdr{http://www.lix.polytechnique.fr/~pilaud/}

\author{Michel Pocchiola}
\address{Institut de Mathématiques de Jussieu (UMR 7586), Universit\'e Pierre et Marie Curie}
\email{pocchiola@math.jussieu.fr}
\urladdr{http://people.math.jussieu.fr/~pocchiola/}

\thanks{
VP was partially supported by grant MTM2008-04699-C03-02 and MTM2011-22792 of the spanish Ministerio de Ciencia e Innovaci\'on, and by European Research Project ExploreMaps (ERC StG 208471). MP was partially supported by the TEOMATRO grant ANR-10-BLAN 0207.\\
\indent An extended abstract of this paper was presented in the \emph{25th European Workshop on Computational Geometry} (Brussels, March 2009)~\cite{pp-mpt-09}. The main results of this paper also appeared in Chapter~$3$ of the first author's PhD dissertation~\cite{p}.
}

\vspace*{-1cm}
\maketitle

\vspace{-.8cm}
\begin{abstract}
We study the set of all pseudoline arrangements with contact points which cover a given support. We define a natural notion of flip between these arrangements and study the graph of these flips. In particular, we provide an enumeration algorithm for arrangements with a given support, based on the properties of certain greedy pseudoline arrangements and on their connection with sorting networks. Both the running time per arrangement and the working space of our algorithm are polynomial.

As the motivation for this work, we provide in this paper a new interpretation of both pseudotriangulations and multitriangulations in terms of pseudoline arrangements on specific supports. This interpretation explains their common properties and leads to a natural definition of \mpt{}s, which generalizes both. We study elementary properties of \mpt{}s and compare them to iterations of pseudotriangulations.

\bigskip
\noindent {\sc Keywords}. pseudoline arrangement $\cdot$ pseudotriangulation $\cdot$ multitriangulation $\cdot$ flip $\cdot$ sorting network $\cdot$ enumeration algorithm
\end{abstract}


\section{Introduction}\label{sec:intro}

The original motivation for this paper is the interpretation of certain families of planar geometric graphs in terms of pseudoline arrangements. As an introductory illustration, we present this interpretation on the family of triangulations of a convex polygon.

\begin{figure}
	\capstart
	\centerline{\includegraphics[scale=1]{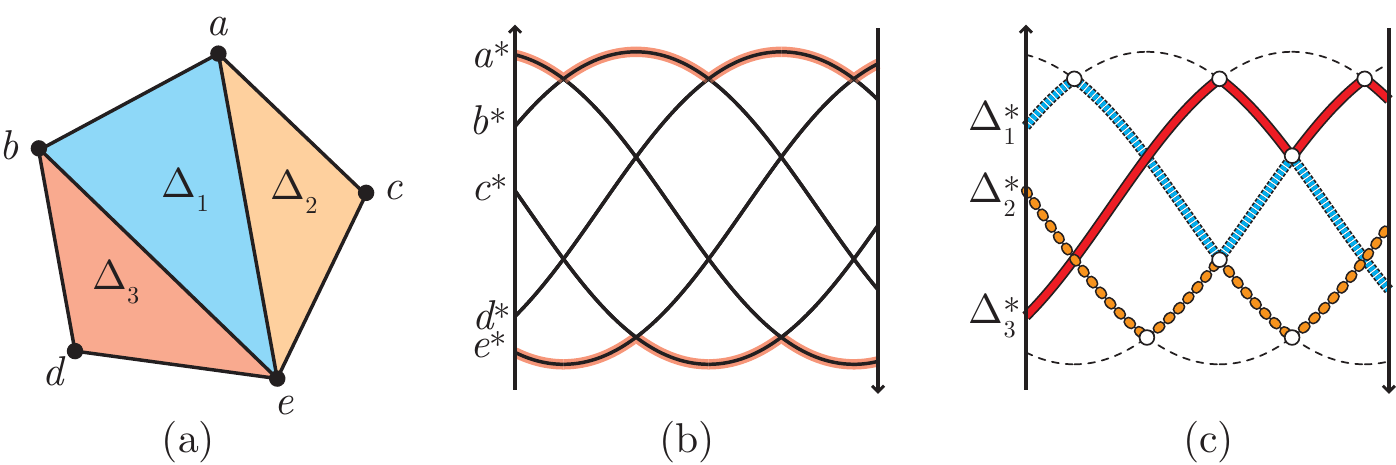}}
	\caption[Dual interpretation of a triangulation]{(a) A triangulation~$T$ of a convex point set~$P$. (b) The dual pseudoline arrangement~$P^*$ of~$P$ (whose first level is highlighted). (c) The dual pseudoline arrangement~$T^*$ of~$T$.}
	\label{fig:triangulation}
\end{figure}

Let~$P$ be a finite point set in convex position in the Euclidean plane~$\R^2$. A \defn{triangulation}~$T$ of~$P$ is a maximal crossing-free set of edges on~$P$, or equivalently, a decomposition of the convex hull of~$P$ into triangles with vertices in~$P$. See \fref{fig:triangulation}(a). A \defn{(strict) bisector} of a triangle~$\Delta$ of~$T$ is a line which passes through a vertex of~$\Delta$ and (strictly) separates its other two vertices. Observe that any triangle has a unique bisector parallel to any direction. Moreover any two triangles of~$T$ have a unique common strict bisector, and possibly an additional non-strict bisector if they share an edge.

We now switch to the line space~$\cM$ of the Euclidean plane~$\R^2$. Remember that~$\cM$ is a M\"obius strip homeomorphic to the quotient space ${\R^2/(\theta,d) \sim (\theta+\pi,-d)}$ via the parametrization of an oriented line by its angle~$\theta$ with the horizontal axis and its algebraic distance~$d$ to the origin (see Section~\ref{sec:duality} for details). For pictures, we represent~$\cM$ as a vertical band whose boundaries are identified in opposite directions. The dual of a point~$p \in P$ is the set of all lines of~$\R^2$ passing through~$p$. It is a \defn{pseudoline} of~$\cM$, \ie a non-separating simple closed curve. The dual of~$P$ is the set~$P^* \eqdef \set{p^*}{p \in P}$. It is a \defn{pseudoline arrangement} of~$\cM$, since any two pseudolines~$p^*,q^*$ of~$P^*$ cross precisely once at the line~$(pq)$. Call \defn{first level} of~$P^*$ the boundary of the external face of the complement of~$P^*$. It corresponds to the supporting lines of the convex hull of~$P$. See \fref{fig:triangulation}(b).

As illustrated on \fref{fig:triangulation}(c), we observe after~\cite{pv-ot-94,pv-ptta-96} that:
\begin{enumerate}[(i)]
\item the set~$\Delta^*$ of all bisectors of a triangle~$\Delta$ of~$T$ is a pseudoline~of~$\cM$;
\item the dual pseudolines~$\Delta_1^*,\Delta_2^*$ of any two triangles~$\Delta_1,\Delta_2$ of~$T$ have a unique crossing point (the unique common strict bisector of~$\Delta_1$ and~$\Delta_2$) and possibly a contact point (when~$\Delta_1$ and~$\Delta_2$ share a common edge);
\item the set~$T^* \eqdef \set{\Delta^*}{\Delta \text{ triangle of } T}$ is a pseudoline arrangement with contact points;~and
\item $T^*$ covers~$P^*$ minus its first level.
\end{enumerate}

We furthermore prove in this paper that this interpretation is bijective: any pseudoline arrangement with contact points supported by~$P^*$ minus its first level is the dual pseudoline arrangement~$T^*$ of a triangulation~$T$ of~$P$.

\medskip
Motivated by this interpretation, we study the set of all \defn{pseudoline arrangements with contact points} which cover a given support in the M\"obius strip. We define a natural notion of \defn{flip} between them, and study the graph of these flips. In particular, we provide an enumeration algorithm for arrangements with a given support (similar to the enumeration algorithm of~\cite{bkps-ceppgfa-06} for pseudotriangulations), based on certain \defn{greedy pseudoline arrangements} and their connection with \defn{primitive sorting networks}~\cite[Section~5.3.4]{k-acp-73}\cite{k-ah-92,b-sms-74}. The running time per arrangement and the working space of our algorithm are both polynomial.

\medskip
We are especially interested in the following particular situation. Let~$L$ be a pseudoline arrangement and $k$ be a positive integer. We call \defn{\pt{k}} of $L$ any pseudoline arrangement with contact points that covers $L$ minus its first $k$ levels. These objects provide dual interpretations for two, until now unrelated, classical generalizations of triangulations of a convex polygon (see Section~\ref{sec:duality} for the definitions and basic properties of these geometric graphs):
\begin{enumerate}[\;\;\;]
\item {\bf Pseudotriangulations} ($k=1$).
Introduced for the study of the visibility complex of a set of disjoint convex obstacles in the plane \cite{pv-tsvcp-96,pv-vc-96}, pseudotriangulations were used in different contexts such as motion planning and rigidity theory \cite{s-ptrmp-05,horsssssw-pmrgpt-05}. Their combinatorial and geometric structure has been extensively studied in the last decade (number of pseudotriangulations \cite{aaks-cmpt-04,aoss-nptcps-08}, polytope of pseudotriangulations \cite{rss-empppt-03}, algorithmic issues \cite{b-eptp-05,bkps-ceppgfa-06,hp-cpcpb-07}, \etc{}). See \cite{rss-pt-06} for a detailed survey on the subject. As far as pseudotriangulations are concerned, this paper has two main applications: it proves the dual characterization of pseudotriangulations in terms of pseudoline arrangements and provides an interpretation of greedy pseudotriangulations in terms of sorting networks, leading to a new proof of the greedy flip property for points~\cite{pv-tsvcp-96,ap-sstvc-03,bkps-ceppgfa-06}. The objects studied in this paper have a further (algorithmic) motivation: as a first step to compute the dual arrangement of a set of disjoint convex bodies described only by its chirotope, Habert and Pocchiola raise in \cite{hp-cpcpb-07} the question to compute efficiently a pseudotriangulation of a pseudotriangulation of the set, \ie a \pt{2}.

\medskip
\item {\bf Multitriangulations} (convex position).
Introduced in the context of extremal theory for geometric graphs~\cite{cp-tttccp-92}, multitriangulations were then studied for their combinatorial structure \cite{n-gdfcp-00,dkkm-2kn..-01,dkm-lahp-02,j-gtdfssp-05}. The study of stars in multitriangulations \cite{ps-mtcsp-09}, generalizing triangles for triangulations, naturally leads to interpret multitriangulations as \mpt{}s of points in convex position. As far as we know, this paper provides the first interpretation of multitriangulations in terms of pseudoline arrangements on the M\"obius strip.
\end{enumerate}

\medskip
The paper is organized as follows. Section~\ref{sec:enumeration} is devoted to the study of all pseudoline arrangements with contact points covering a given support. We define the flip and the greedy pseudoline arrangements~(\ref{subsec:enumeration:flipgraph}) whose properties yield the enumeration algorithm for pseudoline arrangements with a given support (\ref{subsec:enumeration:gfp}).

In Section~\ref{sec:duality}, we prove that the pseudotriangulations of a finite planar point set $P$ in general position correspond to the pseudoline arrangements with contact points supported by the dual pseudoline arrangement of $P$ minus its first level (\ref{subsec:duality:pt}). Similarly, we observe that multitriangulations of a planar point set~$P$ in convex position correspond to pseudoline arrangements with contact points supported by the dual pseudoline arrangement of~$P$ minus its first $k$ levels (\ref{subsec:duality:mt}).

This naturally yields to the definition of \mpt{}s in Section \ref{sec:mpt}. We study the primal of a \mpt{}. We discuss some of its structural properties (\ref{subsec:mpt:pointedcrossing}) which generalize the cases of pseudotriangulations and multitriangulations: number of edges, pointedness, crossing-freeness. We study in particular the faces of \mpt{}s (\ref{subsec:mpt:stars}) which naturally extend triangles in triangulations, pseudotriangles in pseudotriangulations, and stars in multitriangulations.

In Section~\ref{sec:iterated}, we compare \mpt{}s to iterated pseudotriangulations. We give an example of a \ktri{2} which is not a pseudotriangulation of a triangulation (\ref{subsec:iterated:definition}). We prove however that greedy \mpt{}s are iterated greedy pseudotriangulations (\ref{subsec:iterated:greedy}), and we study flips in iterated pseudotriangulations (\ref{subsec:iterated:flips}).

Section~\ref{sec:furthertopics} presents two further topics. The first one (\ref{subsec:furthertopics:horizon}) is a pattern avoiding characterization of greedy \mpt{}s related to horizon trees. The second one (\ref{subsec:furthertopics:dpl}) is a discussion on \mpt{}s of double pseudoline arrangements, which extend pseudotriangulations of convex bodies in the plane. 

Finally, we discuss in Section~\ref{sec:open} some related problems and open questions concerning in particular the primal of a \mpt{}, the diameter and the polytopality of the graph of flips, and the number of \mpt{}s. Since the submission of this paper, some of these questions were partially answered in~\cite{s-npkt-11, ss-mfmpscsp-10, ps-bp-12} based on a framework similar to the material presented in this paper.


\newpage
\section{Pseudoline arrangements with the same support}\label{sec:enumeration}


\subsection{Pseudoline arrangements in the M\"obius strip}\label{subsec:enumeration:pa}

Let~$\cM$ denote the \defn{M\"obius strip} (without boundary), defined as the quotient set of the plane~$\R^2$ under the map $\tau: \R^2\to\R^2,\; (x,y)\mapsto(x+\pi,-y)$. The induced canonical projection will be denoted by~$\pi:\R^2\to\cM$.

A \defn{pseudoline} is the image~$\lambda$ under the canonical projection~$\pi$ of the graph~$\set{(x,f(x))}{x\in\R}$ of a continuous and \defn{\piantiperiodic} function $f:\R\to\R$ (that is, which satisfies~$f(x+\pi)=-f(x)$ for all~$x\in\R$). We say that~$f$ represents the pseudoline~$\lambda$.

When we consider two pseudolines, we always assume that they have a finite number of intersection points. Thus, these intersection points can only be either \defn{crossing points} or \defn{contact points}. See \fref{fig:flip}(a). Any two pseudolines always have an odd number of crossing points (in particular, at least one). When~$\lambda$ and~$\mu$ have exactly one crossing point, we denote it by~$\lambda\wedge\mu$.

\begin{figure}[b]
	\capstart
	\centerline{\includegraphics[scale=1]{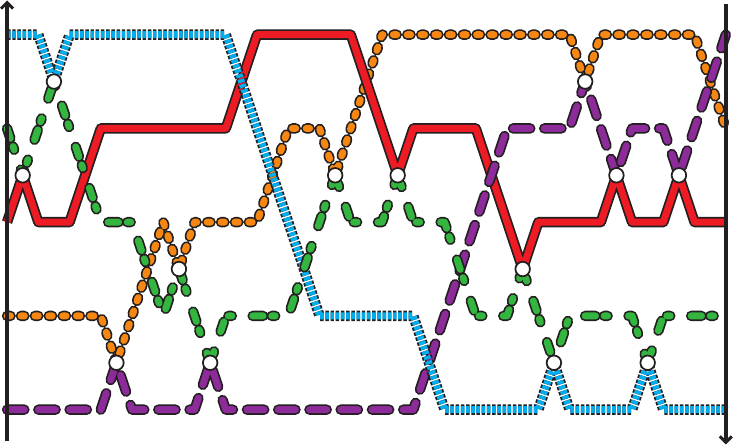}}
	\caption[A pseudoline arrangement in the M\"obius strip]{A pseudoline arrangement in the M\"obius strip. Its contact points are represented by white circles.}
	\label{fig:arrangement}
\end{figure}

A \defn{pseudoline arrangement with contact points} is a finite set~$\Lambda$ of pseudolines such that any two of them have exactly one crossing point and possibly some contact points. See \fref{fig:arrangement}. In this paper, we are only interested in \defn{simple} arrangements, that is, where no three pseudolines meet in a common point. The \defn{support} of~$\Lambda$ is the union of its pseudolines. Observe that~$\Lambda$ is completely determined by its support together with its set of contact points. The \defn{first level} of~$\Lambda$ is the external hull of the support of~$\Lambda$, \ie the boundary of the external face of the complement of the support of~$\Lambda$. We define inductively the \defn{$k$th level} of~$\Lambda$ as the external hull of the support of~$\Lambda$ minus its first $k-1$ levels.

\begin{remark}
The usual definition of pseudoline arrangements does not allow contact points. In this paper, they play a crucial role since we are interested in all pseudoline arrangements which share a common support, and which only differ by their sets of contact points. To simplify the exposition, we omit throughout the paper to specify that we work with pseudoline arrangements \emph{with contact points}.

Pseudoline arrangements are also classically defined on the projective plane rather than the M\"obius strip. The projective plane is obtained from the M\"obius strip by adding a point at infinity.

For more details on pseudoline arrangements, we refer to the broad literature on the topic \cite{g-as-72,bvswz-om-99,k-ah-92,g-pa-97}.
\end{remark}


\subsection{Flip graph and greedy pseudoline arrangements}\label{subsec:enumeration:flipgraph}

\subsubsection{Flips}\label{subsubsec:flipgraph}

\begin{figure}[b]
	\capstart
	\centerline{\includegraphics[scale=1]{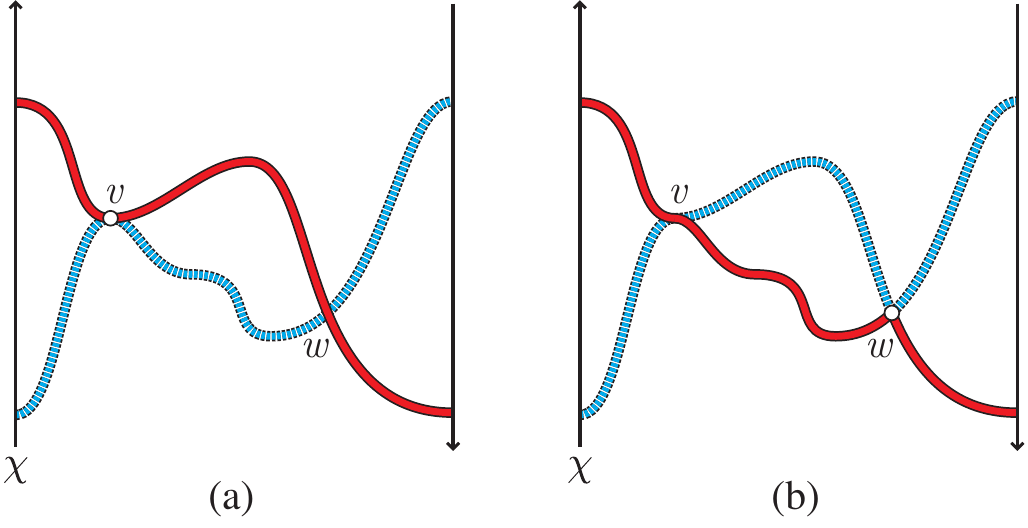}}
	\caption[Flipping a contact point in a pseudoline arrangement]{(a)~A pseudoline arrangement with one contact point~$v$ and one crossing point~$w$. (b)~Flipping~$v$ in the pseudoline arrangement of~(a).}
	\label{fig:flip}
\end{figure}

In the following lemma, we use the symbol~$\triangle$ for the symmetric difference:~$X\triangle Y \eqdef (X\ssm Y)\cup(Y\ssm X)$. We refer to \fref{fig:flip} for an illustration of this lemma.

\begin{lemma}\label{lem:flip}
Let~$\Lambda$ be a pseudoline arrangement,~$\cS$ be its support, and~$V$ be the set of its contact points. Let~$v\in V$ be a contact point of two pseudolines of~$\Lambda$, and~$w$ denote their unique crossing point. Then~$V\triangle\{v,w\}$ is also the set of contact points of a pseudoline arrangement~$\Lambda'$ supported~by~$\cS$.
\end{lemma}

\begin{proof}
Let~$f$ and~$g$ represent the two pseudolines~$\lambda$ and~$\mu$ of~$\Lambda$ in contact at~$v$. Let~$x$~and~$y$ be such that~$v=\pi(x,f(x))$, $w=\pi(y,f(y))$~and $x<y<x+\pi$. We define two functions~$f'$~and~$g'$~by
$$ f' \eqdef 
	\begin{cases}
		f & \text{on } [x,y]+\Z\pi, \\
		g & \text{otherwise,}
	\end{cases}
\quad \text{and} \quad g' \eqdef 
	\begin{cases}
		g & \text{on } [x,y]+\Z\pi, \\
		f & \text{otherwise.}
	\end{cases}
$$
These two functions are continuous and \piantiperiodic, and thus define two pseudolines~$\lambda'$~and~$\mu'$. These two pseudolines have a contact point at~$w$ and a unique crossing point at~$v$, and they cross any pseudoline~$\nu$ of~$\Lambda\ssm\{\lambda,\mu\}$ exactly once (since~$\nu$ either crosses both $\lambda$~and~$\mu$ between $v$~and~$w$, or crosses both $\lambda$~and~$\mu$ between~$w$~and~$v$). Consequently, $\Lambda' \eqdef \Lambda\triangle\{\lambda,\mu,\lambda',\mu'\}$ is a pseudoline arrangement, with support~$\cS$, and contact points~$V\triangle\{v,w\}$.
\end{proof}

\begin{definition}
Let~$\Lambda$ be a pseudoline arrangement with support~$\cS$ and contact points~$V$, let~$v \in V$ be a contact point between two pseudolines of~$\Lambda$ which cross at~$w$, and let~$\Lambda'$ be the pseudoline arrangement with support~$\cS$ and contact points~$V \triangle \{v,w\}$. We say that we obtain~$\Lambda'$ by \defn{flipping}~$v$ in~$\Lambda$.
\end{definition}

Note that the starting point of a flip is always a contact point. To this contact point corresponds precisely one crossing point. In contrast, it would be incorrect to try to flip a crossing point: the two pseudolines which cross at this point might have zero or many contact points.

Observe also that the pseudoline arrangements~$\Lambda$~and~$\Lambda'$ are the only two pseudoline arrangements supported by~$\cS$ whose sets of contact points contain ${V\ssm\{v\}}$.

\begin{definition}
Let~$\cS$ be the support of a pseudoline arrangement. The \defn{flip graph} of~$\cS$, denoted by~$G(\cS)$, is the graph whose vertices are all the pseudoline arrangements supported by~$\cS$, and whose edges are flips between them.
\end{definition}

In other words, there is an edge in the graph~$G(\cS)$ between two pseudoline arrangements if and only if their sets of contact points differ by exactly two points.

Observe that the graph~$G(\cS)$ is regular: there is one edge adjacent to a pseudoline arrangement~$\Lambda$ supported by~$\cS$ for each contact point of~$\Lambda$, and two pseudoline arrangements with the same support have the same number of contact points.

\begin{example}
The flip graph of the support of an arrangement of two pseudolines with~$p$ contact points is the complete graph on~$p+1$~vertices.
\end{example}

\subsubsection{Acyclic orientations}\label{subsubsec:orientations}

Let~$\cS$ be the support of a pseudoline arrangement and~$\bar{\cS}$ denote its preimage under the projection~$\pi$. We orient the graph~$\bar{\cS}$ along the abscissa increasing direction, and the graph~$\cS$ by projecting the orientations of the edges of~$\bar{\cS}$. We denote by~$\cle$ the induced partial order on the vertex set of~$\bar{\cS}$ (defined by~$z\cle z'$ if there exists an oriented path on~$\bar{\cS}$ from~$z$~to~$z'$).

A \defn{filter} of~$\bar{\cS}$ is a proper set~$F$ of vertices of~$\bar{\cS}$ such that~$z\in F$ and~$z\cle z'$ implies~$z'\in F$. The corresponding \defn{antichain} is the set of all edges and faces of~$\bar{\cS}$ with one vertex in~$F$ and one vertex not in~$F$. This antichain has a linear structure, and thus, can be seen as the set of edges and faces that cross a vertical curve~$\bar \chi$ of~$\R^2$. The projection~$\chi \eqdef \pi(\bar \chi)$ of such a curve is called a \defn{cut} of~$\cS$. We see the fundamental domain located between the two curves~$\bar \chi$ and~$\tau(\bar \chi)$ as the result of cutting the M\"obius strip along the cut~$\chi$. For example, we use such a cut to represent pseudoline arrangements in all figures of this paper. See for example \fref{fig:flip}.

The cut~$\chi$ defines a partial order~$\cle_\chi$ on the vertex set of~$\cS$: for all vertices~$v$~and~$w$ of~$\cS$, we write~$v\cle_\chi w$ if there is an oriented  path in~$\cS$ which does not cross~$\chi$. In other words,~$v\cle_\chi w$ if~$\bar v\cle\bar w$, where~$\bar v$ (resp.~$\bar w$) denotes the unique preimage of~$v$ (resp.~$w$) between~$\bar\chi$ and~$\tau(\bar\chi)$. For example, in the arrangements of \fref{fig:flip}, we have~${v\cl_\chi w}$.

Let~$\Lambda$ be a pseudoline arrangement supported by~$\cS$, $v$~be a contact point between two pseudolines of~$\Lambda$ and~$w$~denote their crossing point. Since~$v$~and~$w$ lie on a same pseudoline on~$\cS$, they are comparable for~$\cle_\chi$. We say that the flip of~$v$ is \defn{\increasing{\chi}} if~$v\cl_\chi w$ and \mbox{\defn{\decreasing{\chi}}} otherwise. For example, the flip from~(a) to~(b) in \fref{fig:flip} is \increasing{\chi}. We denote by~$G_\chi(\cS)$ the directed graph of \increasing{\chi} flips on pseudoline arrangements supported by~$\cS$.

\begin{lemma}\label{lem:acyclic}
The directed graph~$G_\chi(\cS)$ of \increasing{\chi} flips is acyclic.
\end{lemma}

\begin{proof}
If~$\Lambda$~and~$\Lambda'$ are two pseudoline arrangements supported by~$\cS$, we write $\Lambda\tle_\chi \Lambda'$ if there exists a bijection~$\phi$ between their sets of contact points such that~$v\cle_\chi \phi(v)$ for any contact point~$v$ of~$\Lambda$. It is easy to see that this relation is a partial order on the vertices of~$G_\chi(\cS)$. Since the edges of~$G_\chi(\cS)$ are oriented according to~$\tle_\chi$, the graph~$G_\chi(\cS)$ is acyclic.
\end{proof}

Theorem~\ref{theo:sorting} below states that this acyclic graph~$G_\chi(\cS)$ has in fact a unique source, and thus is connected.

\subsubsection{Sorting networks}\label{subsubsec:sorting}

\begin{figure}
	\capstart
	\centerline{\includegraphics[scale=1]{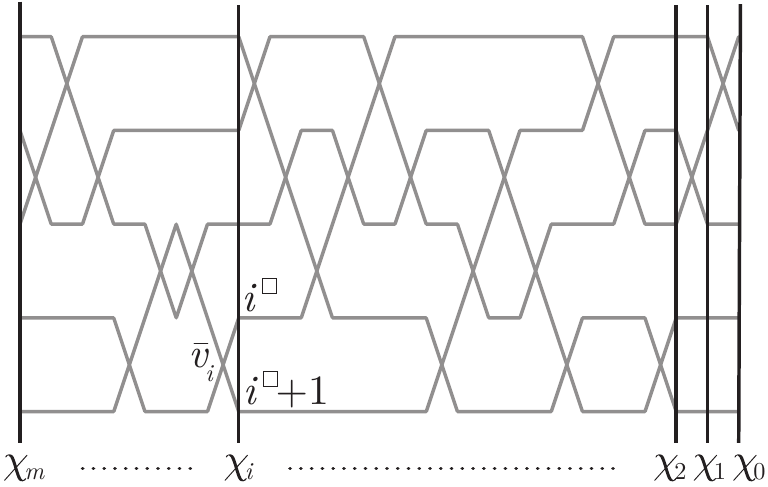}}
	\caption[Sweeping a pseudoline arrangement]{A (backward) sweep of the support of the pseudoline arrangement of \fref{fig:arrangement}.}
	\label{fig:sweep}
\end{figure}

Let~$n$ denote the number of pseudolines of the arrangements supported by~$\cS$ and $m\ge {n \choose 2}$ their number of intersection points (crossing points plus contact points). We consider a chain $F=F_m\supset F_{m-1}\supset \dots\supset F_1\supset F_0=\tau(F)$ of filters of~$\bar{\cS}$ such that two consecutive of them~$F_i$~and~$F_{i+1}$ only differ by a single element:~$\{\bar v_i\} \eqdef F_{i+1}\ssm F_i$. This corresponds to a (backward) \defn{sweep}~$\chi=\chi_0,\chi_1,\dots,\chi_m=\chi$ of the M\"obius strip, where each cut~$\chi_{i+1}$ is obtained from the cut~$\chi_i$ by sweeping a maximal vertex~$v_i \eqdef \pi(\bar v_i)$ of~$\cS$ (for the partial order~$\cle_\chi$). For all~$i$, let~$e_i^1,e_i^2,\dots,e_i^n$ denote the sequence of edges of~$\bar{\cS}$ with exactly one vertex in~$F_i$, ordered from top to bottom, and let~$i^\square$ be the index such that~$\bar v_i$ is the common point of edges~$e_i^{i^\square}$ and~$e_i^{i^\square+1}$~(see \fref{fig:sweep}).

Let~$\Lambda \eqdef \{\lambda_1,\dots,\lambda_n\}$ be a pseudoline arrangement supported by~$\cS$. For all~$i$, we denote by~$\sigma_i^\Lambda$ the permutation of~$\{1,\dots,n\}$ whose $j$th entry~$\sigma_i^\Lambda(j)$ is the index of the pseudoline supporting~$e_i^j$, \ie such that~$\pi(e_i^j)\subset \lambda_{\sigma_i^\Lambda(j)}$. Up to reindexing the pseudolines of~$\Lambda$, we can assume that~$\sigma_0^\Lambda$ is the inverted permutation~$\sigma_0^\Lambda \eqdef [n,n-1,\dots,2,1]$, and consequently that~$\sigma_m^\Lambda$ is the identity permutation. Observe that for all~$i$:
\begin{enumerate}[(i)]
\item if~$v_i$ is a contact point of~$\Lambda$, then~$\sigma_i^\Lambda=\sigma_{i+1}^\Lambda$;
\item otherwise,~$\sigma_{i+1}^\Lambda$~is obtained from~$\sigma_i^\Lambda$ by inverting its $i^\square$th and \mbox{$(i^\square+1)$th} entries.
\end{enumerate}

\begin{figure}
	\capstart
	\centerline{\includegraphics[scale=1]{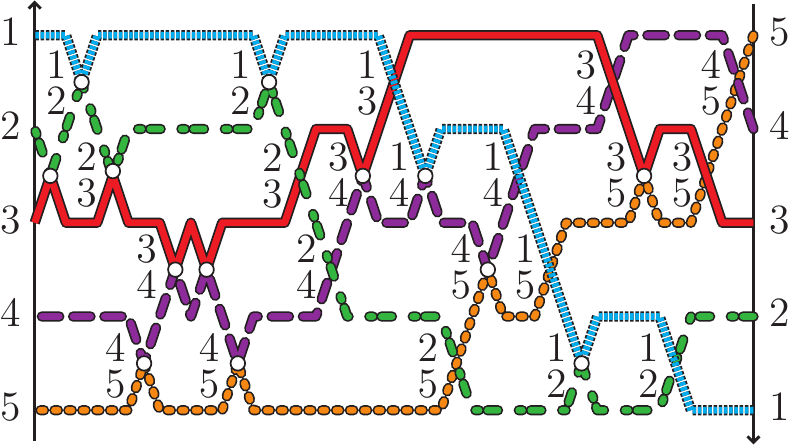}}
	\caption[The greedy pseudoline arrangement obtained by sorting]{The greedy pseudoline arrangement on the support of \fref{fig:arrangement}, obtained by sorting the permutation $[5,4,3,2,1]$. The result of each comparator is written after it (\ie to its left since we sweep backwards).}
	\label{fig:sorting}
\end{figure}

The following theorem is illustrated on \fref{fig:sorting}.

\begin{theorem}\label{theo:sorting}
The directed graph~$G_\chi(\cS)$ has a unique source~$\Gamma$, characterized by the property that the permutation~$\sigma_{i+1}^\Gamma$ is obtained from~$\sigma_i^\Gamma$ by sorting its $i^\square$th and $(i^\square+1)$th entries, for all~$i$.
\end{theorem}

\begin{proof}
If~$\Gamma$ satisfies the above property, then it is obviously a source of the directed graph~$G_\chi(\cS)$: any flip of~$\Gamma$ is \increasing{\chi} since two of its pseudolines cannot touch before they cross.

Assume conversely that~$\Gamma$ is a source of~$G_\chi(\cS)$. Let~$a \eqdef \sigma_i^\Gamma(i^\square)$ and~$b \eqdef \sigma_i^\Gamma(i^\square+1)$. We have two possible situations:
\begin{enumerate}[(i)]
\item If~$a<b$, then the two pseudolines~$\lambda_a$~and~$\lambda_b$ of~$\Gamma$ already cross before~$v_i$. Consequently,~$v_i$ is necessarily a contact point of~$\Gamma$, which implies that~$\sigma_{i+1}^\Gamma(i^\square)=a$ and~$\sigma_{i+1}^\Gamma(i^\square+1)=b$.
\item If~$a>b$, then the two pseudolines~$\lambda_a$~and~$\lambda_b$ of~$\Gamma$ do not cross before~$v_i$. Since~$\Gamma$ is a source of~$G_\chi(\cS)$, $v_i$~is necessarily a crossing point of~$\Gamma$. Thus~$\sigma_{i+1}^\Gamma(i^\square)=b$ and~$\sigma_{i+1}^\Gamma(i^\square+1)=a$.
\end{enumerate}
In both cases,~$\sigma_{i+1}^\Gamma$ is obtained from~$\sigma_i^\Gamma$ by sorting its $i^\square$th and $(i^\square+1)$th entries.
\end{proof}

\begin{corollary}\label{coro:uniquesource}
The graphs of flips~$G(\cS)$ and~$G_\chi(\cS)$ are connected.
\end{corollary}

\begin{definition}
The unique source of the directed graph~$G_\chi(\cS)$ is denoted by~$\Gamma_\chi(\cS)$ and called the \defn{\greedy{\chi} pseudoline arrangement} on~$\cS$.
\end{definition}

Let us reformulate Theorem~\ref{theo:sorting} in terms of sorting networks (see \cite[Section~5.3.4]{k-acp-73} for a detailed presentation; see also~\cite{b-sms-74}). Let~$i<j$ be two integers. A \defn{comparator}~$[i:j]$ transforms a sequence of numbers~$(a_1,\dots,a_n)$ by sorting~$(a_i,a_j)$, \ie replacing~$a_i$ by~$\min(a_i,a_j)$ and~$a_j$ by~$\max(a_i,a_j)$. A comparator~$[i:j]$ is \defn{primitive} if~$j=i+1$. A \defn{sorting network} is a sequence of comparators that sorts any sequence~$(a_1,\dots,a_n)$.

The support~$\cS$ of an arrangement of~$n$ pseudolines together with a sweep $F_m\supset\dots\supset F_0$ corresponds to the primitive sorting network $[1^\square:1^\square+1],\dots,[m^\square:m^\square+1]$ (see~\cite[Section 8]{k-ah-92}). Theorem~\ref{theo:sorting} affirms that sorting the permutation $[n,n-1,\dots,2,1]$ according to this sorting network provides a pseudoline arrangement supported~by~$\cS$, which depends only upon the support~$\cS$ and the filter~$F_0$ (not on the total order given by the sweep).

\subsubsection{Greedy set of contact points}\label{subsubsec:greedycontact}

The following proposition provides an alternative construction of the greedy pseudoline arrangement~$\Gamma_\chi(\cS)$ of the support~$\cS$.

\begin{proposition}
Let~$v_1,\dots,v_q$ be a sequence of vertices of~$\cS$ constructed recursively by choosing, as long as possible, a remaining vertex~$v_i$ of~$\cS$ minimal (for the partial order~$\cle_\chi$) such that~$\{v_1,\dots,v_i\}$ is a subset of the set of contact points of a pseudoline arrangement supported by~$\cS$. Then the resulting set~$\{v_1,\dots,v_q\}$ is exactly the set of contact points of~$\Gamma_\chi(\cS)$.
\end{proposition}

\begin{proof}
First of all,~$\{v_1,\dots,v_q\}$~is by construction the set of contact points of a pseudoline arrangement~$\Lambda$ supported by~$\cS$. If~$\Lambda$ is not the (unique) source~$\Gamma_\chi(\cS)$ of the oriented graph~$G_\chi(\cS)$, then there is a contact point~$v_i$ of~$\Lambda$ whose flip is \decreasing{\chi}. Let~$w$ denote the corresponding crossing point, and~$\Lambda'$ the pseudoline arrangement obtained from~$\Lambda$  by flipping $v_i$. This implies that~$\{v_1,\dots,v_{i-1},w\}$ is a subset of the set of contact points of~$\Lambda'$ and that~$w\cl_\chi v_i$, which contradicts the minimality of~$v_i$.
\end{proof}

Essentially, this proposition affirms that we obtain the same pseudoline arrangement when:
\begin{enumerate}[(i)]
\item sweeping~$\cS$ decreasingly and place crossing points as long as possible; or
\item sweeping~$\cS$ increasingly and place contact points as long as possible.
\end{enumerate}

\subsubsection{Constrained flip graph}\label{subsubsec:constrained}
We now need to extend the previous results to constrained pseudoline arrangements on~$\cS$, in which we force a set~$V$ of vertices of~$\cS$ to be contact points.

\begin{theorem}\label{theo:constrained}
Let~$V$~be a subset of vertices of the support~$\cS$, and let $G_\chi(\cS\,|\,V)$ be the subgraph of~$G_\chi(\cS)$ induced by the pseudoline arrangements with support~$\cS$, whose set of contact points contains~$V$. Then this directed graph~$G_\chi(\cS\,|\,V)$ is either empty or an acyclic connected graph with a unique source~$\Gamma$ characterized by the property that for all~$i$:
\begin{enumerate}[(i)]
\item if~$v_i\in V$, then~$\sigma_{i+1}^\Gamma=\sigma_i^\Gamma$;
\item if~$v_i\notin V$, then~$\sigma_{i+1}^\Gamma$ is obtained from~$\sigma_i^\Gamma$ by sorting its $i^\square$th and $(i^\square+1)$th entries.
\end{enumerate}
\end{theorem}

\begin{proof}
We transform our support~$\cS$ into another one~$\cS'$ by opening all intersection points of~$V$ (the local picture of this transformation is~\includegraphics[scale=1]{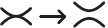}). If~$\cS'$ supports at least one pseudoline arrangement, we apply the result of Theorem~\ref{theo:sorting}: a pseudoline arrangement supported by~$\cS'$ corresponds to a pseudoline arrangement with support~$\cS$ whose set of contact points contains~$V$.
\end{proof}

We denote by~$\Gamma_\chi(\cS\,|\,V)$ the unique source of the constrained flip graph~$G_\chi(\cS\,|\,V)$.

In terms of sorting networks,~$\Gamma_\chi(\cS\,|\,V)$~is the result of the sorting of the inverted permutation $[n,n-1,\dots,2,1]$ by the restricted primitive network~$([i^\square:i^\square+1])_{i\in I}$, where~$I \eqdef \set{i}{v_i\notin V}$.

Observe also that we can obtain, like in the previous subsection, the contact points of~$\Gamma_\chi(\cS\,|\, V)$ by an iterative procedure: we start from the set~$V$ and add recursively  a minimal (for the partial order~$\cle_\chi$) remaining vertex~$v_i$ of~$\cS$ such that~$V\cup\{v_1,\dots,v_i\}$ is a subset of the set of contact points of a pseudoline arrangement supported by~$\cS$. The vertex set produced by this procedure is the set of contact points of the \greedy{\chi} constrained pseudoline arrangement~$\Gamma_\chi(\cS\,|\, V)$.


\subsection{Greedy flip property and enumeration}\label{subsec:enumeration:gfp}

\subsubsection{Greedy flip property}\label{subsubsec:gfp}

\begin{figure}
	\capstart
	\centerline{\includegraphics[scale=1]{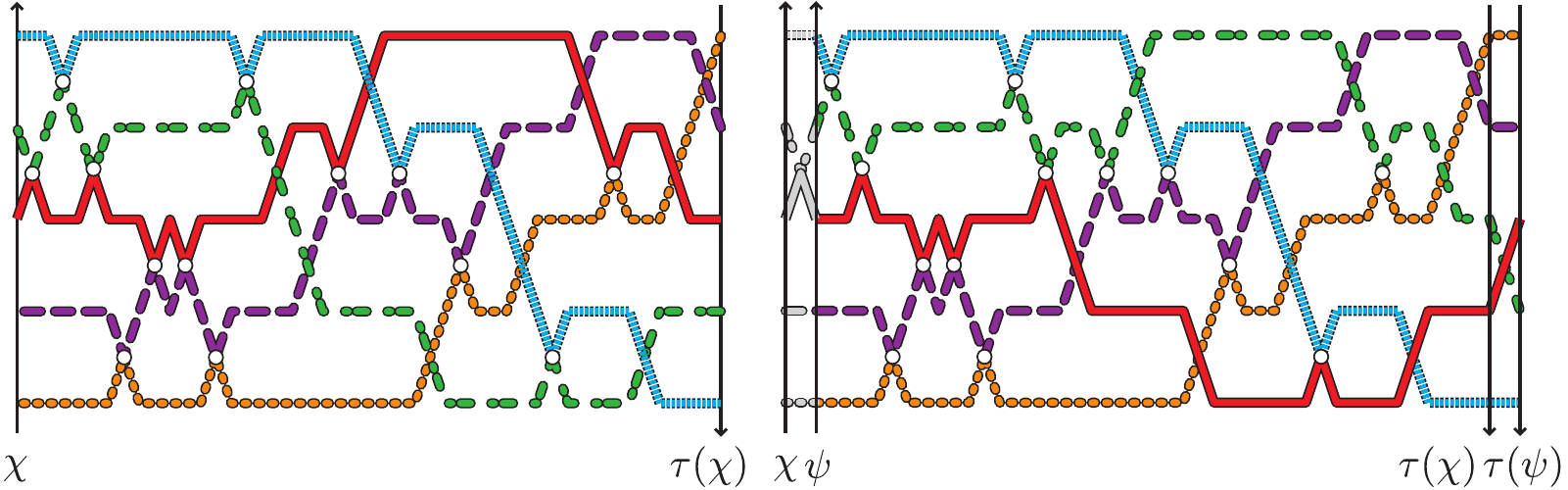}}
	\caption[The greedy flip property]{The greedy flip property.}
	\label{fig:gfp}
\end{figure}

We are now ready to state the \defn{greedy flip property} (see \fref{fig:gfp}) that says how to update the greedy pseudoline arrangement~$\Gamma_\chi(\cS\,|\,V)$ when either~$\chi$ or~$V$ are slightly perturbed.

\begin{theorem}[Greedy flip property]\label{theo:gfp}
Let~$\cS$ be the support of a pseudoline arrangement. Let~$\chi$ be a cut of~$\cS$, $v$~be a minimal (for the order~$\cle_\chi$) vertex of~$\cS$, and~$\psi$~denote the cut obtained from~$\chi$ by sweeping~$v$. Let~$V$ be a set of vertices of~$\cS$ (such that~$G(\cS\,|\,V)$ is not empty), and~${W \eqdef V\cup\{v\}}$. Then:
\begin{enumerate}
\item If~$v$ is a contact point of~$\Gamma_\chi(\cS\,|\,V)$ which is not in~$V$, then $\Gamma_{\psi}(\cS\,|\,V)$ is obtained from $\Gamma_\chi(\cS\,|\,V)$ by flipping~$v$. Otherwise, $\Gamma_{\psi}(\cS\,|\,V)=\Gamma_\chi(\cS\,|\,V)$.
\item If~$v$ is a contact point of~$\Gamma_\chi(\cS\,|\,V)$, then~$\Gamma_{\psi}(\cS\,|\,W)=\Gamma_\chi(\cS\,|\,V)$. Otherwise,~$G(\cS\,|\,W)$ is empty.
\end{enumerate}
\end{theorem}

\begin{proof}
We consider a sweep
$$F_{m+1}=F\supset F_m=F'\supset F_{m-1}\supset\dots\supset F_2\supset F_1=\tau(F)\supset F_0=\tau(F')$$
such that~$F$ (resp.~$F'$) is a filter corresponding to the cut~$\chi$ (resp.~$\psi$). Define the points~${\bar v_i \eqdef F_{i+1}\ssm F_i}$ and ${v_i \eqdef \pi(\bar v_i)}$, and the index~$i^\square$ as previously.
Let $\sigma_1,\dots,\sigma_{m+1}$ denote the sequence of permutations corresponding to~$\Gamma_\chi(\cS\,|\, V)$ on the sweep~$F_1\subset\dots\subset F_{m+1}$. In other words:
\begin{enumerate}[(i)]
\item $\sigma_1$~is the inverted permutation~$[n,n-1,\dots,2,1]$;
\item if~$v_i\in V$, then~$\sigma_{i+1}=\sigma_i$;
\item otherwise,~$\sigma_{i+1}$~is obtained from~$\sigma_i$ by sorting its $i^\square$th and \mbox{$(i^\square+1)$th} entries.
\end{enumerate}
Similarly, let~$\rho_0,\dots,\rho_m$ and $\omega_0,\dots,\omega_m$ denote the sequences of permutations corresponding to~$\Gamma_{\psi}(\cS\,|\ V)$ and~$\Gamma_{\psi}(\cS\,|\ W)$ respectively on the sweep~$F_0\subset\dots\subset F_m$.

Assume first that~$v$ is a contact point of~$\Gamma_\chi(\cS\,|\, V)$, but is not in~$V$. Let~$j$ denote the integer such that~$v_j$ is the crossing point of the two pseudolines of~$\Gamma_\chi(\cS\,|\, V)$ that are in contact at~$v$. We claim that in this case~$\Gamma_{\psi}(\cS\,|\, V)$ is obtained from~$\Gamma_\chi(\cS\,|\, V)$ by flipping~$v$, \ie that:
\begin{enumerate}[(i)]
\item for all~$1\le i\le j$, $\rho_i$ is obtained from~$\sigma_i$ by exchanging~$m^\square$ and~$m^\square+1$;
\item for all~$j<i\le m$,~$\rho_i=\sigma_i$.
\end{enumerate}
Indeed,~$\rho_1$ is obtained by exchanging~$m^\square$ and~$m^\square+1$ in the sequence $\rho_0=[n,n-1,\dots,2,1]=\sigma_1$ (since~$m^\square$ and~$m^\square+1$ are respectively the $(0^\square+1)$th and $0^\square$th entries of~$\rho_0$). Then, any comparison between two consecutive entries give the same result in~$\rho_i$ and in~$\sigma_i$, until~$m^\square$ and~$m^\square+1$ are compared again, \ie until~$i=j$. At this stage,~$m^\square$~and~$m^\square+1$ are already sorted in~$\rho_j$ but not in~$\sigma_j$. Consequently, we have to exchange~$m^\square$~and~$m^\square+1$ in~$\sigma_j$ and not in~$\rho_j$, and we obtain~${\sigma_{j+1}=\rho_{j+1}}$. After this, all the comparisons give the same result in~$\rho_i$~and~$\sigma_i$, and~$\rho_i=\sigma_i$ for all~${j<i\le m}$.

We prove similarly that:
\begin{itemize}
\item When~$v$ is not a crossing point of~$\Gamma_\chi(\cS\,|\, V)$, or is in~$V$, $\rho_i=\sigma_i$~for all~$i\in[m]$, and $\Gamma_{\psi}(\cS\,|\, V)=\Gamma_\chi(\cS\,|\, V)$.
\item When~$v$ is a contact point of~$\Gamma_\chi(\cS\,|\, V)$, $\omega_i=\sigma_i$~for all~$i\in[m]$, and ${\Gamma_{\psi}(\cS\,|\, W)=\Gamma_\chi(\cS\,|\, V)}$.
\end{itemize}

Finally, we prove that~$G(\cS\,|\, W)$ is empty when~$v$ is not a contact point of~$\Gamma_\chi(\cS\,|\, V)$. For this, assume that~$G(\cS\,|\, W)$ is not empty, and consider the greedy arrangement ${\Gamma=\Gamma_\chi(\cS\,|\, W)}$. The flip of any contact point of~$\Gamma$ not in~$W$ is \increasing{\chi}. Furthermore, since~$v$ is a minimal element for~$\cle_\chi$, the flip of~$v$ is also \increasing{\chi}. Consequently,~$\Gamma$ is a source in the graph~$G_\chi(\cS\,|\, V)$, which implies that~$\Gamma_\chi(\cS\,|\, V)=\Gamma$, and thus,~$v$ is a contact point of~$\Gamma_\chi(\cS\,|\, V)$.
\end{proof}

\subsubsection{Enumeration}\label{subsubsec:enumeration}

From the greedy flip property, we derive a binary tree structure on colored pseudoline arrangements supported by~$\cS$, whose left-pending leaves are precisely the pseudoline arrangements supported by~$\cS$. A pseudoline arrangement is \defn{colored} if its contact points are colored in blue, green or red. Green and red contact points are considered to be fixed, while blue ones can be flipped.

\begin{theorem}\label{theo:gfa}
Let~$\cT$ be the binary tree on colored pseudoline arrangements supported by~$\cS$ defined~as~follows:
\begin{enumerate}[(i)]
\item The root of the tree is the \greedy{\chi} pseudoline arrangement on~$\cS$, entirely colored in blue.
\item Any node~$\Lambda$ of~$\cT$ is a leaf of~$\cT$ if either it contains a green contact point or it only contains red contact points.
\item Otherwise, choose a minimal blue point~$v$ of~$\Lambda$. The right child of~$\Lambda$ is obtained by flipping~$v$ and coloring it in blue if the flip is \increasing{\chi} and in green if the flip is \decreasing{\chi}. The left child of~$\Lambda$ is obtained by changing the color of~$v$ into red.
\end{enumerate}
Then the set of pseudoline arrangements supported by~$\cS$ is exactly the set of red-colored leafs of~$\cT$.
\end{theorem}

\begin{proof}
The proof is similar to that of Theorem~$9$ in~\cite{bkps-ceppgfa-06}.

We define inductively a cut~$\chi(\Lambda)$ for each node~$\Lambda$ of~$\cT$: the cut of the root is~$\chi$, and for each node~$\Lambda$ the cut of its children is obtained from~$\chi$ by sweeping the contact point~$v$. We also denote~$V(\Lambda)$ the set of red contact points of~$\Lambda$. With these notations, the greedy flip property (Theorem~\ref{theo:gfp}) ensures that~$\Lambda=\Gamma_{\chi(\Lambda)}(\cS\,|\,V(\Lambda))$, for each node~$\Lambda$ of~$\cT$.

The fact that any red-colored leaf of~$\cT$ is a pseudoline arrangement supported by~$\cS$ is obvious. Reciprocally, let us prove that any pseudoline arrangement supported by~$\cS$ is a red leaf of~$\cT$. Let~$\Lambda$ be a pseudoline arrangement supported by~$\cS$. We define inductively a path~$\Lambda_0,\dots,\Lambda_p$ in the tree~$\cT$ as follows:~$\Lambda_0$~is the root of~$\cT$ and for all~$i\ge 0$,~$\Lambda_{i+1}$~is the left child of~$\Lambda_i$ if the minimal blue contact point of~$\Lambda_i$ is a contact point of~$\Lambda$, and its right child otherwise (we stop when we reach a leaf). We claim that for all~$0\le i\le p$:
\begin{itemize}
\item the set~$V(\Lambda_i)$ is a subset of contact points of~$\Lambda$;
\item the contact points of~$\Lambda$ not in~$V(\Lambda_i)$ are not located between~$\chi(\Lambda)$ and~$\chi$;
\end{itemize}
from which we derive that~$\Lambda=\Lambda_p$ is a red-colored leaf.
\end{proof}

Visiting the tree~$\cT$ provides an algorithm to enumerate all pseudoline arrangements with a given support. In the next section, we will see the connection between this algorithm and the enumeration algorithm of~\cite{bkps-ceppgfa-06} for pseudotriangulations of a point~set.

Let us briefly discuss the complexity of this algorithm. We assume that the input of the algorithm is a pseudoline arrangement and we consider a flip as an elementary operation. Then, this algorithm requires a polynomial running time per pseudoline arrangement supported by~$\cS$. As for many enumeration algorithms, the crucial point of this algorithm is that its working space is also polynomial (while the number of pseudoline arrangements supported by~$\cS$ is exponential).


\section{Dual pseudoline arrangements}\label{sec:duality}

In this section, we prove that both the graph of flips on ``(pointed) pseudotriangulations of a point set'' and the graph of flips on ``multitriangulations of a convex polygon'' can be interpreted as graphs of flips on ``pseudoline arrangements with a given support''. This interpretation is based on the classical duality that we briefly recall in the first subsection, and leads to a natural definition of ``\mpt{}s of a pseudoline arrangement'' that we present in Section~\ref{sec:mpt}.


\subsection{Dual pseudoline arrangement of a point set}\label{subsec:duality:points}

To a given oriented line in the Euclidean plane, we associate its angle~${\theta\in\R/2\pi\Z}$ with the horizontal axis and its algebraic distance~$d\in\R$ to the origin (\ie the value of~$\dotprod{(-v,u)}{.}$ on the line, where~$(u,v)$ is its unitary direction vector). Since the same line oriented in the other direction gives an angle~$\theta+\pi$ and a distance~$-d$, this parametrization naturally associates a point of the M\"obius strip~${\cM \eqdef \R^2/(\theta,d) \sim (\theta+\pi,-d)}$ to each line of the Euclidean plane. In other words, the line space of the Euclidean plane is (isomorphic to) the M\"obius strip.

\begin{figure}
	\capstart
	\centerline{\includegraphics[scale=1]{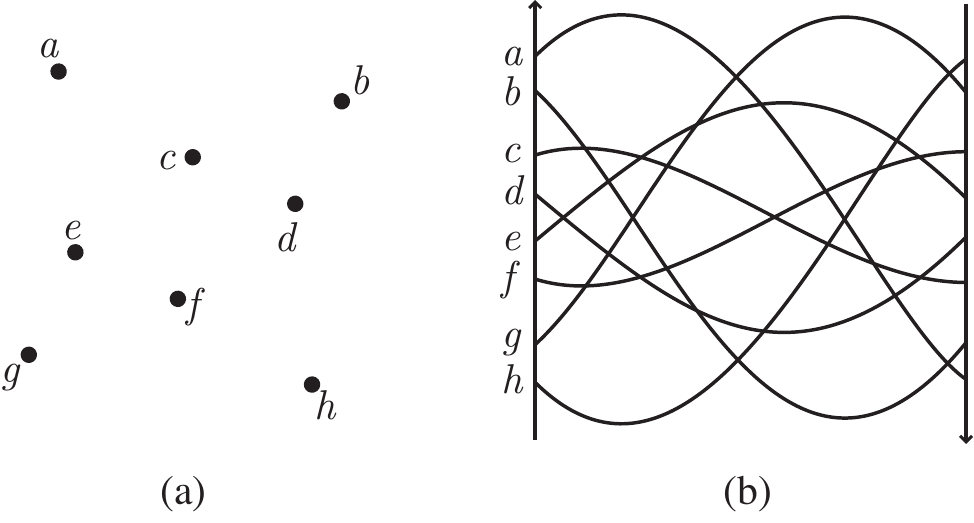}}
	\caption[Dual arrangement of a point set]{A point set~$P$ in general position (a) and (a representation of) its dual arrangement~$P^*$~(b).}
	\label{fig:dual}
\end{figure}

Via this parametrization, the set of lines passing through a point~$p$ forms a pseudoline~$p^*$. The pseudolines~$p^*$~and~$q^*$ dual to two distinct points~$p$~and~$q$ have a unique crossing point, namely the line~$(pq)$. Thus, for a finite point set~$P$ in the Euclidean plane, the set~$P^* \eqdef \set{p^*}{p\in P}$ is a pseudoline arrangement without contact points (see \fref{fig:dual}). Again, we always assume that the point set~$P$ is in general position (no three points lie in a same line), so that the arrangement~$P^*$ is simple (no three pseudolines pass through the same point).

This elementary duality also holds for any topological plane (or \mbox{$\R^2$-plane}, see~\cite{sbghls-cpp-95} for a definition), not only for the Euclidean plane~$\R^2$. That is to say, the line space of a topological plane is (isomorphic to) the M\"obius strip and the dual of a finite set of points in a topological plane is a pseudoline arrangement without contact points. Let us also recall that any pseudoline arrangement of the M\"obius strip without contact points is the dual arrangement of a finite set of points in a certain topological plane~\cite{hp-adp-08}. Thus, in the rest of this paper, we deal with sets of points and their duals without restriction to the Euclidean plane.


\subsection{Dual pseudoline arrangement of a pseudotriangulation}\label{subsec:duality:pt}

We refer to~\cite{rss-pt-06} for a detailed survey on pseudotriangulations, and just recall here some basic definitions.

\begin{definition}\label{def:pseudotriangulation}
A \defn{pseudotriangle} is a polygon~$\Delta$ with only three convex vertices (the \defn{corners} of~$\Delta$), joined by three concave polygonal chains (\fref{fig:pseudotriangle}). A line is said to be \defn{tangent} to~$\Delta$ if:
\begin{enumerate}[(i)]
\item either it passes through a corner of~$\Delta$ and separates the two edges incident to it;
\item or it passes through a concave vertex of~$\Delta$ and does not separate the two edges incident~to~it.
\end{enumerate}

\begin{figure}[b]
	\capstart
	\centerline{\includegraphics[scale=1.2]{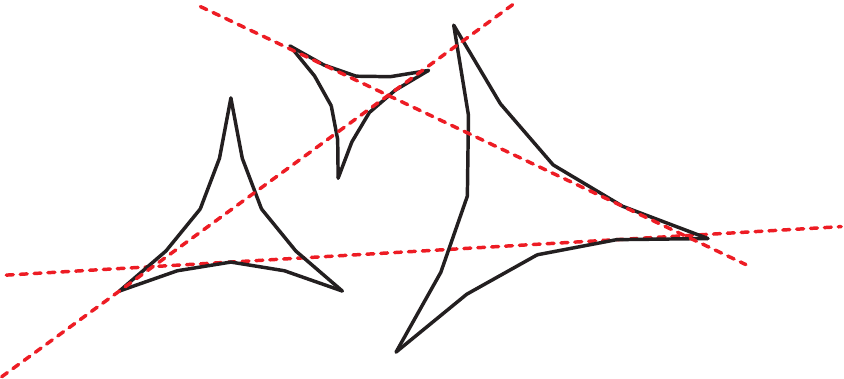}}
	\caption[Three pseudotriangles and their common tangents]{Three pseudotriangles and their common tangents.}
	\label{fig:pseudotriangle}
\end{figure}

A \defn{pseudotriangulation} of a point set~$P$ in general position is a set of edges of~$P$ which decomposes the convex hull of~$P$ into pseudotriangles. We moreover always assume that pseudotriangulations are \defn{pointed}, meaning that there exists a line passing through any point~$p \in P$ and defining a half-plane containing all the edges incident to~$p$.
\end{definition}

The results of this paper only concern pointed pseudotriangulations. Therefore we omit to always specify that pseudotriangulations are pointed. Historically, pseudotriangulations were introduced for families of smooth convex bodies~\cite{pv-vc-96} and were therefore automatically pointed. Pseudotriangulations of points, pointed or not, can be regarded as limits of pseudotriangulations of infinitesimally small convex bodies. Note that pointed pseudotriangulations are edge-minimal pseudotriangulations.

Under the pointedness assumption, any two pseudotriangles of a pseudotriangulation have a unique common tangent. This leads to the following observation (see \fref{fig:pseudotriangulationpoints}):

\begin{observation}[\cite{pv-ot-94,pv-ptta-96}]\label{observation:pseudotriangulations}
Let~$T$ be a pseudotriangulation of a point set~$P$ in general position.~Then:
\begin{enumerate}[(i)]
\item the set~$\Delta^*$ of all tangents to a pseudotriangle~$\Delta$ of~$T$ is a \mbox{pseudoline};
\item the dual pseudolines~$\Delta_1^*, \Delta_2^*$ of any two pseudotriangles~$\Delta_1,\Delta_2$ of~$T$ have a unique crossing point (the unique common tangent to~$\Delta_1$ and~$\Delta_2$) and possibly a contact point (when~$\Delta_1$ and~$\Delta_2$ share a common edge);
\item the set~$T^* \eqdef \set{\Delta^*}{\Delta\text{ pseudotriangle of } T}$ is a pseudoline arrangement (with contact points);~and
\item $T^*$ is supported by~$P^*$ minus its first level (see \fref{fig:pseudotriangulationpoints}(b)).
\end{enumerate}
\begin{figure}
	\capstart
	\centerline{\includegraphics[scale=1]{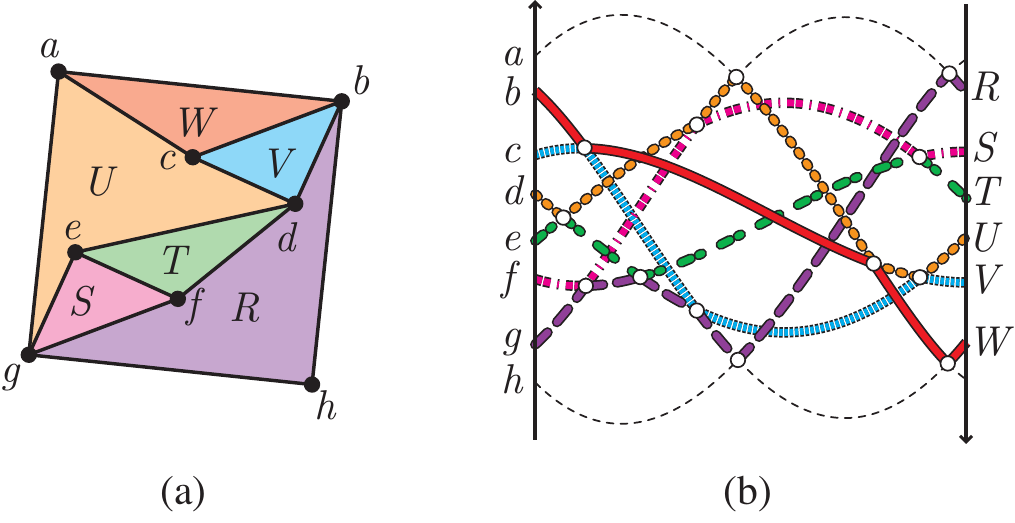}}
	\caption[A pseudotriangulation of a point set and its dual pseudoline arrangement]{(a) A pseudotriangulation~$T$ of the point set~$P$ of \fref{fig:dual}(a). (b) The dual arrangement~$T^*$ of~$T$, drawn on the dual arrangement~$P^*$ of~$P$ of \fref{fig:dual}(b). Each pseudoline of~$T^*$ corresponds to a pseudotriangle of~$T$; each contact point in~$T^*$ corresponds to an edge in~$T$; each crossing point in~$T^*$ corresponds to a common tangent in~$T$.}
	\label{fig:pseudotriangulationpoints}
\end{figure}
\end{observation}

In fact, this covering property characterizes pseudotriangulations:

\begin{theorem}\label{theo:dualitypt}
Let~$P$ be a finite point set in general position in the plane, and~$P^{*1}$ denote the support of its dual pseudoline arrangement minus its first level. Then:
\begin{enumerate}[(i)]
\item The dual arrangement~$T^* \eqdef \set{\Delta^*}{\Delta\text{ pseudotriangle of }T}$ of a pseudotriangulation~$T$ of~$P$ is supported by~$P^{*1}$.
\item The primal set of edges $$\quad E \eqdef \set{[p,q]}{p,q\in P,\; p^*\wedge q^* \text{ contact point of } \Lambda}$$ of a pseudoline arrangement~$\Lambda$ supported by~$P^{*1}$ is a pseudotriangulation of~$P$.
\end{enumerate}
\end{theorem}

\begin{figure}
	\capstart
	\centerline{\includegraphics[scale=1]{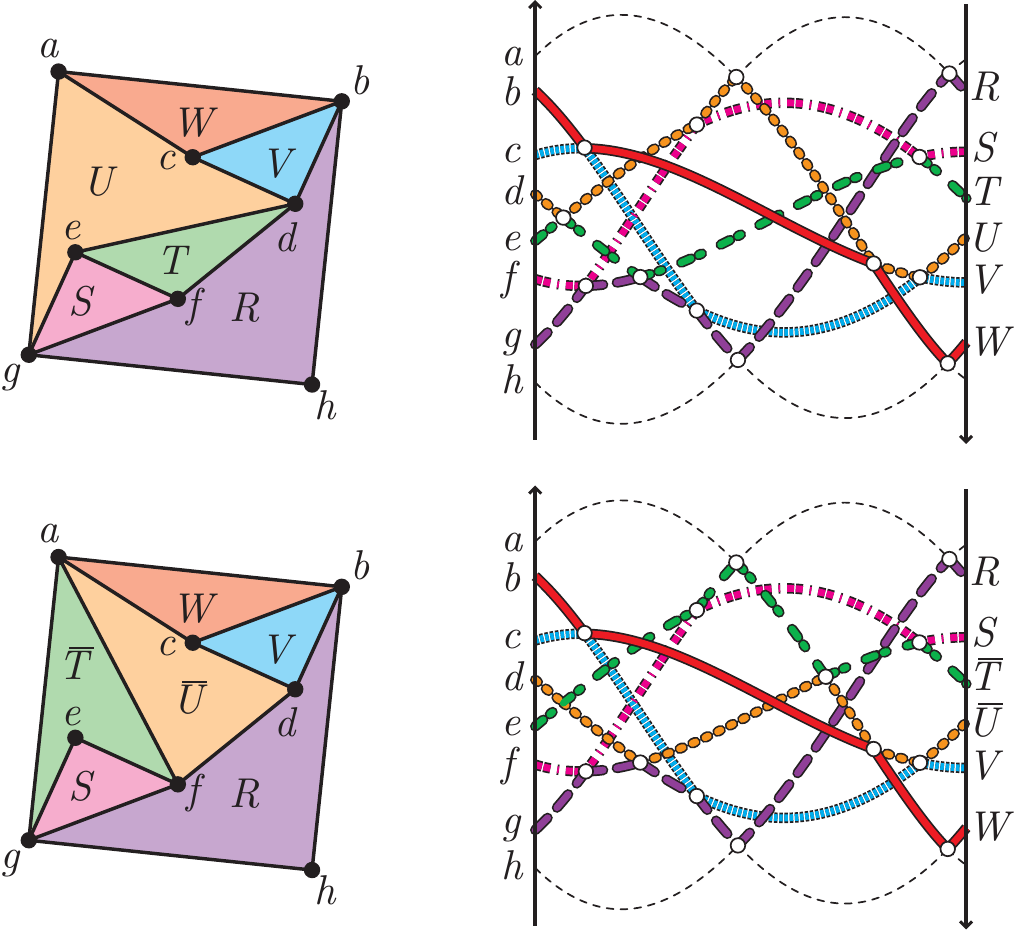}}
	\caption[A flip in a pseudotriangulation]{A flip in the pseudotriangulation of \fref{fig:pseudotriangulationpoints}(a) and the corresponding flip in its dual pseudoline arrangement of \fref{fig:pseudotriangulationpoints}(b).}
	\label{fig:pseudotriangulationpointsflip}
\end{figure}

In this section, we provide three proofs of Part~(ii) of this result. The first proof is based on flips. First, remember that there is also a simple flip operation on pseudotriangulations of~$P$: replacing any internal edge~$e$ in a pseudotriangulation of~$P$ by the common tangent of the two pseudotriangles containing~$e$ produces a new pseudotriangulation of~$P$. For example, \fref{fig:pseudotriangulationpointsflip} shows two pseudotriangulations of the point set of \fref{fig:dual}(a), related by a flip, together with their dual pseudoline arrangements. We denote by~$G(P)$ the graph of flips on pseudotriangulations of~$P$, whose vertices are pseudotriangulations of~$P$ and whose edges are flips between them. In other words, there is an edge in~$G(P)$ between two pseudotriangulations of~$P$ if and only if their symmetric difference is reduced to a pair.

\begin{proof}[Proof~$1$ of Theorem~\ref{theo:dualitypt}(ii)]
The two notions of flips (the primal notion on pseudotriangulations of~$P$ and the dual notion on pseudoline arrangements supported by~$P^{*1}$) coincide via duality: an internal edge~$e$ of a pseudotriangulation~$T$ of~$P$ corresponds to a contact point~$e^*$ of the dual pseudoline arrangement~$T^*$; the two pseudotriangles~$\Delta_1$~and~$\Delta_2$ of~$T$ containing~$e$ correspond to the two pseudolines~$\Delta_1^*$~and~$\Delta_2^*$ of~$T^*$ in contact at~$e^*$; and the common tangent~$f$ of~$\Delta_1$~and~$\Delta_2$ corresponds to the crossing point~$f^*$ of~$\Delta_1^*$~and~$\Delta_2^*$.

Thus, the graph~$G(P)$ is a subgraph of~$G(\cS)$. Since both are connected and regular of degree~$|P|-3$, they coincide. In particular, any pseudoline arrangement supported by~$P^{*1}$ is the dual of a pseudotriangulation of~$P$.
\end{proof}

\begin{remark}
Observe that this duality matches our greedy pseudoline arrangement supported by~$P^{*1}$ with the greedy pseudotriangulation of~\cite{bkps-ceppgfa-06}.
In particular, the greedy flip property and the enumeration algorithm of Subsection~\ref{subsec:enumeration:gfp} are generalizations of results in~\cite{bkps-ceppgfa-06}.
\end{remark}

Our second proof of Theorem~\ref{theo:dualitypt} is slightly longer but more direct, and it introduces a ``witness method'' that we will repeatedly use throughout this paper. It is based on the following characterization of pseudotriangulations:

\begin{lemma}[\cite{s-ptrmp-05}]\label{lem:streinu}
A graph~$T$ on~$P$ is a pointed pseudotriangulation of~$P$ if and only if it is crossing-free, pointed and has~${2|P|-3}$ edges.
\end{lemma}

\begin{proof}[Proof~$2$ of Theorem~\ref{theo:dualitypt}(ii)]
We check that~$E$ is crossing-free, pointed and has~$2|P|-3$ edges:

\smallskip
\paragraph{\bf Cardinality.} First, the number of edges of~$E$ equals the difference between the number of crossing points of~$P^*$ and of~$\Lambda$:
$$|E|={\left|P^*\right| \choose 2}-{| \Lambda| \choose 2}={|P| \choose 2}-{|P|-2 \choose 2}=2|P|-3.$$

\begin{figure}[b]
	\capstart
	\centerline{\includegraphics[scale=1]{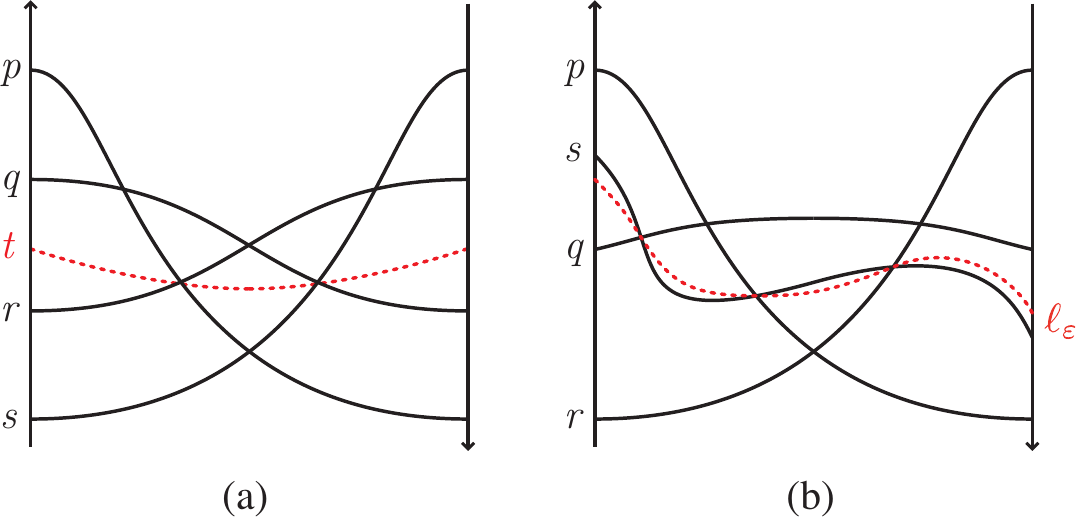}}
	\caption[Dual arrangements of forbidden configurations]{(a) Four points~$p,q,r,s$ in convex position with the intersection~$t$ of~$[p,r]$ and~$[q,s]$. (b)~A point~$s$ inside a triangle~$pqr$ with the witness pseudoline~$\ell_\varepsilon$.}
	\label{fig:pqrs}
\end{figure}

\paragraph{\bf Crossing-free.} Let~$p,q,r,s$ be four points of~$P$ in convex position. Let~$t$ be the intersection of~$[p,r]$~and~$[q,s]$ (see \fref{fig:pqrs}(a)). We use the pseudoline~$t^*$ as a \defn{witness} to prove that~$[p,r]$ and~$[q,s]$ cannot both be in~$E$. For this, we count crossings of~$t^*$ with~$P^*$ and~$\Lambda$ respectively:
\begin{enumerate}[(i)]
\item Since~$P$ is in general position, the point $t$~is not in~$P$. Therefore ${P^*\cup\{t^*\}=(P\cup\{t\})^*}$ is a (non-simple) pseudoline arrangement, and~$t^*$ crosses~$P^*$ exactly~$|P|$ times.
\item Since~$t^*$ is a pseudoline, it crosses each pseudoline of~$\Lambda$ at least once. Thus, it crosses~$\Lambda$ at least~$|\Lambda|=|P|-2$ times.
\item For each of the points~$p^*\wedge r^*$~and~$q^*\wedge s^*$, replacing the crossing point by a contact point removes two crossings with~$t^*$.
\end{enumerate}
Thus,~$[p,r]$~and~$[q,s]$ cannot both be in~$E$, and~$E$ is crossing-free.

\medskip
\paragraph{\bf Pointed.} Let~$p,q,r,s$ be four points of~$P$ such that~$s$ lies inside the convex hull of~$\{p,q,r\}$. We first construct a witness pseudoline (see \fref{fig:pqrs}(b)) that we use to prove that~$[p,s]$,~$[q,s]$~and~$[r,s]$ cannot all be in~$E$. Let~$f_p,f_q,f_r$~and~$f_s$ represent~$p^*,q^*,r^*$~and~$s^*$ respectively. Let ${x,y,z\in\R}$ be such that~$f_p(x)=f_s(x)$, $f_q(y)=f_s(y)$~and~$f_r(z)=f_s(z)$. Let~$g$ be a continuous and \piantiperiodic{} function vanishing exactly on ${\{x,y,z\}+\Z\pi}$ and changing sign each time it vanishes; say for example $g(t) \eqdef \sin(t-x)\sin(t-y)\sin(t-z)$. For all~$\varepsilon> 0$, we define the function $h_\varepsilon:\R\to\R$ by~$h_\varepsilon(t)=f_s(t)+\varepsilon g(t)$. The function~$h_\varepsilon$ is continuous and \piantiperiodic. The corresponding pseudoline~$\ell_\varepsilon$ crosses~$s^*$ three times. It is also easy to see that if~$\varepsilon$ is sufficiently small, then~$\ell_\varepsilon$ crosses the pseudolines of~$(P\ssm\{s\})^*$ exactly as~$s^*$ does (see \fref{fig:pqrs}(b)). For such a small~$\varepsilon$, we count the crossings of~$\ell_\varepsilon$ with~$P^*$~and~$\Lambda$ respectively:
\begin{enumerate}[(i)]
\item $\ell_\varepsilon$~crosses~$P^*$ exactly~$|P|+2$ times (it crosses~$s^*$ three times and any other pseudoline of~$P^*$ exactly once).
\item Since~$\ell_\varepsilon$ is a pseudoline, it crosses~$\Lambda$ at least~$|\Lambda|=|P|-2$ times.
\item For each of the points~$p^*\wedge s^*$, $q^*\wedge s^*$~and~$r^*\wedge s^*$, replacing the crossing point by a contact point removes two crossings with~$\ell_\varepsilon$.
\end{enumerate}
Thus,~$[p,r]$,~$[q,s]$~and~$[r,s]$ cannot all be in~$E$, and~$E$ is pointed.
\end{proof}

Observe that once we know that~$E$ is crossing-free, we could also argue its pointedness observing that the pseudolines of~$\Lambda$ would cover the dual pseudoline of a non-pointed vertex twice.

\smallskip

Our third proof of Theorem~\ref{theo:dualitypt} focusses on pseudotriangles. For every pseudoline~$\lambda$ of the pseudoline arrangement~$\Lambda$, we denote by $${S(\lambda) \eqdef \set{[p,q]}{p,q\in P,\; p^*\wedge q^*\text{ contact point of }\lambda}}$$ the polygonal cycle formed by the edges primal to the contact points of~$\lambda$. For any point~$q$ in the plane, we denote by~$\sigma_\lambda(q)$ the \defn{winding number} of~$S(\lambda)$ around~$q$, \ie the number of rounds made by~$S(\lambda)$ around the point~$q$.

\begin{proof}[Proof~$3$ of Theorem~\ref{theo:dualitypt}(ii)]
Consider a point~$q$ inside the convex hull of our point set~$P$, and such that~$P\cup\{q\}$ be in general position. Hence, its dual pseudoline~$q^*$ has exactly~$|P|$ crossings with~$P^*$, none of which are on the first level of~$P^*$. For any pseudoline~$\lambda\in\Lambda$, let~$\tau_\lambda(q)$ denote the number of intersection points between~$q^*$ and~$\lambda$ (that is, the number of tangents to~$S(\lambda)$ passing through~$q$). Then~$\sigma_\lambda(q)=(\tau_\lambda(q)-1)/2$ and
$$|P|=|q^*\cap P^{*1}|=\sum_{\lambda\in\Lambda} \tau_\lambda(q)=|\Lambda|+2\sum_{\lambda\in\Lambda} \sigma_\lambda(q)=|P|-2+2\sum_{\lambda\in\Lambda} \sigma_\lambda(q).$$
Consequently, $\sum_{\lambda\in\Lambda} \sigma_\lambda(q)=1$. Hence $\sigma_\lambda(q)=0$ for all~$\lambda\in\Lambda$, except for precisely one pseudoline~$\mu\in\Lambda$ which satisfies $\sigma_\mu(q)=1$.

As a consequence, all polygons $S(\lambda)$ for $\lambda\in\Lambda$ are pseudotriangles (otherwise, we would have points such that $\sigma_\lambda(q)>1$ for some $\lambda$) and they cover the convex hull of~$P$. Consequently, these $|P|-2$ pseudotriangles form a pointed pseudotriangulation of~$P$.
\end{proof}


\subsection{Dual pseudoline arrangement of a multitriangulation}\label{subsec:duality:mt}

\begin{figure}
	\capstart
	\centerline{\includegraphics[scale=1]{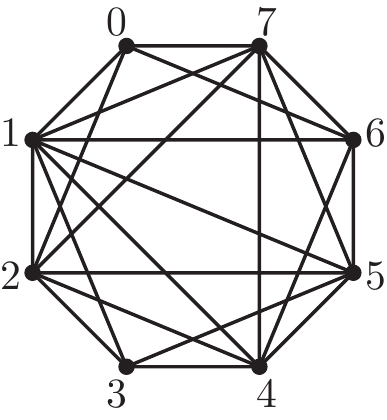}}
	\caption[A multitriangulation]{A \ktri{2} of the octagon.}
	\label{fig:2triang8points}
\end{figure}

Let $C_n$ denote the set of vertices of the convex regular \gon{n}. We are interested in the following generalization of triangulations, introduced by Capoyleas and Pach \cite{cp-tttccp-92} in the context of extremal theory for geometric graphs (see \fref{fig:2triang8points}).

\begin{definition}
For $\ell\in\N$, an \defn{\kcross{\ell}} is a set of $\ell$ mutually crossing edges of $C_n$.
A \defn{\ktri{k}} of the \gon{n} is a maximal set of edges of $C_n$ with no \kcross{(k+1)}.
\end{definition}

Observe that an edge of $C_n$ can be involved in a \kcross{(k+1)} only if there remain at least $k$ vertices on each side. Such an edge is called \defn{\krel{k}}. An edge with exactly  (resp.~strictly less than) $k-1$ vertices on one side is a \defn{\kbound{k}} edge (resp.~a \defn{\kirrel{k}} edge).
By maximality, every \ktri{k} consists of all the $nk$ \kirrel{k} plus \kbound{k} edges and some \krel{k} edges.

In \cite{ps-mtcsp-09}, the triangles and their bisectors are generalized for \ktri{k}s as follows (see \fref{fig:2triang8pointsstars}):

\begin{figure}[b]
	\capstart
	\centerline{\includegraphics[scale=1]{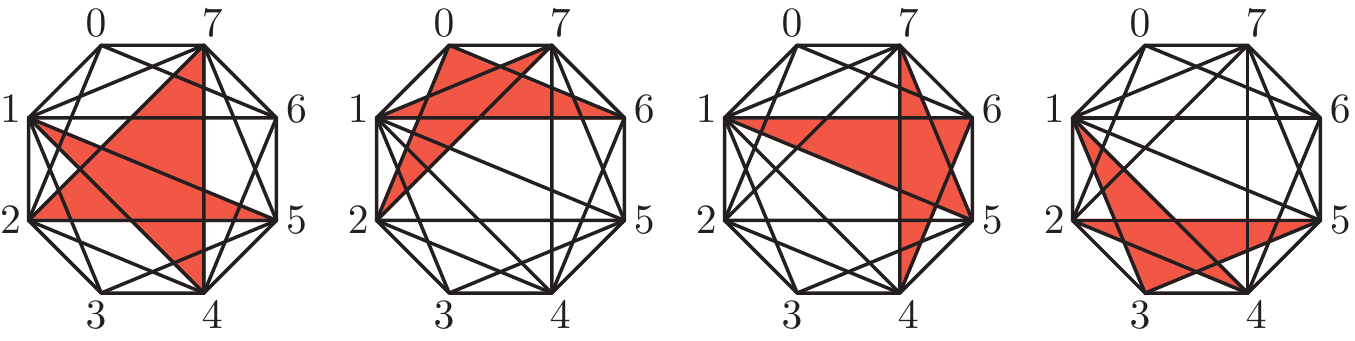}}
	\caption[\kstar{k}s in a \ktri{k}]{The four \kstar{2}s of the \ktri{2} of \fref{fig:2triang8points}.}
	\label{fig:2triang8pointsstars}
\end{figure}

\begin{definition}[\cite{ps-mtcsp-09}]
A \defn{\kstar{k}} is a star polygon of type $\left\{2k+1/k\right\}$, that is, a set of edges of the form $\set{s_js_{j+k}}{j\in \mathbb{Z}_{2k+1}}$, where $s_0,s_1,\ldots,s_{2k}$ are cyclically ordered around the unit circle. A \defn{(strict) bisector} of a \kstar{k} is a (strict) bisector of one of its angles $s_{j-k}s_js_{j+k}$. 
\end{definition}

\begin{figure}
	\capstart
	\centerline{\includegraphics[scale=1]{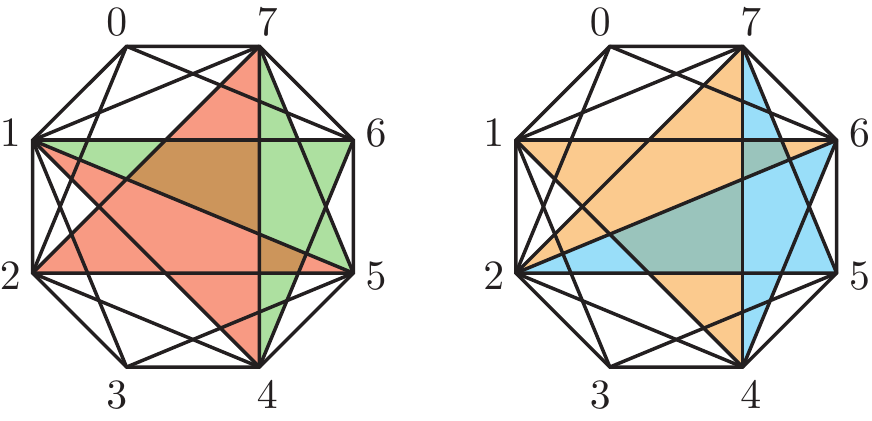}}
	\caption[A flip in a multitriangulation]{A flip in the \ktri{2} of \fref{fig:2triang8points}.}
	\label{fig:2triang8pointsflip}
\end{figure}

As for $k=1$, where triangles provide a powerful tool to study triangulations, \kstar{k}s are useful to understand \ktri{k}s. In the following theorem, we point out five properties of stars proved in \cite{ps-mtcsp-09}. 
\fref{fig:2triang8pointsstars} and  \fref{fig:2triang8pointsflip} illustrate these results on the \ktri{2} of \fref{fig:2triang8points}.

\begin{theorem}[\cite{ps-mtcsp-09}]\label{theo:stars}
Let $T$ be a \ktri{k} of the \gon{n}. Then
\begin{enumerate}[(i)]
\item $T$ contains exactly $n-2k$ \kstar{k}s and $k(n-2k-1)$ \krel{k} edges.
\item Each edge of $T$ belongs to zero, one, or two \kstar{k}s, depending on whether it is \kirrel{k}, \kbound{k}, or \krel{k}.
\item Every pair of \kstar{k}s of $T$ has a unique common strict bisector.
\item Flipping any \krel{k} edge $e$ of $T$ into the common strict bisector $f$ of the two \kstar{k}s containing $e$ produces a new \ktri{k} $T\triangle\{e,f\}$ of the \gon{n}. $T$ and $T\triangle\{e,f\}$ are the only two \ktri{k}s of the \gon{n} containing $T\ssm\{e\}$.
\item The flip graph~$G_{n,k}$ on \ktri{k}s of the \gon{n} is connected and regular of degree~$k(n-2k-1)$.
\end{enumerate}
\end{theorem}

Similarly to Observation~\ref{observation:pseudotriangulations}, we can interpret these properties of the stars of the multitriangulations in the dual space (see \fref{fig:2triang8pointsdual}):

\begin{observation}
Let~$T$ be a \ktri{k} of a convex \gon{n}. Then:
\begin{enumerate}[(i)]
\item the set~$S^*$ of all bisectors of a \kstar{k}~$S$ of~$T$ is a pseudoline of the M\"obius strip;
\item the dual pseudolines~$S_1^*, S_2^*$ of any two \kstar{k}s~$S_1,\Delta_2$ of~$T$ have a unique crossing point (the unique common strict bisector of~$S_1$ and~$S_2$) and possibly some contact points (when~$S_1$ and~$S_2$ share common edges);
\item the set~$T^* \eqdef \set{S^*}{S\;k\text{-star of } T}$ of dual pseudolines of \kstar{k}s of~$T$ is a pseudoline arrangement (with contact points);~and
\item $T^*$ is supported by the dual pseudoline arrangement~$C_n^*$ of $C_n$ minus its first $k$~levels (see \fref{fig:2triang8pointsdual}(b)).
\end{enumerate}
\begin{figure}
	\capstart
	\centerline{\includegraphics[scale=1]{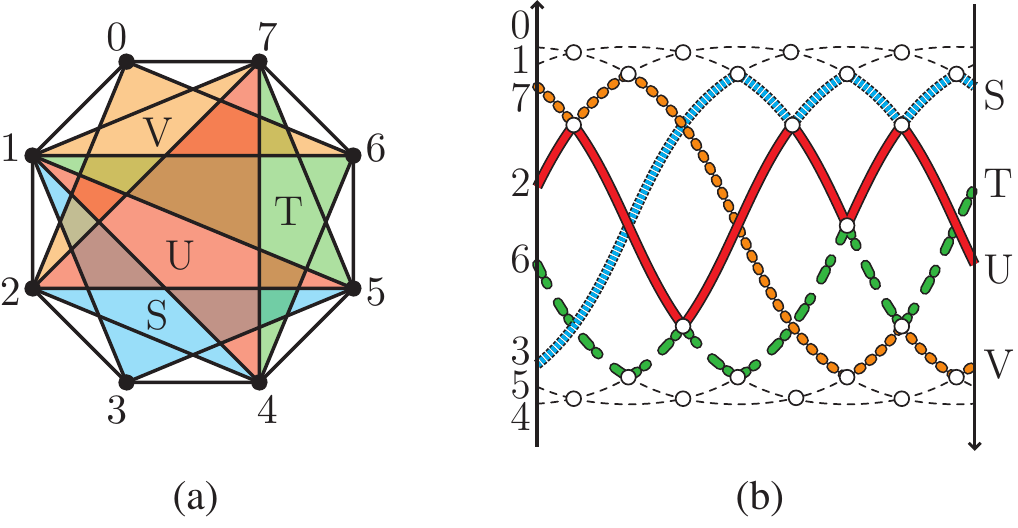}}
	\caption[A \ktri{2} of the octagon and its dual pseudoline arrangement]{A \ktri{2}~$T$ of the octagon~(a) and its dual pseudoline arrangement~$T^*$~(b). Each thick pseudoline of~$T^*$ corresponds to a \kstar{2} of~$T$; each contact point in~$T^*$ corresponds to an edge in~$T$; each crossing point in~$T^*$ corresponds to a common bisector in~$T$.}
	\label{fig:2triang8pointsdual}
\end{figure}
\end{observation}

Again, it turns out that this observation provides a characterization of multitriangulations of a convex polygon:

\begin{theorem}\label{theo:dualitymt}
Let~$C_n$ denote the set of vertices of a convex \gon{n}, and~$C_n^{*k}$ denote the support of its dual pseudoline arrangement minus its first $k$~levels. Then:
\begin{enumerate}[(i)]
\item The dual pseudoline arran\-gement~$T^* \eqdef \set{S^*}{S\;k\text{-star of }T}$ of a \ktri{k}~$T$ of the \gon{n} is supported by~$C_n^{*k}$.
\item The primal set of edges $$\qquad E \eqdef \set{[p,q]}{p,q\in C_n,\; p^*\wedge q^*\text{ contact point of } \Lambda}$$ of a pseudoline arrangement~$\Lambda$ supported by~$C_n^{*k}$ is a \ktri{k} of the \gon{n}.
\end{enumerate}
\end{theorem}

We provide two proofs of this theorem.

\begin{proof}[Proof~$1$ of Theorem~\ref{theo:dualitymt}(ii)]
The two notions of flips (the primal notion on \ktri{k}s of the \gon{n} and the dual notion on pseudoline arrangements supported by~$C_n^{*k}$) coincide. Thus, the flip graph~$G_{n,k}$ on \ktri{k}s of the \gon{n} is a subgraph of~$G(C_n^{*k})$. Since they are both connected and regular of degree~$k(n-2k-1)$, these two graphs coincide. In particular, any pseudoline arrangement supported by~$C_n^{*k}$ is the dual of a \ktri{k} of the \gon{n}.
\end{proof}

\begin{proof}[Proof~$2$ of Theorem~\ref{theo:dualitymt}(ii)]
We follow the method of our second proof of Theorem~\ref{theo:dualitypt}(ii). Since~$E$ has the right number of edges (namely $k(2n-2k-1)$), we only have to prove that it is \kcross{(k+1)}-free. We consider~$2k+2$ points~$p_0,\dots, p_k, q_0,\dots, q_k$ cyclically ordered around the unit circle. Since the definition of crossing (and thus, of \kcross{\ell}) is purely combinatorial, \ie depends only on the cyclic order of the points and not on their exact positions, we can move all the vertices of our \gon{n} on the unit circle while preserving their cyclic order. In particular, we can assume that the lines~$(p_iq_i)_{i\in\{0,\dots,k\}}$ all contain a common point~$t$. Its dual pseudoline~$t^*$ crosses~$C_n^*$ exactly~$n$ times and~$\Lambda$ at least~$|\Lambda|=n-2k$ times. Furthermore, for any point~${p_i^*\wedge q_i^*}$, replacing the crossing point by a contact point removes two crossings with~$t^*$. Thus, the pseudoline~$t^*$ provides a witness which proves that the edges~$[p_i,q_i]$,~$i\in\{0,\dots,k\}$, cannot be all in~$E$, and thus ensures that~$E$ is \kcross{(k+1)}-free.
\end{proof}


\section{Multipseudotriangulations}\label{sec:mpt}

Motivated by Theorems~\ref{theo:dualitypt} and~\ref{theo:dualitymt}, we define in terms of pseudoline arrangements a natural generalization of both pseudotriangulations and multitriangulations. We then study elementary properties of the corresponding set of edges in the primal space.


\subsection{Definition}\label{subsec:mpt:definition}

We consider the following generalizations of both pseudotriangulations and multitriangulations:

\begin{definition}
Let~$L$ be a pseudoline arrangement supported by~$\cS$. Define its \defn{\kkernel{k}}~$\cS^k$ to be its support minus its first~$k$ levels (which are the iterated external hulls of~$\cS$). Denote by~$V^k$ the set of contact points of~$L$ in~$\cS^k$. A \defn{\pt{k}} of~$L$ is a pseudoline arrangement whose support is~$\cS^k$ and whose set of contact points contains~$V^k$. 
\end{definition}

Pseudotriangulations of a point set~$P$ correspond via duality to \pt{1}s of the dual pseudoline arrangement~$P^*$. Similarly, \ktri{k}s of the \gon{n} correspond to \pt{k}s of the pseudoline arrangement~$C_n^*$ in convex position.
If~$L$ is a pseudoline arrangement with no contact point, then any pseudoline arrangement supported by~$\cS^k$ is a \pt{k} of~$L$. In general, the condition that the contact points of~$L$ in its \kkernel{k} should be contact points of any \pt{k} of~$L$ is a natural assumption for iterating \mpt{}s (see Section~\ref{sec:iterated}).

Let~$\Lambda$ be a \pt{k} of~$L$. We denote by~$V(\Lambda)$ the union of the set of contact points of~$\Lambda$ with the set of intersection points of the first~$k$ levels of~$L$. In other words,~$V(\Lambda)$ is the set of intersection points of~$L$ which are not crossing points of~$\Lambda$. As for pseudoline arrangements, the set~$V(\Lambda)$ completely determines~$\Lambda$.

Flips for \mpt{}s are defined as in Lemma~\ref{lem:flip}, with the restriction that the contact points in~$V^k$ cannot be flipped. In other words, the flip graph on \pt{k}s of~$L$ is exactly the graph~$G(\cS^k\,|\, V^k)$. Section~\ref{sec:enumeration} asserts that the graph of flips is regular and connected, and provides an enumeration algorithm for \mpt{}s~of~$L$.

\begin{figure}
	\capstart
	\centerline{\includegraphics[scale=1]{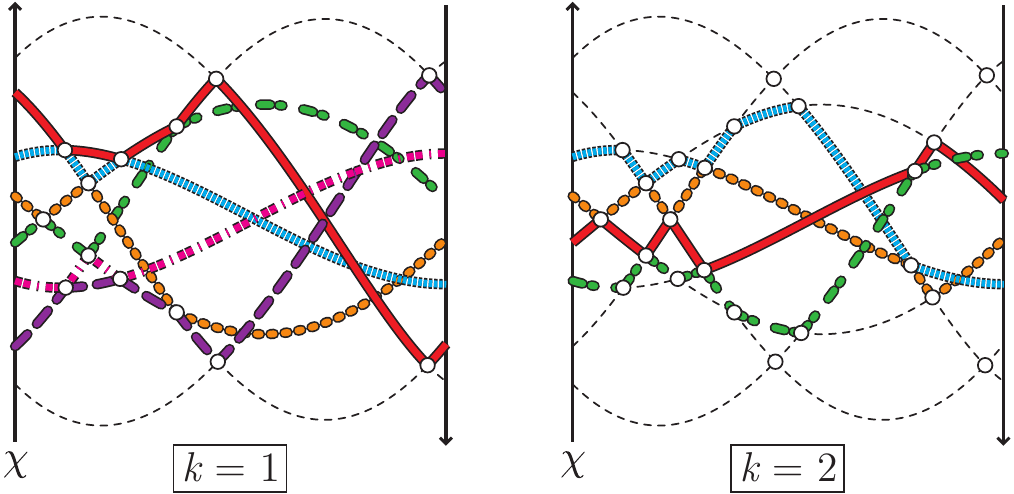}}
	\caption[Greedy \mpt{}s]{The \greedy{\chi} \pt{1} and the \greedy{\chi} \pt{2} of the pseudoline arrangement of \fref{fig:dual}(b).}
	\label{fig:greedy}
\end{figure}

Let~$\chi$ be a cut of (the support of)~$L$. It is also a cut of the \kkernel{k}~$\cS^k$ of~$L$. A particularly interesting example of \pt{k} of~$L$ is the source of the graph of \increasing{\chi} flips on \pt{k}s~of~$L$ (see \fref{fig:greedy} for an illustration):

\begin{definition}
The \defn{\greedy{\chi} \pt{k}} of~$L$, denoted $\Gamma_\chi^k(L)$, is the greedy pseudoline arrangement~$\Gamma_\chi(\cS^k\,|\, V^k)$.
\end{definition}

\subsection{Pointedness and crossings}\label{subsec:mpt:pointedcrossing}

Let~$P$ be a point set in general position. Let~$\Lambda$ be a \pt{k} of~$P^*$ and~$V(\Lambda)$ be the set of crossing points of~$P^*$ which are not crossing points of~$\Lambda$. We call \defn{primal of~$\Lambda$} the set $$E \eqdef \set{[p,q]}{p,q\in P,\; p^*\wedge q^*\in V(\Lambda)}$$ of edges of~$P$ primal to~$V(\Lambda)$ (see \fref{fig:2-pseudotriangulation}). Here, we discuss general properties of primals of \mpt{}s. We start with elementary properties that we already observed for the special cases of pseudotriangulations and multitriangulations in the proofs of Theorems~\ref{theo:dualitypt} and~\ref{theo:dualitymt}:

\begin{figure}
	\capstart
	\centerline{\includegraphics[scale=1]{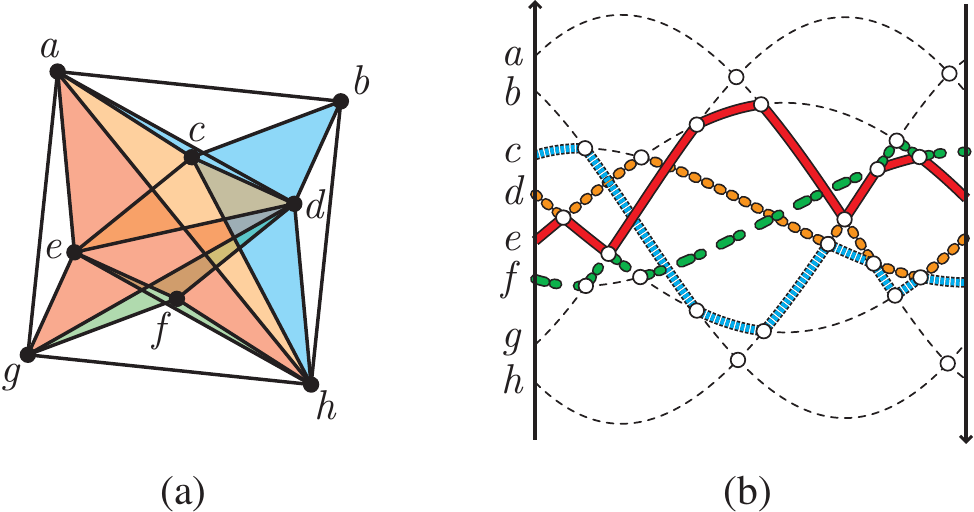}}
	\caption[The primal of a \pt{2}]{The primal set of edges~(a) of a \pt{2} (b)~of the dual pseudoline of the point set of~\fref{fig:dual}(a).}
	\label{fig:2-pseudotriangulation}
\end{figure}

\begin{lemma}
The primal~$E$ of~$\Lambda$ has~$k(2|P|-2k-1)$ edges.
\end{lemma}

\begin{proof}
The number of edges of~$E$ is the difference between the number of crossing points in the pseudoline arrangements~$P^*$ and~$\Lambda$:
$$|E|={\left|P^*\right| \choose 2}-{|\Lambda| \choose 2}={|P| \choose 2}-{|P|-2k \choose 2}=k(2|P|-2k-1).\qedhere$$
\end{proof}

We now discuss pointedness of~$E$. We call \defn{\kalter{k}} any set $\set{f_i}{i\in\Z_{2k+1}}$ of~$2k+1$ edges all incident to a common vertex and whose cyclic order around it is given by $$f_0\cl \bar f_{1+k} \cl f_1 \cl \bar f_{2+k} \cl \dots \cl f_{2k} \cl \bar f_k \cl f_0,$$ where~$\bar f_i$ denotes the opposite direction of the edge~$f_i$.

\begin{lemma}\label{lem:pointed}
The primal~$E$ of~$\Lambda$ cannot contain a \kalter{k}.
\end{lemma}

\begin{proof}
We simply mimic the proof of pointedness in Theorem~\ref{theo:dualitypt}. Let $p_0,\dots,p_{2k}$ and~$q$ be~${2k+2}$ points of~$P$ such that~$F \eqdef \set{[p_i,q]}{i\in\Z_{2k+1}}$ is a \kalter{k}. We prove that~$F$ cannot be a subset of~$E$ by constructing a witness pseudoline~$\ell$ that separates all the crossing points~${p_i^*\wedge q^*}$ corresponding to~$F$, while crossing~$q^*$ exactly~$2k+1$ times and the other pseudolines of~$P^*$ exactly as~$q^*$ does. (We skip the precise construction, since it is exactly the same as in the proof of Theorem~\ref{theo:dualitypt}.) Counting the crossings of~$\ell$ with~$P^*$ and~$\Lambda$, we obtain:
\begin{enumerate}[(i)]
\item $\ell$~crosses~$P^*$ exactly~$|P|+2k$ times;
\item $\ell$~crosses~$\Lambda$ at least~$|\Lambda|=|P|-2k$ times;
\item for each of the points~$p_i^*\wedge q^*$, replacing the crossing point by a contact point removes two crossings with~$\ell$.
\end{enumerate}
Thus the edges $[p_i,q]$ cannot all be contained in~$E$.
\end{proof}

\begin{figure}[h]
	\capstart
	\centerline{\includegraphics[scale=1]{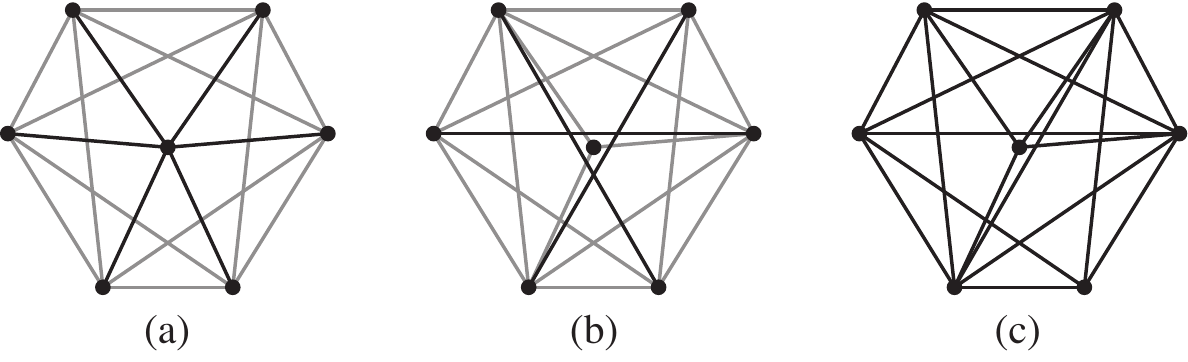}}
	\caption[A non-\kpointed{2} \pt{2}, a \pt{2} containing a \kcross{3} and a \kcross{3}-free \kalter{2}-free set not contained in a \pt{2}]{(a) A non-\kpointed{2} (but \kalter{2}-free) \pt{2}. (b)~A \pt{2} containing a \kcross{3}. (c) A \kcross{3}-free \kalter{2}-free set not contained in a \pt{2}.}
	\label{fig:3crossing}
\end{figure}

\begin{remark}
Observe that a set of edges is pointed if and only if it is \kalter{1}-free.
In contrast, we want to observe the difference between \kalter{k}-freeness and the following natural notion of \kpointed{k}ness: we say that a set~$F$ of edges with vertices in~$P$ is \defn{\kpointed{k}} if for all~$p$ in~$P$, there exists a line which passes through~$p$ and defines a half-plane that contains at most~$k-1$ segments of~$F$ adjacent to~$p$.
Observe that a \kpointed{k} set is automatically \kalter{k}-free but that the converse statement does not hold (see \fref{fig:3crossing}(a)).
\end{remark}

Finally, contrarily to pseudotriangulations ($k=1$) and multitriangulations (convex position), the condition of avoiding \kcross{(k+1)}s does not hold for \pt{k}s in general:

\begin{remark}
There exist \pt{k}s with \kcross{(k+1)}s (see \fref{fig:3crossing}(b)) as well as \kcross{(k+1)}-free \kalter{k}-free sets of edges that are not subsets of \pt{k}s (see \fref{fig:3crossing}(c)).
\end{remark}

\subsection{Stars in \mpt{}s}\label{subsec:mpt:stars}

To complete our understanding of the primal of \mpt{}s, we need to generalize pseudotriangles of pseudotriangulations and \kstar{k}s of \ktri{k}s: both pseudotriangles and \kstar{k}s correspond to pseudolines of the covering pseudoline arrangement.

We keep the notations of the previous section:~$P$~is a point set in general position, $\Lambda$~is a \pt{k} of~$P^*$ and~$E$~is the primal set of edges of~$\Lambda$. 

\begin{definition}
We call \defn{star} of a pseudoline~$\lambda \in \Lambda$ the set of edges $${S(\lambda) \eqdef \set{[p,q]}{p,q\in P,\; p^*\wedge q^*\text{ contact point of }\lambda}}$$ primal to the contact points of~$\lambda$.
\end{definition}

\begin{lemma}
For any~$\lambda\in \Lambda$, the star~$S(\lambda)$ is non-empty.
\end{lemma}

\begin{proof}
We have to prove that any pseudoline~$\lambda$ of~$\Lambda$ supports at least one contact point. If it is not the case, then~$\lambda$ is also a pseudoline of~$P^*$, and all the~$|P|-1$ crossing points of~$\lambda$ with~$P^*\ssm\{\lambda\}$ should be crossing points of~$\lambda$ with the arrangement~$\Lambda\ssm\{\lambda\}$. This is impossible since $|\Lambda\ssm\{\lambda\}|=|P|-2k-1$.
\end{proof}

Similarly to the case of \ktri{k}s of the \gon{n}, we say that an edge~$[p,q]$ of~$E$ is a \defn{\krel{k}} (resp.~\defn{\kbound{k}}, resp.~\defn{\kirrel{k}}) edge if there remain strictly more than (resp.~exactly, resp.~strictly less than)~$k-1$~points of~$P$ on each side (resp.~one side) of the line~$(pq)$. In other words,~$p^*\wedge q^*$ is located in the \kkernel{k} (resp.~in the intersection of the $k$th level and the \kkernel{k}, resp.~in the first $k$ levels) of the pseudoline arrangement~$P^*$. Thus, the edge~$[p,q]$ is contained in~$2$ (resp.~$1$, resp.~$0$) stars of~$\Lambda$.

\begin{figure}
	\capstart
	\centerline{\includegraphics[scale=1]{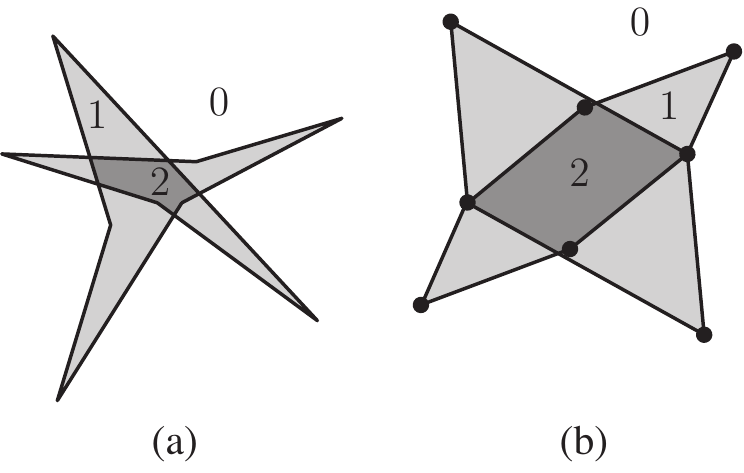}}
	\caption[Winding number of a star; depth in a point set]{(a) Winding number of a star. (b)~\kdepth{2} in the point set of \fref{fig:dual}(a).}
	\label{fig:depth_winding_number}
\end{figure}

The edges of a star~$S(\lambda)$ are cyclically ordered by the order of their dual contact points on~$\lambda$, and thus~$S(\lambda)$ forms a (not-necessarily simple) polygonal cycle. For any point~$q$ in the plane, let~$\sigma_\lambda(q)$ denote the \defn{winding number} of~$S(\lambda)$ around~$q$, that is, the number of rounds made by~$S(\lambda)$ around the point~$q$ (see \fref{fig:depth_winding_number}(a)). For example, the winding number of a point in the external face is~$0$.

We call \defn{\kdepth{k}} of a point~$q$ the number~$\delta^k(q)$ of \kbound{k} edges of~$P$ crossed by any (generic continuous) path from~$q$ to the external face, counted positively when passing from the ``big'' side (the one containing at least~$k$ vertices of~$P$) to the ``small side'' (the one containing~$k-1$ vertices of~$P$), and negatively otherwise (see \fref{fig:depth_winding_number}(b)). That this number is independent from the path can be seen by mutation. For example,~$\delta^1(q)$ is~$1$ if~$q$ is in the convex hull of~$P$ and~$0$ otherwise.

\begin{proposition}\label{prop:decomposition}
Any point~$q$ of the plane is covered~$\delta^k(q)$ times by the stars~$S(\lambda)$, $\lambda\in \Lambda$, of the \pt{k}~$\Lambda$ of~$P^*$:
$$ \delta^k(q)=\sum_{\lambda\in\Lambda} \sigma_\ell(q).$$
\end{proposition}

The proposition is intuitively clear: let us walk on a continuous path from the external face to the point~$q$. Initially, the winding numbers of all stars of~$\Lambda$ are zero (we start outside all stars of~$\Lambda$). Then, each time we cross an edge~$e$:
\begin{enumerate}[(i)]
\item If~$e$ is \kirrel{k}, it is not contained in any star of~$\Lambda$, and we do not change the winding numbers of the stars of~$\Lambda$.
\item If~$e$ is a \kbound{k} edge, and if we cross it positively, we increase the winding number of the star~$S$ of~$\Lambda$ containing~$e$; if we cross~$e$ negatively, we decrease the winding~number~of~$S$.
\item If~$e$ is \krel{k}, then we decrease the winding number of one star of~$\Lambda$ containing~$e$ and increase the winding number of the other star of~$\Lambda$ containing~$e$.
\end{enumerate}
Let us give a formal proof in the dual:

\begin{proof}[Proof of Proposition~\ref{prop:decomposition}]
Both~$\sigma_\lambda(q)$ and~$\delta^k(q)$ can be read on the pseudoline~$q^*$:
\begin{enumerate}[(i)]
\item If~$\tau_\lambda(q)$ denotes the number of intersection points between~$q^*$ and~$\lambda$ (that is, the number of tangents to~$S(\lambda)$ passing through~$q$), then~$\sigma_\lambda(q)=(\tau_\lambda(q)-1)/2$.
\item If~$\gamma^k(q)$ denotes the number of intersection points between~$q^*$ and the first~$k$ levels of~$P^*$, then~$\delta^k(q)=k-\gamma^k(q)/2$.
\end{enumerate}
The pseudoline~$q^*$ has exactly~$|P|$ crossings with~$P^*$ (since~$P^*\cup\{q^*\}$ is an arrangement), which are crossings either with the pseudolines of~$\Lambda$ or with the first~$k$ levels of~$P^*$. Hence,
$$|P|=\gamma^k(q)+\sum_{\lambda\in \Lambda} \tau_\lambda(q)=2k-2\delta^k(q)+|\Lambda|+2\sum_{\lambda\in \Lambda}\sigma_\lambda(q),$$
and we get the aforementioned result since~$|\Lambda|=|P|-2k$.
\end{proof}

\begin{remark}
As a consequence of Proposition~\ref{prop:decomposition}, we obtain that the \kdepth{k} of any point~$q$ in any point set~$P$ is always non-negative (as a sum of non-negative numbers). It is interesting to notice that Welzl proved in~\cite{w-elf-01} that this non-negativity property is actually equivalent to the Lower Bound Theorem for $d$-dimensional polytopes with $d+3$ vertices.
\end{remark}

\begin{figure}
	\capstart
	\centerline{\includegraphics[scale=1]{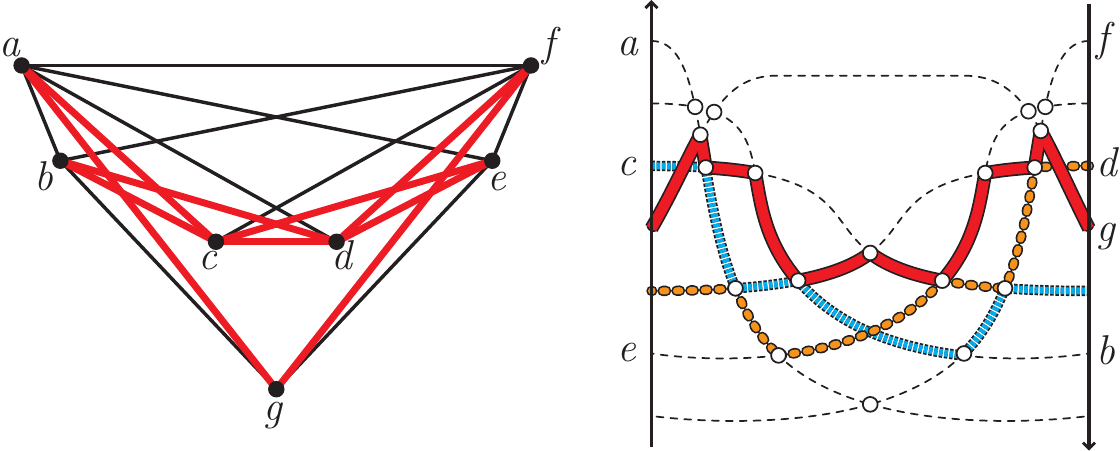}}
	\caption[A star of a \pt{2} with~$2|P|-5$ corners]{A star of a \pt{2} of $7$ points with $9$ corners. This example can be generalized to a \pt{2} of $n$ points with $2n-5$ corners.}
	\label{fig:2nstar}
\end{figure}

A \defn{corner} of the star~$S(\lambda)$ is an internal convex angle of it (see \fref{fig:star}(a)). In the following proposition, we are interested in the number of corners of~$S(\lambda)$. For points in convex position (Section~\ref{subsec:duality:mt}), the number of corners of a star is always~$2k+1$. This is not true anymore in general position as illustrated in \fref{fig:2nstar}. The following proposition gives our best bounds on the number of corners of the star of a \mpt{}.

\begin{proposition}\label{prop:corners}
The number of corners of a star~$S(\lambda)$ of a \pt{k} of~$P^*$ is odd and between~$2k+1$~and $2(k-1)|P|+2k+1$.
\end{proposition}

\begin{proof}
We read convexity of internal angles of~$S(\lambda)$ on the preimage~$\bar\lambda$ of the pseudoline~$\lambda$ under the projection~$\pi$. Let~$pqr$ be an internal angle, let~$v=p^*\wedge q^*$~and~$w=q^*\wedge r^*$ denote the contact points corresponding to the two edges~$[p,q]$~and~$[q,r]$ of this angle, and let~$\bar v$ and~$\bar w$ denote two consecutive preimages of~$v$~and~$w$ on~$\bar\lambda$ (meaning that~$\bar w$ is located between~$\bar v$~and~$\tau(\bar v)$). The angle~$pqr$ is a corner if and only if~$\bar v$~and~$\bar w$ lie on opposite sides of~$\bar\lambda$, meaning that the other curves touching~$\bar\lambda$ at~$\bar v$~and~$\bar w$ lie on opposite sides, one above and one below~$\bar\lambda$ (see \fref{fig:star}(b)).

\begin{figure}
	\capstart
	\centerline{\includegraphics[scale=1]{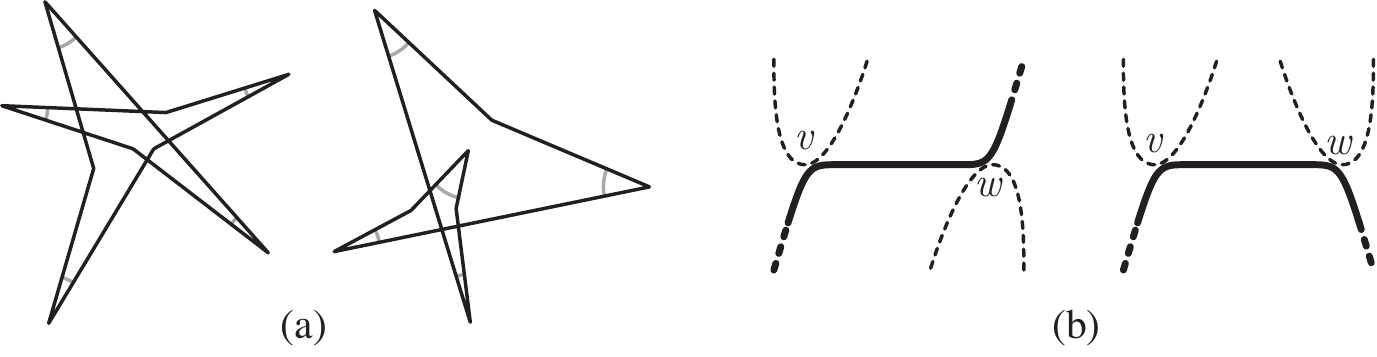}}
	\caption[Two stars with~$5$ corners; The two possible configurations of two consecutive contact points on~$\lambda$]{(a) Two stars with~$5$ corners. (b)~The two possible configurations of two consecutive contact points on~$\lambda$: convex (left) and concave (right).}
	\label{fig:star}
\end{figure}

\medskip

In particular, the number~$c(\lambda)=c$ of corners of~$S(\lambda)$ is the number of opposite consecutive contact points on~$\bar\lambda$ between two versions~$\bar v$~and~$\tau(\bar v)$ of a contact point~$v$ of~$\lambda$. To see that~$c$ is odd, imagine that we are discovering the contact points of~$\lambda$ one by one. The first contact point~$v$ that we see corresponds to two opposite contact points~$\bar v$~and~$\tau(\bar v)$ on~$\bar\lambda$. Then, at each stage, we insert a new contact point~$\bar w$ between two old contact points that can be:
\begin{enumerate}[(i)]
\item either on opposite sides and then we are not changing~$c$;
\item or on the same side and we are adding to~$c$ either~$0$ (if~$\bar w$ is also on the same side) or~$2$ (if~$\bar w$ is on the opposite side).
\end{enumerate}
Thus,~$c$~remains odd in any case.

\medskip

To prove the lower bound, we use our witness method. We perturb~$\lambda$ a little bit to obtain a pseudoline~$\mu$ that passes on the opposite side of each contact point (this is possible since~$c$ is odd). This pseudoline~$\mu$ crosses~$\lambda$ between each pair of opposite contact points and crosses the other pseudolines of~$\Lambda$ exactly as~$\lambda$ does. Thus,~$\mu$ crosses~$\Lambda$ exactly $|\Lambda|-1+c$ times. But since~$\mu$ is a pseudoline, it has to cross all the pseudolines of~$P^*$ at least once. Thus, $|P| \le |\Lambda|-1+c=|P|-2k-1+c$ and~$c\ge 2k+1$.

From this lower bound, we derive automatically the upper bound. Indeed, we know that the number of corners around one point~$p$ is at most~$\deg(p)-1$. Consequently,
$$\sum_{p\in P} (\deg(p)-1) \ge \sum_{\nu\in \Lambda} c(\nu) = c(\lambda)+\sum_{\substack{\nu\in \Lambda\\ \nu\ne\lambda}} c(\nu).$$
The left sum equals~$2k(2|P|-2k-1)-|P|$ while, according to the previous lower bound, the right one is at least~$c+(|P|-2k-1)(2k+1)$. Thus we get~$c\le 2(k-1)|P|+2k+1$.
\end{proof}


\newpage
\section{Iterated \mpt{}s}\label{sec:iterated}

By definition, a \pt{k} of an \pt{m} of a pseudoline arrangement~$L$ is a \pt{(k+m)} of~$L$. In this section, we study these iterated sequences of \mpt{}s. In particular, we compare \mpt{}s with iterated sequences of \pt{1}s.


\subsection{Definition and examples}\label{subsec:iterated:definition}

Let~$L$ be a pseudoline arrangement. An \defn{iterated \mpt{}} of~$L$ is a sequence $\Lambda_1,\dots,\Lambda_r$ of pseudoline arrangements such that~$\Lambda_i$ is a \mpt of~$\Lambda_{i-1}$ for all~$i$ (by convention,~$\Lambda_0=L$). We call \defn{signature} of~$\Lambda_1,\dots,\Lambda_r$ the sequence~$k_1<\cdots<k_r$ of integers such that~$\Lambda_i$ is a \pt{k_i} of~$L$ for all~$i$. Observe that the assumption that contact points of a pseudoline arrangement~$L$ should be contact points of any \mpt of~$L$ is natural in this setting: iterated \mpt{}s correspond to decreasing sequences of sets of crossing points.

A \defn{decomposition} of a \mpt~$\Lambda$ of a pseudoline arrangement~$L$ is an iterated \mpt~$\Lambda_1,\dots,\Lambda_r$ of~$L$ such that~$\Lambda_r=\Lambda$ and~$r>1$. We say that~$\Lambda$ is \defn{decomposable} if such a decomposition exists, and \defn{irreducible} otherwise. The decomposition is \defn{complete} if its signature is~$1,2,\dots,r$.

It is tempting to believe that all \mpt{}s are completely decomposable. This would allow to focus only on pseudotriangulations. However, we start by showing that not even all multitriangulations are decomposable. The following example is due to Francisco Santos.

\begin{example}[An irreducible \ktri{2} of the \gon{15}]\label{exm:irred}
We consider the geometric graph~$T$ of \fref{fig:15gonctrexm}. The edges are:
\begin{enumerate}[(i)]
\item all the \kirrel{2} and \kbound{2} edges of the \gon{15}, and
\item the five zigzags~${Z_a=\{[3a,3a+6], [3a+6,3a+1], [3a+1,3a+5],}$ $[3a+5,3a+2]\}$, for~$a\in\{0,1,2,3,4\}$.
\end{enumerate}
Thus,~$T$ has~$50$ edges and is \kcross{3}-free (since the only \krel{2} edges of~$T$ that cross a zigzag~$Z_a$ are edges of~$Z_{a-1}$ and~$Z_{a+1}$). Consequently,~$T$ is a \ktri{2} of the \gon{15}.

\begin{figure}
	\capstart
	\centerline{\includegraphics[scale=1]{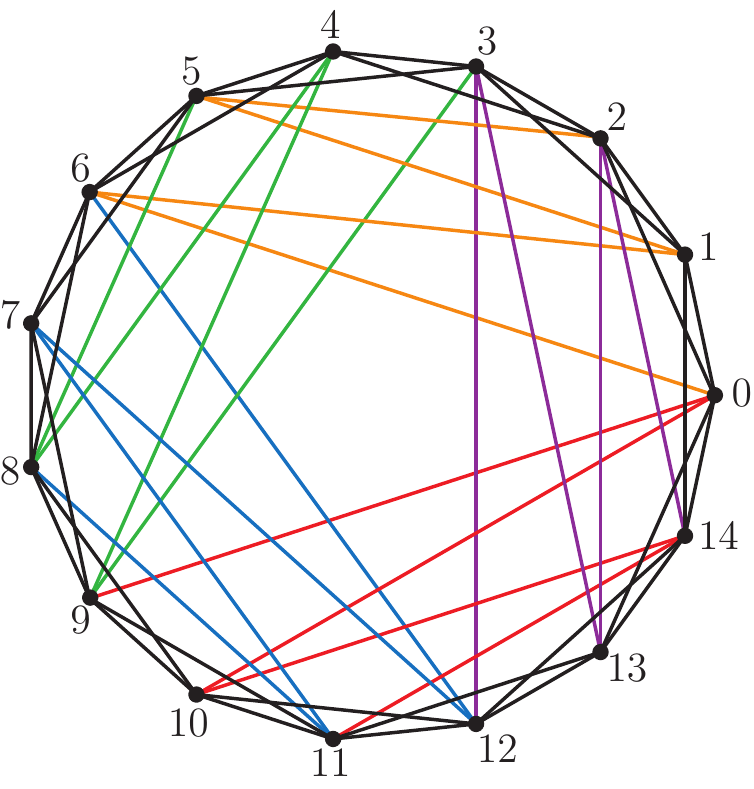}}
	\caption[An irreducible \ktri{2} of the \gon{15}]{An irreducible \ktri{2} of the \gon{15}: it contains no triangulation.}
	\label{fig:15gonctrexm}
\end{figure}

Let us now prove that~$T^*$ is irreducible, that is, that~$T$ contains no triangulation. Observe first that the edge~$[0,6]$ cannot be an edge of a triangulation contained in~$T$ since none of the triangles~$06i$,~${i\in\{7,\dots,14\}}$, is contained in~$T$. Thus, we are looking for a triangulation contained in~$T\ssm\{[0,6]\}$. Repeating the argument successively for the edges~$[1,6]$, $[1,5]$~and~$[2,5]$, we prove that the zigzag~$Z_0$ is disjoint from any triangulation contained in~$T$. By symmetry, this proves the irreducibility of~$T^*$.
\end{example}


\subsection{Iterated greedy pseudotriangulations}\label{subsec:iterated:greedy}

Greedy \mpt{}s provide interesting examples of iteration of pseudotriangulations. Let~$L$ be a pseudoline arrangement, and~$\chi$ be a cut of~$L$. 

\begin{theorem}
For any positive integers~$a$ and~$b$, $\Gamma_\chi^{a+b}(L)=\Gamma_\chi^b(\Gamma_\chi^a(L))$. Consequently, for any integer~$k$, $\Gamma_\chi^k(L)=\Gamma_\chi^1\circ\Gamma_\chi^1\circ\cdots\circ\Gamma_\chi^1(L)$, where~$\Gamma_\chi^1(.)$ is iterated~$k$ times.
\end{theorem}

\begin{proof}
Since~$\chi$ is a cut of~$L$, it is also a cut of~$\Gamma^a_\chi(L)$ and thus~$\Gamma^b_\chi(\Gamma^a_\chi(L))$ is well defined. Observe also that we can assume that~$L$ has no contact point (otherwise, we can open them). Let~$n \eqdef |L|$ and~$m \eqdef {n \choose 2}$.

Let~$\chi=\chi_0,\dots,\chi_m=\chi$ be a backward sweep of~$L$. For all~$i$, let~$v_i$ denote the vertex of~$L$ swept when passing from~$\chi_i$ to $\chi_{i+1}$, and~$i^\square$ denote the integer such that the pseudolines that cross at~$v_i$ are the $i^\square$th and $(i^\square+1)$th pseudolines of~$L$ on~$\chi_i$.

Let~$\sigma_0,\dots,\sigma_m$ denote the sequence of permutations corresponding to~$\Gamma^a_\chi(L)$ on the sweep $\chi_0,\dots,\chi_m$. In other words,~$\sigma_0$ is the permutation
$$[1,\dots,a,n-a,n-a-1,\dots,a+2,a+1,n-a+1,\dots,n],$$
whose first~$a$ and last~$a$ entries are preserved, while its~$n-2a$ intermediate entries are inverted. Then, for all~$i$, the permutation~$\sigma_{i+1}$ is obtained from~$\sigma_i$ by sorting its $i^\square$th and $(i^\square+1)$th entries.

Similarly, let~$\rho_0,\dots,\rho_m$ and~$\omega_0,\dots,\omega_m$ denote the sequences of permutations corresponding to~$\Gamma_\chi^{a+b}(L)$ and~$\Gamma^b_\chi(\Gamma^a_\chi(L))$ respectively: both $\rho_0$ and~$\omega_0$ equal the permutation whose first and last~$a+b$ entries are preserved and whose~$n-2a-2b$ intermediate entries are inverted,~and:
\begin{itemize}
\item $\rho_{i+1}$ is obtained from~$\rho_i$ by sorting its $i^\square$th and $(i^\square+1)$th entries;
\item if~$v_i\notin\Gamma^a_\chi(L)$, then~$\omega_{i+1}$ is obtained from~$\omega_i$ by sorting its $i^\square$th and $(i^\square+1)$th entries; otherwise,~$\omega_{i+1}=\omega_i$.
\end{itemize}

We claim that for all~$i$,
\begin{enumerate}[(A)]
\item all the inversions of~$\rho_i$ are also inversions of~$\sigma_i$: $\rho_i(p)>\rho_i(q)$ implies~$\sigma_i(p)>\sigma_i(q)$ for all~$1\le p<q\le n$; and
\item $\rho_i=\omega_i$.
\end{enumerate}

We prove this claim by induction on~$i$. It is clear for~$i=0$. Assume it holds for~$i$ and let us prove it for~$i+1$. We have two possible situations:

\begin{enumerate}
\item \textbf{First case:}~$\sigma_i(i^\square)<\sigma_i(i^\square+1)$. Then,~$\sigma_{i+1}=\sigma_i$ and ${v_i\in\Gamma^a_\chi(L)}$. Thus,~${\omega_{i+1}=\omega_i}$. Furthermore, using Property~(A) at rank~$i$, we know that~$\rho_i(i^\square)<\rho_i(i^\square+1)$, and thus ${\rho_{i+1}=\rho_i}$. To summarize, $\sigma_{i+1}=\sigma_i$, $\omega_{i+1}=\omega_i$, and $\rho_{i+1}=\rho_i$, which trivially implies that Properties (A) and (B) remain true.

\item \textbf{Second case:}~$\sigma_i(i^\square)>\sigma_i(i^\square+1)$. Then,~$\sigma_{i+1}$ is obtained from~$\sigma_i$ by exchanging the $i^\square$th and $(i^\square+1)$th entries, and $v_i\notin\Gamma^a_\chi(L)$. Consequently, $\rho_{i+1}$~and~$\omega_{i+1}$ are both obtained from~$\rho_i$~and~$\omega_i$ respectively by sorting their $i^\square$th and $(i^\square+1)$th entries. Thus, Property~(B) obviously remains true. As far as Property~(A) is concerned, the result is obvious if~$p$~and~$q$ are different from~$i^\square$~and~$i^\square+1$. By symmetry, it suffices to prove that for any~$p<i^\square$, $\rho_{i+1}(p)>\rho_{i+1}(i^\square)$~implies ${\sigma_{i+1}(p)>\sigma_{i+1}(i^\square)}$. We have to consider two subcases:
\begin{enumerate}
\item \textbf{First subcase:}~$\rho_i(i^\square)<\rho_i(i^\square+1)$. Then~$\rho_{i+1}=\rho_i$. Thus, if~$p<i^\square$ is such that~$\rho_{i+1}(p)>\rho_{i+1}(i^\square)$, then we have~$\rho_i(p)>\rho_i(i^\square)$. Consequently, we obtain that ${\sigma_{i+1}(p)=\sigma_i(p)>\sigma_i(i^\square)>\sigma_i(i^\square+1)=\sigma_{i+1}(i^\square)}$.
\item \textbf{Second subcase:}~$\rho_i(i^\square)>\rho_i(i^\square+1)$. Then~$\rho_{i+1}$ is obtained from~$\rho_i$ by exchanging its $i^\square$th and $(i^\square+1)$th entries. If~$p<i^\square$ is such that~$\rho_{i+1}(p)>\rho_{i+1}(i^\square)$, then we have ${\rho_i(p)>\rho_i(i^\square+1)}$. Consequently, we obtain that $\sigma_{i+1}(p)=\sigma_i(p)>\sigma_i(i^\square+1)=\sigma_{i+1}(i^\square)$.
\end{enumerate}
\end{enumerate}

Obviously, Property~(B) of our claim proves the theorem.
\end{proof}


\subsection{Flips in iterated \mpt{}s}\label{subsec:iterated:flips}

Let~$\Lambda_1,\dots,\Lambda_r$  be an iterated \mpt of a pseudoline arrangement~$L$, with signature~$k_1<\dots<k_r$.
Let~$v$ be a contact point of~$\Lambda_r$ (which is not a contact point of~$L$), and let~$i$ denote the first integer for which~$v$ is a contact point of~$\Lambda_i$ (thus,~$v$ is a contact point of~$\Lambda_j$ if and only if $i\le j\le r$). For all~$i\le j\le r$, let~$\Lambda'_j$ denote the pseudoline arrangement obtained from~$\Lambda_j$ by flipping~$v$, and let~$w_j$ denote the new contact point of~$\Lambda'_j$. Let~$j$ denote the biggest integer such that~$w_j=w_i$.
There are three possibilities:
\begin{enumerate}[(i)]
\item If~$j=r$, then~$\Lambda_1,\dots,\Lambda_{i-1},\Lambda'_i,\dots,\Lambda'_r$ is an iterated \mpt of~$L$. We say that it is obtained from~$\Lambda_1,\dots,\Lambda_r$ by a \defn{complete flip} of~$v$ .
\item If~$j<r$, and~$w_i=w_j$ is a contact point of~$\Lambda_{j+1}$, then
$$\Lambda_1,\dots,\Lambda_{i-1},\Lambda'_i,\dots,\Lambda'_j,\Lambda_{j+1},\dots,\Lambda_r$$ is an iterated \mpt of~$L$. We say that it is obtained from~$\Lambda_1,\dots,\Lambda_r$ by a \defn{partial flip} of~$v$.
\item If~$j<r$, and~$w_i=w_j$ is a crossing point in~$\Lambda_{j+1}$, then we cannot flip~$v$ in~$\Lambda_i$ maintaining an iterated \mpt of~$L$.
\end{enumerate}

\begin{figure}[b]
	\capstart
	\centerline{\includegraphics[scale=1]{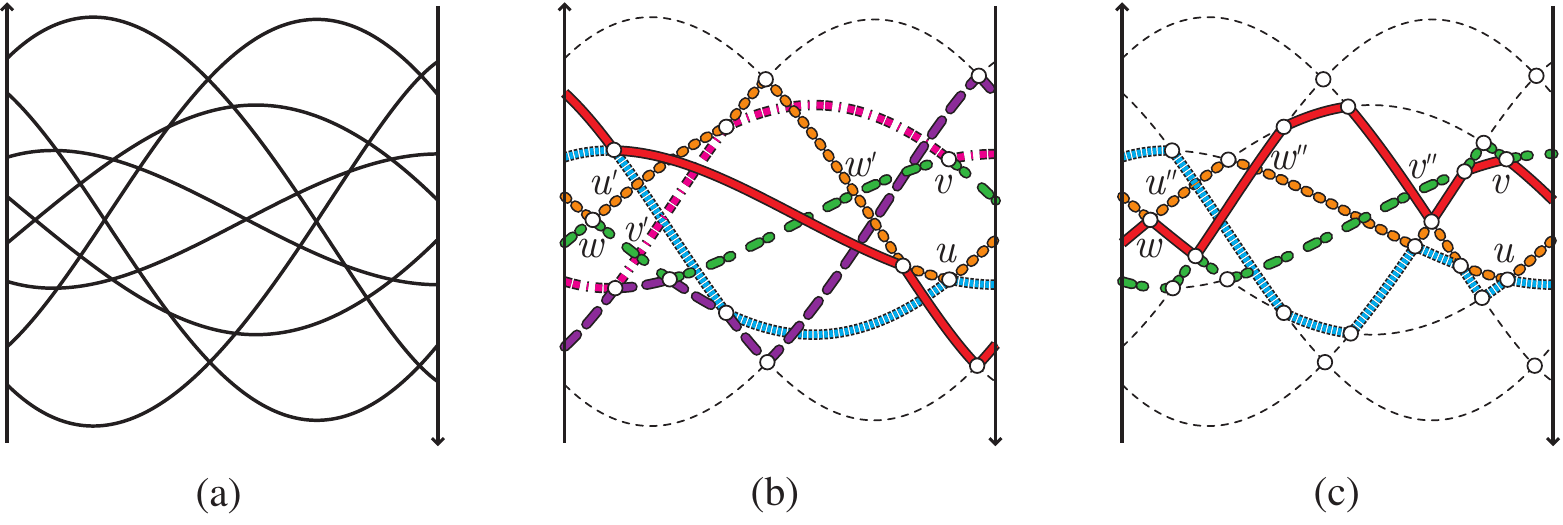}}
	\caption[The three possible situations for flipping a contact point in an iterated \mpt{}]{An iterated multipseudotriangulation: a pseudoline arrangement~$L$ (a), a \pt{1}~$\Lambda_1$ of~$L$ (b), and a \pt{1}~$\Lambda_2$ of~$\Lambda_1$~(c). The contact points~$u,v,w$ illustrate the three possible situations for flipping a contact point.}
	\label{fig:iteratedflip}
\end{figure}

To illustrate these three possible cases, we have labeled on \fref{fig:iteratedflip} some intersection points of an iterated pseudotriangulation. We have chosen three contact points~$u,v,w$ in~(b). For~$z\in\{u,v,w\}$, we label $z'$ (resp.~$z''$) the crossing point corresponding to~$z$ in~(b) (resp.~in~(c)). Observe that:
\begin{enumerate}[(i)]
\item points~$u'$~and~$u''$ coincide. Thus we can flip simultaneously point~$u$ in~(b) and~(c) (complete flip);
\item points~$v'$ is different from~$v''$ but is a contact point in~(c). Thus, we can just flip~$v$ in~(b), without changing~(c) and we preserve an iterated pseudotriangulation (partial flip);
\item point~$w'$ is a crossing point in~(c), different from~$w''$. Thus, we cannot flip~$w$ in~(b) maintaining an iterated pseudotriangulation.
\end{enumerate}

Let~$G^{k_1,\dots,k_r}(L)$ be the graph whose vertices are the iterated \mpt{}s of~$L$ with signature~$k_1<\dots<k_r$, and whose edges are the pairs of iterated \mpt{}s linked by a (complete or partial) flip.

\begin{theorem}\label{theo:iteratedflip}
The graph of flips~$G^{k_1,\dots,k_r}(L)$ is connected.
\end{theorem}

To prove this proposition, we need the following lemma:

\begin{lemma}\label{lem:exist}
Any intersection point~$v$ in the \kkernel{k} of a pseudoline arrangement is a contact point in a \pt{k} of it.
\end{lemma}

\begin{proof}
The result holds when~$k=1$. We obtain the general case by iteration.
\end{proof}

\begin{proof}[Proof of Theorem~\ref{theo:iteratedflip}]
We prove the result by induction on~$r$ ($L$~is fixed). When~$r=1$, we already know that the flip graph is connected. Now, let~$A_-$~and~$A_+$ be two iterated \mpt{}s of~$L$ with signature~$k_1<\dots<k_r$, that we want to join by flips. Let~$B_-$~and~$B_+$ be iterated \mpt{}s of~$L$ with signature~$k_1<\dots<k_{r-1}$, and~$\Lambda_-$~and~$\Lambda_+$ be \pt{k_r}s of~$L$ such that~$A_-=B_-,\Lambda_-$ and~$A_+=B_+,\Lambda_+$.

By induction,~$G^{k_1,\dots,k_{r-1}}(L)$ is connected: let $$B_-=B_1,B_2,\dots,B_{p-1},B_p=B_+$$ be a path from~$B_-$ to~$B_+$ in~$G^{k_1,\dots,k_{r-1}}(L)$. For all~$j$, let~$v_j$ be such that~$B_{j+1}$ is obtained from~$B_j$ by flipping~$v_j$  and let~$w_j$ be such that~$B_j$ is obtained from~$B_{j+1}$ by flipping~$w_j$. Let~$\Lambda_j$ be a \pt{k_r} of~$L$ containing the contact points of~$B_j$ plus~$w_j$ (it exists by Lemma~\ref{lem:exist}), and let~$C_j=B_j,\Lambda_j$. Let~$D_j$ be the iterated \mpt of~$L$ obtained from the iterated pseudotriangulation~$C_j$ by a partial flip of~$v_j$. Finally, since~$G^{k_r}(B_j)$ is connected, there is a path of complete flips from~$D_{j-1}$ to~$C_j$.

Merging all these paths, we obtain a global path from~$A_-$ to~$A_+$: we transform~$A_-$ into~$C_1$ via a path of complete flips; then~$C_1$ into~$D_1$ by the partial flip of~$v_1$; then~$D_1$ into~$C_2$ via a path of complete flips; then~$C_2$ into~$D_2$ by the partial flip of~$v_2$; and so on.
\end{proof}


\section{Further topics}\label{sec:furthertopics}

We discuss here the extensions in the context of \mpt{}s of two known results on pseudotriangulations:
\begin{enumerate}
\item The first one concerns the connection between the greedy pseudotriangulation of a point set and its horizon trees.
\item The second one extends to arrangements of double pseudolines the definition and properties of \mpt{}s.
\end{enumerate}


\subsection{Greedy \mpt{}s and horizon graphs}\label{subsec:furthertopics:horizon}

We have seen in previous sections that the greedy \pt{k} of a pseudoline arrangement~$L$ can be seen as:
\begin{enumerate}
\item the unique source of the graph of increasing flips;
\item a greedy choice of crossing points given by a sorting network;
\item a greedy choice of contact points;
\item an iteration of greedy \pt{1}s.
\end{enumerate}
In this section, we provide a ``pattern avoiding'' characterization of the crossing points of the greedy \pt{k} of~$L$.

\svs
Let~$L$ be a pseudoline arrangement, and~$\chi$ be a cut of~$L$. We index by~$\ell_1,\dots,\ell_n$ the pseudolines of~$L$ in the order in which they cross~$\chi$ (it is well defined, up to a complete inversion).

We define the \defn{$k$-upper $\chi$-horizon set} of~$L$ to be the set~$\UU^k_\chi(L)$ of crossing points~$\ell_\alpha\wedge\ell_\beta$, with~$1\le \alpha<\beta\le n$, such that there is no~$\gamma_1,\dots,\gamma_k$ satisfying~$\alpha<\gamma_1<\dots<\gamma_k$ and $\ell_\alpha\wedge\ell_{\gamma_i}\cle_\chi\ell_\alpha\wedge\ell_\beta$ for all~$i\in[k]$. In other words, on each pseudoline~$\ell_\alpha$ of~$L$, the set~$\UU^k_\chi(L)$ consists of the smallest~$k$ crossing points of the form~$\ell_\alpha\wedge\ell_\beta$, with~$\alpha<\beta$.

Similarly, define the \defn{$k$-lower $\chi$-horizon set} of~$L$ to be the set~$\LL^k_\chi(L)$ of crossing points~${\ell_\alpha\wedge\ell_\beta}$, with ${1\le \alpha<\beta\le n}$, such that there is no $\delta_1,\dots,\delta_k$ satisfying~$\delta_1<\dots<\delta_k<\beta$ and ${\ell_\beta\wedge\ell_{\delta_i}\cle_\chi\ell_\alpha\wedge\ell_\beta}$ for all $i\in[k]$. On each pseudoline~$\ell_\beta$ of~$L$, the set~$\LL^k_\chi(L)$ consists of the smallest~$k$ crossing points of the form~$\ell_\alpha\wedge\ell_\beta$, with~$\alpha<\beta$.

Finally, we define the set~$\GG^k_\chi(L)$ to be the set of crossing points~$\ell_\alpha\wedge\ell_\beta$, with $1\le \alpha<\beta\le n$, such that there is no~$\gamma_1,\dots,\gamma_k$ and~$\delta_1,\dots,\delta_k$ satisfying:
\begin{enumerate}[(i)]
\item $\alpha<\gamma_1<\dots<\gamma_k$, $\delta_1<\dots<\delta_k<\beta$, and~$\delta_k<\gamma_1$;~and
\item $\ell_\alpha\wedge\ell_{\gamma_i}\cle_\chi\ell_\alpha\wedge\ell_\beta$ and~$\ell_\beta\wedge\ell_{\delta_i}\cle_\chi\ell_\alpha\wedge\ell_\beta$ for all $i\in[k]$.
\end{enumerate}

Obviously, the sets~$\UU^k_\chi(L)$ and~$\LL^k_\chi(L)$ are both contained in~$\GG^k_\chi(L)$.

\begin{example}
In \fref{fig:greedyhorizon}, we have labeled the vertices of the pseudoline arrangement~$L$ of \fref{fig:dual}(b) with different geometric tags according to their status: 
\begin{itemize}
\item[($\triangle$)] crossing points of the $k$-upper $\chi$-horizon set~$\UU_\chi^k(L)$ are represented by up triangles~$\triangle$;
\item[($\vartriangledown$)] crossing points of the $k$-lower $\chi$-horizon set~$\LL_\chi^k(L)$ are represented by down triangles~$\vartriangledown$;
\item[({\Large $\davidsstar$})] crossing points in both~$\UU_\chi^k(L)$ and~$\LL_\chi^k(L)$) are represented by up and down triangles~{\Large $\davidsstar$};
\item[($\Box$)] crossing points of~$\GG_\chi^k(L)$ but neither in~$\UU_\chi^k(L)$, nor in~$\LL_\chi^k(L)$ are represented by squares~$\Box$.
\end{itemize}
Observe that the remaining vertices are exactly the crossing points of the \greedy{\chi} \pt{k} of~$L$.

\begin{figure}
	\capstart
	\centerline{\includegraphics[scale=1]{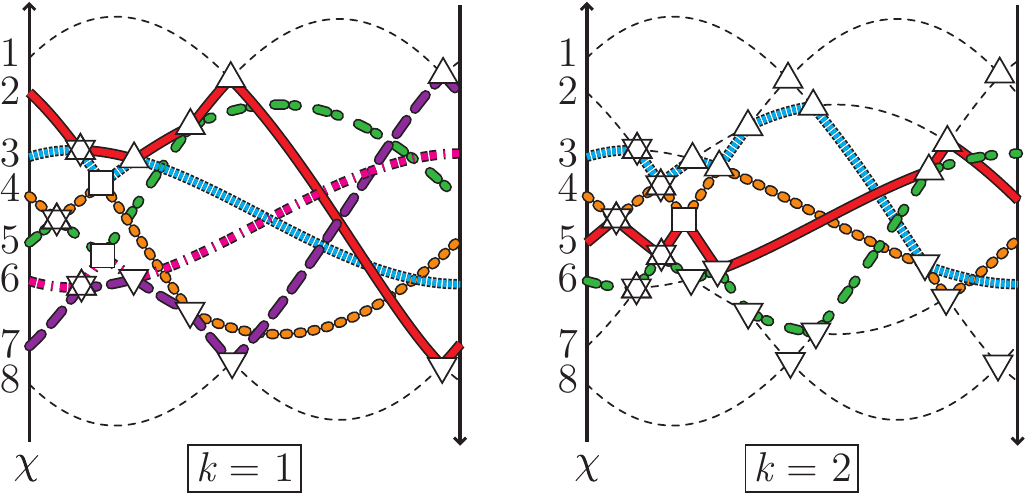}}
	\caption[The sets~$\UU^k_\chi(L)$,~$\LL^k_\chi(L)$, and~$\GG^k_\chi(L)$]{The sets $\UU^k_\chi(L)$, $\LL^k_\chi(L)$, and~$\GG^k_\chi(L)$ for the pseudoline arrangement of \fref{fig:dual}(b) and~${k\in\{1,2\}}$. The underlying \pt{k} is the greedy \pt{k} of~$L$.}
	\label{fig:greedyhorizon}
\end{figure}

\end{example}

\begin{example}\label{exm:convex}
We consider the arrangement~$C_n^*$ of~$n$ pseudolines in convex position. Let~$z$ be a vertex on the upper hull of its support,~$F \eqdef \set{z'}{z\cle z'}$ denote the filter generated by~$z$, and~$\chi$ denote the corresponding cut (see \fref{fig:convex}). It is easy to check that:
\begin{enumerate}[(i)]
\item $\UU^k_\chi(C_n^*)=\set{\ell_\alpha\wedge\ell_\beta}{1\le \alpha\le n\text{ and }\alpha<\beta\le \alpha+k}$;
\item $\LL^k_\chi(C_n^*)=\set{\ell_\alpha\wedge\ell_\beta}{1\le \alpha\le k\text{ and }\alpha< j\le n}$;~and
\item $\GG^k_\chi(C_n^*)=\UU^k_\chi(\cC_n)\cup\LL^k_\chi(\cC_n)$.
\end{enumerate}
Observe again that the remaining vertices are exactly the crossing points of the \greedy{\chi} \pt{k} of~$C_n^*$.

\begin{figure}
	\capstart
	\centerline{\includegraphics[scale=1]{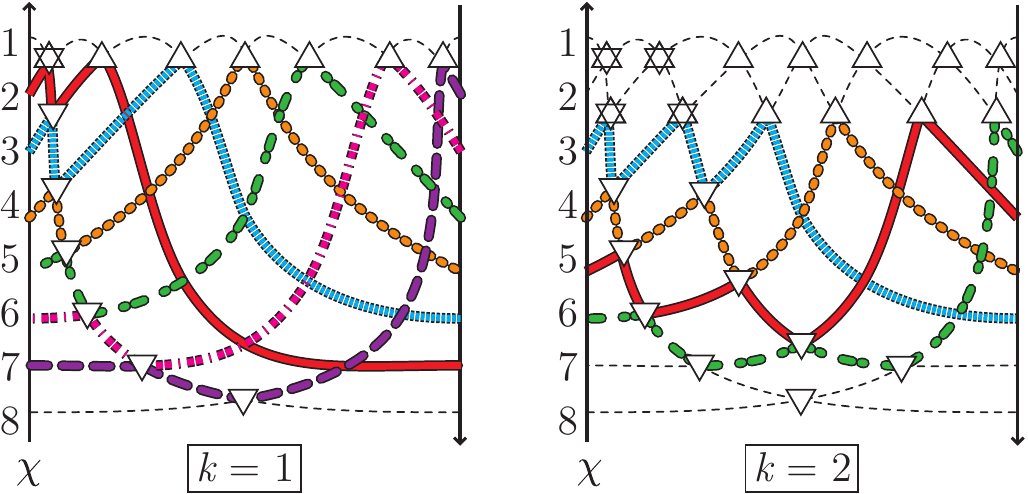}}
	\caption[The sets $\UU^k_\chi(C_8^*)$, $\LL^k_\chi(C_8^*)$, and $\GG^k_\chi(C_8^*)$]{The sets $\UU^k_\chi(C_8^*)$, $\LL^k_\chi(C_8^*)$, and $\GG^k_\chi(C_8^*)$ for the arrangement~$C_8^*$ of~$8$ pseudolines in convex position, and~$k\in\{1,2\}$. The underlying \pt{k} is the greedy \pt{k} of~$C_8^*$.}
	\label{fig:convex}
\end{figure}

\end{example}

Theorem~\ref{theo:horizon} extends this observation to all pseudoline arrangements, using convex position as a starting point for a proof by mutation.

\begin{theorem}\label{theo:horizon}
For any pseudoline arrangement~$L$ with no contact point, and any cut~$\chi$ of~$L$, the sets~$V(\Gamma_\chi^k(L))$ and~$\GG^k_\chi(L)$ coincide.
\end{theorem}

\begin{figure}[b]
	\capstart
	\centerline{\includegraphics[scale=1]{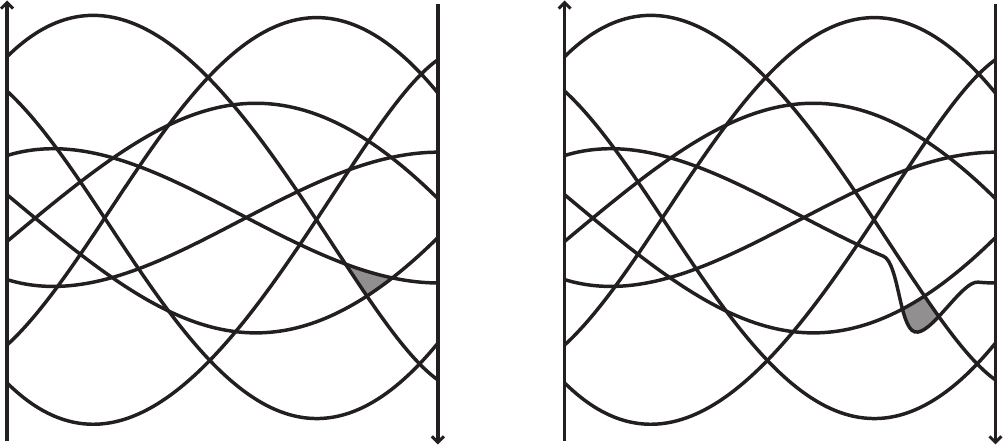}}
	\caption[A mutation in a pseudoline arrangement]{A mutation in the arrangement of \fref{fig:dual}(b).}
	\label{fig:mutation}
\end{figure}

The proof of this theorem works by mutation. A \defn{mutation} is a local transformation of an arrangement~$L$ that only inverts one triangular face of~$L$. More precisely, it is a homotopy of arrangements during which only one curve~$\ell\in L$ moves, sweeping a single vertex of the remaining arrangement~$L\ssm\{\ell\}$ (see \fref{fig:mutation}).

If~$P$ is a point set of a topological plane, mutating an empty triangle~$p^*q^*r^*$ of~$P^*$ by sweeping the vertex~$q^*\wedge r^*$ with the pseudoline~$p^*$ corresponds in the primal to moving~$p$ a little bit such that only the orientation of the triangle~$pqr$ changes.

The graph of mutations on pseudoline arrangements is known to be connected: any two pseudoline arrangements (with no contact points and the same number of pseudolines) are homotopic via a finite sequence of mutations (followed by a homeomorphism). In fact, one can even avoid mutations of triangles that cross a given cut of~$L$:

\begin{proposition}\label{prop:mutation}
Let~$L$ and~$L'$ be two pseudoline arrangements of~$\cM$ (with no contact points and the same number of pseudolines) and~$\chi$ be a cut of both of~$L$ and~$L'$. There is a finite sequence of mutations of triangles disjoint from~$\chi$ that transforms~$L$ into~$L'$.
\end{proposition}

\begin{proof}
We prove that any arrangement~$L$ of~$n$ pseudolines can be transformed into the arrangement~$C_n^*$ of $n$ pseudolines in convex position (see \fref{fig:convex}).

Let~$\ell_1,\dots,\ell_n$ denote the pseudolines of~$L$ (ordered by their crossings with~$\chi$). Let~$\Delta$ denote the triangle formed by~$\chi$, $\ell_1$~and~$\ell_2$. If there is a vertex of the arrangement~$L\ssm\{\ell_1,\ell_2\}$ inside~$\Delta$, then there is a triangle of the arrangement~$L$ inside~$\Delta$ and adjacent to~$\ell_1$ or~$\ell_2$. Mutating this triangle reduces the number of vertices of~$L\ssm\{\ell_1,\ell_2\}$ inside~$\Delta$ such that after some mutations, there is no more vertex inside~$\Delta$. If~$\Delta$ is intersected by pseudolines of~$L\ssm\{\ell_1,\ell_2\}$, then there is a triangle inside~$\Delta$ formed by~$\ell_1$,~$\ell_2$~and one of these intersecting pseudolines (the one closest to~$\ell_1\wedge\ell_2$). Mutating this triangle reduces the number of pseudolines intersecting~$\Delta$. Thus, after some mutations,~$\Delta$~is a triangle of the arrangement~$L$.

Repeating these arguments, we prove that for all ${i\in\{2,\dots,n-1\}}$ and after some mutations, $\ell_i$, $\ell_1$, $\ell_{i+1}$ and~$\chi$ delimit a face of the arrangement~$L$. Thus, one of the two topological disk delimited by~$\chi$ and~$\ell_1$ contains no more vertex of~$L$, and the proof is then straightforward by induction.
\end{proof}

\begin{figure}[b]
	\capstart
	\centerline{\includegraphics[scale=1]{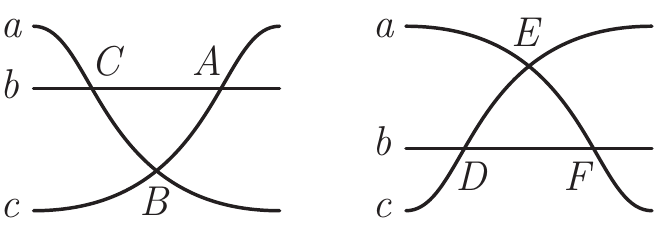}}
	\caption[A local image of a mutation]{A local image of a mutation.}
	\label{fig:mutationlocal}
\end{figure}

Let~$\nabla$ be a triangle of $L$ not intersecting~$\chi$. Let~$L'$ denote the pseudoline arrangement obtained from~$L$ by mutating the triangle~$\nabla$ into the inverted triangle~$\Delta$.
Let~$a<b<c$ denote the indices of the pseudolines~$\ell_a,\ell_b$ and~$\ell_c$ that form~$\nabla$ and~$\Delta$. In~$\nabla$, we denote~$A=\ell_b\wedge\ell_c$, $B=\ell_a\wedge\ell_c$ and~$C=\ell_a\wedge\ell_b$; similarly, in~$\Delta$, we denote~$D=\ell_b\wedge\ell_c$, $E=\ell_a\wedge\ell_c$ and~$F=\ell_a\wedge\ell_b$ (see \fref{fig:mutationlocal}).

\begin{lemma}\label{lem:mutation}
With these notations, the following properties hold:
\begin{enumerate}[(i)]
\item $B\in\UU_\chi^k(L)\Leftrightarrow C\in\UU_\chi^{k-1}(L)\Leftrightarrow E\in\UU_\chi^{k-1}(L')\Leftrightarrow F\in\UU_\chi^k(L')$\\
$A\in\UU_\chi^k(L)\Leftrightarrow D\in\UU_\chi^k(L')$\\
$E\in\LL_\chi^k(L')\Leftrightarrow D\in\LL_\chi^{k-1}(L')\Leftrightarrow B\in\LL_\chi^{k-1}(L)\Leftrightarrow A\in\LL_\chi^k(L)$\\
$F\in\LL_\chi^k(L')\Leftrightarrow C\in\LL_\chi^k(L)$
\item $C\in\UU_\chi^k(L)\Rightarrow A\in\UU_\chi^k(L)$\\
$E\in\UU_\chi^k(L')\Rightarrow D\in\UU_\chi^k(L')$\\
$D\in\LL_\chi^k(L')\Rightarrow F\in\LL_\chi^k(L')$\\
$B\in\LL_\chi^k(L)\Rightarrow C\in\LL_\chi^k(L)$\\
$B\in\GG_\chi^k(L)\Rightarrow C\in\GG_\chi^k(L)\wedge D\in\GG_\chi^k(L')\wedge F\in\GG_\chi^k(L')$\\
$E\in\GG_\chi^k(L')\Rightarrow D\in\GG_\chi^k(L')\wedge C\in\GG_\chi^k(L)\wedge A\in\GG_\chi^k(L)$
\item $C\in\GG_\chi^k(L)\Rightarrow A\in\UU_\chi^k(L)\vee C\in\LL_\chi^k(L)$\\
$D\in\GG_\chi^k(L')\Rightarrow F\in\LL_\chi^k(L')\vee D\in\UU_\chi^k(L')$\\
$A\in\GG_\chi^k(L)\Rightarrow A\in\UU_\chi^k(L)\vee C\in\LL_\chi^{k-1}(L)$\\
$F\in\GG_\chi^k(L')\Rightarrow F\in\LL_\chi^k(L')\vee D\in\UU_\chi^{k-1}(L')$
\item $C\in\GG_\chi^k(L)\wedge E\notin\GG_\chi^k(L')\Rightarrow A\notin\GG_\chi^k(L)$\\
$D\in\GG_\chi^k(L')\wedge B\notin\GG_\chi^k(L)\Rightarrow F\notin\GG_\chi^k(L')$
\end{enumerate}
\end{lemma}

\begin{proof}
By symmetry, it is enough to prove the first line of each of the four points of the lemma.

Properties of point~(i) directly come from the definitions. For example, all the assertions of the first line are false if and only if there exist~$\gamma_1,\dots,\gamma_{k-1}$ with~$a<\gamma_1<\dots<\gamma_{k-1}$ and, for all~$i\in[k-1]$, $\ell_a\wedge\ell_{\gamma_i}\cle_\chi C$ (or equivalently~$\ell_a\wedge\ell_{\gamma_i}\cle_\chi E$).

We derive point~(ii) from the following observation: if~$\gamma>b$ and if~$\ell_b\wedge\ell_\gamma\cle_\chi C$, then~$\gamma>a$ and~$\ell_a\wedge\ell_\gamma\cle_\chi B$.

For point~(iii), assume that~$A\notin\UU_\chi^k(L)$ and~$C\notin\LL_\chi^k(L)$. Then there exist~$\gamma_1,\dots,\gamma_k$ and~$\delta_1,\dots,\delta_k$ such that~${\delta_1<\dots<\delta_k<b<\gamma_1<\dots<\gamma_k}$ and, for all $i\in[k]$, 
$\ell_b\wedge\ell_{\gamma_i}\cle_\chi A$ (and therefore $\ell_a\wedge\ell_{\gamma_i}\cle_\chi C$) and $\ell_b\wedge\ell_{\delta_i}\cle_\chi C$. Thus~$C\notin\GG_\chi^k(L)$.

Finally, assume that~$C\in\GG_\chi^k(L)$ and~$E\notin\GG_\chi^k(L')$. Then, there exist~$\gamma_1,\dots,\gamma_k$ and $\delta_1,\dots,\delta_k$ such that $a<\gamma_1<\dots<\gamma_k$, $\delta_1<\dots<\delta_k<c$, $\delta_k<\gamma_1$, and for all $i\in[k]$, $\ell_a\wedge\ell_{\gamma_i}\cle_\chi E$ and $\ell_c\wedge\ell_{\delta_i}\cle_\chi E$. Since $C\in\GG_\chi^k(L)$, we have $\delta_k>b$. Thus $b<\gamma_1<\dots<\gamma_k$ and for all $i\in[k]$, $\ell_b\wedge\ell_{\gamma_i}\cle_\chi A$ and $\ell_c\wedge\ell_{\delta_i}\cle_\chi A$. This implies that~$A\notin\GG_C^k(L)$.
\end{proof}

\begin{figure}
	\capstart
	\centerline{\includegraphics[scale=1]{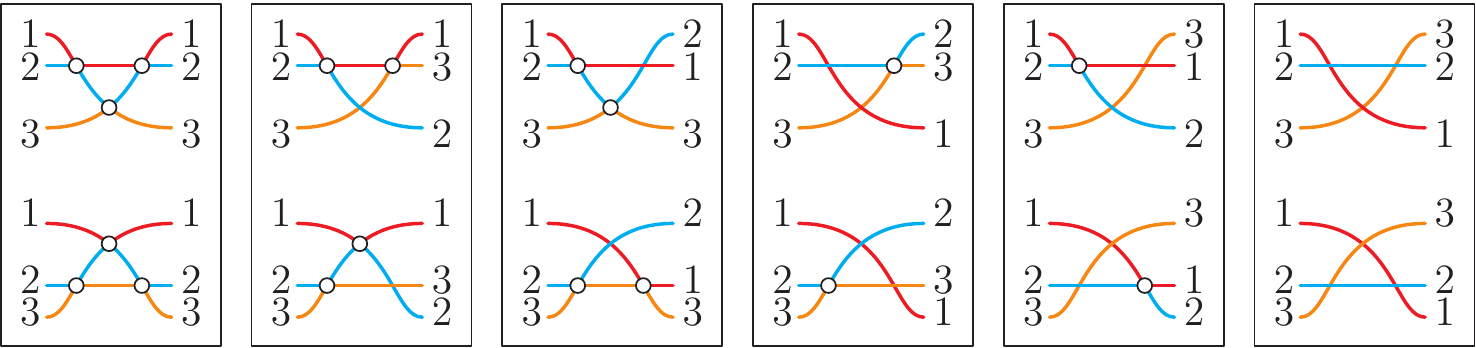}}
	\caption[The six possible cases for the mutation of the greedy \pt{k}]{The six possible cases for the mutation of the greedy \pt{k}.}
	\label{fig:mutationlocal6sol}
\end{figure}

We are now ready to establish the proof of Theorem~\ref{theo:horizon}:

\begin{proof}[Proof of Theorem~\ref{theo:horizon}]
The proof works by mutation. We already observed the result when the pseudoline arrangement is in convex position (see Example~\ref{exm:convex} and \fref{fig:convex}). Proposition~\ref{prop:mutation} ensures that any pseudoline arrangement can be reached from this convex configuration by mutations of triangles not intersecting~$\chi$. Thus, it is sufficient to prove that such a mutation preserves the property.

Assume that~$L$ is a pseudoline arrangement and~$\chi$ is a cut of~$L$, for which the result holds. Let~$\nabla$ be a triangle of~$L$ not intersecting~$\chi$. Let~$L'$ denote the pseudoline arrangement obtained from~$L$ by mutating the triangle~$\nabla$ into the inverted triangle~$\Delta$. Let~$A,B,C$~and~$D,E,F$ denote the vertices of~$\nabla$~and~$\Delta$ as indicated in \fref{fig:mutationlocal}. 

If~$v$ is a vertex of the arrangement~$L'$ different from~$D,E,F$, then:
$$v\in V(\Gamma^k_\chi(L'))\Leftrightarrow v\in V(\Gamma^k_\chi(L))\Leftrightarrow v\in\GG^k_\chi(L)\Leftrightarrow v\in\GG^k_\chi(L').$$

Thus, we only have to prove the equivalence when~$v\in\{D,E,F\}$. The proof is a (computational) case analysis: using the properties of Lemma~\ref{lem:mutation} as boolean equalities relating the boolean variables defined by ``$X\in\mathbb{Y}_\chi^p(L)$'' (where~$X\in\{A,B,C,D,E,F\}$, $\mathbb{Y}\in\{\UU,\LL,\GG\}$, and ${p\in\{k-1,k\}}$), we have written a short boolean satisfiability program which affirms that:
\begin{enumerate}[(i)]
\item either $\{A,B,C\}\subset\GG_\chi^k(L)$ and $\{D,E,F\}\subset\GG_\chi^k(L')$;
\item or $\{A,B,C\}\cap\GG_\chi^k(L)=\{A,C\}$ and $\{D,E,F\}\cap\GG_\chi^k(L')=\{D,E\}$;
\item or $\{A,B,C\}\cap\GG_\chi^k(L)=\{B,C\}$ and $\{D,E,F\}\cap\GG_\chi^k(L')=\{D,F\}$;
\item or $\{A,B,C\}\cap\GG_\chi^k(L)=\{A\}$ and $\{D,E,F\}\cap\GG_\chi^k(L')=\{D\}$;
\item or $\{A,B,C\}\cap\GG_\chi^k(L)=\{C\}$ and $\{D,E,F\}\cap\GG_\chi^k(L')=\{F\}$;
\item or $\{A,B,C\}\cap\GG_\chi^k(L)=\emptyset$ and $\{D,E,F\}\cap\GG_\chi^k(L')=\emptyset$.
\end{enumerate}

It is easy to check that these six cases correspond to sorting the six possible permutations of~$\{1,2,3\}$ on~$\nabla$ and~$\Delta$ (see \fref{fig:mutationlocal6sol}). Consequently, if~$V(\Gamma_\chi^k(L))\cap\{A,B,C\}=\GG_\chi^k(L)\cap\{A,B,C\}$, then we have $V(\Gamma_\chi^k(L'))\cap\{D,E,F\}=\GG_\chi^k(L')\cap\{D,E,F\}$, which finishes the proof.
\end{proof}

\begin{figure}[h]
	\capstart
	\centerline{\includegraphics[scale=1]{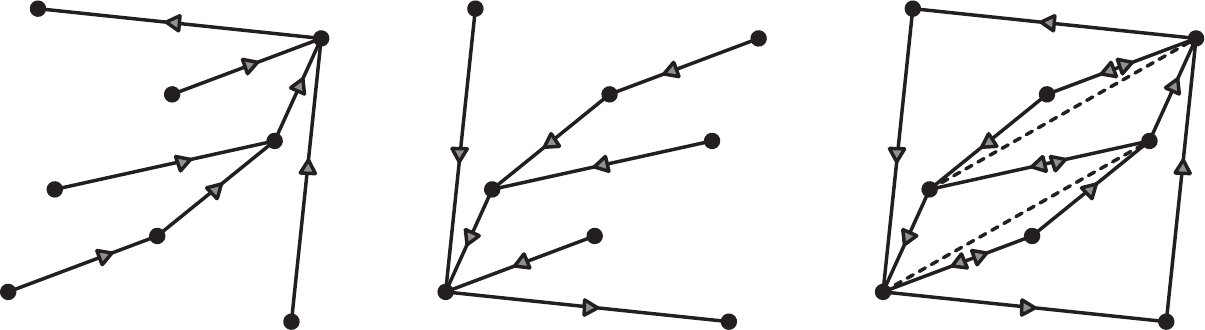}}
	\caption[Horizon trees]{The upper horizon tree (left), the lower horizon tree (middle), and the greedy pseudotriangulation (right) of the point set of \fref{fig:dual}(a). The dashed edges in the greedy pseudotriangulation are not in the horizon trees.}
	\label{fig:horizon}
\end{figure}

Let us finish this discussion by recalling the interpretation of the horizon sets when~${k=1}$.
Let~$P$ be a finite point set. For any~$p\in P$, let~$u(p)$ denote the point~$q$ that minimizes the angle~$(Ox,pq)$ among all points of~$P$ with~$y_p<y_q$ (by convention, for the higher point~$p$ of~$P$,~$u(p)=p$). The \defn{upper horizon tree} of~$P$ is the set~$\UU(P)=\set{pu(p)}{p\in P}$. The \defn{lower horizon tree}~$\LL(P)$ of $P$ is defined symmetrically. See \fref{fig:horizon}.

Choosing a cut~$\chi$ of~$P^*$ corresponding to the point at infinity~$(-\infty,0)$ makes coincide primal and dual definitions of horizon sets: we have $\UU^1_\chi(P^*)=\UU(P)^*$~and~$\LL^1_\chi(P^*)=\LL(P)^*$.

In~\cite{p-htvpt-97}, Pocchiola observed that the set~$\UU(P)\cup\LL(P)$ of edges can be completed into a pseudotriangulation of~$P$ just by adding the sources of the faces of~$P^*$ intersected by the cut~$\chi$. The obtained pseudotriangulation is our \greedy{\chi} \pt{1}~$\Gamma^1_\chi(P^*)$.


\subsection{Multipseudotriangulations of double pseudoline arrangements}\label{subsec:furthertopics:dpl}

In this section, we deal with double pseudoline arrangements, \ie duals of sets of disjoint convex bodies. Definitions and properties of \mpt{}s naturally extend to these objects.

\subsubsection{Definitions}

\begin{figure}
	\capstart
	\centerline{\includegraphics[scale=1]{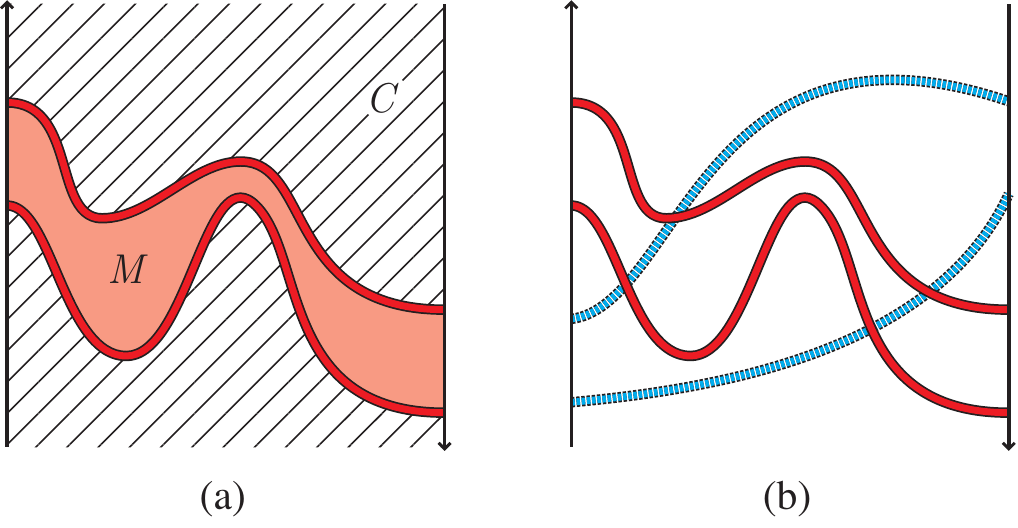}}
	\caption[A double pseudoline and a double pseudoline arrangement in the M\"obius strip]{(a) A double pseudoline. (b) An arrangement of~$2$ double pseudolines.}
	\label{fig:dpl}
\end{figure}

A simple closed curve in the M\"obius strip can be:
\begin{enumerate}[(i)]
\item either contractible (homotopic to a point);
\item or non separating, or equivalently homotopic to a generator of the fundamental group of~$\cM$: it is a pseudoline;
\item or separating and non-contractible, or equivalently homotopic to the double of a generator of the fundamental group of~$\cM$: it is called a \defn{double pseudoline} (see \fref{fig:dpl}(a)).
\end{enumerate}

The complement of a double pseudoline~$\ell$ has two connected components: the bounded one is a M\"obius strip~$M_\ell$ and the unbounded one is an open cylinder~$C_\ell$ (see \fref{fig:dpl}(a)).

\begin{figure}[b]
	\capstart
	\centerline{\includegraphics[scale=1]{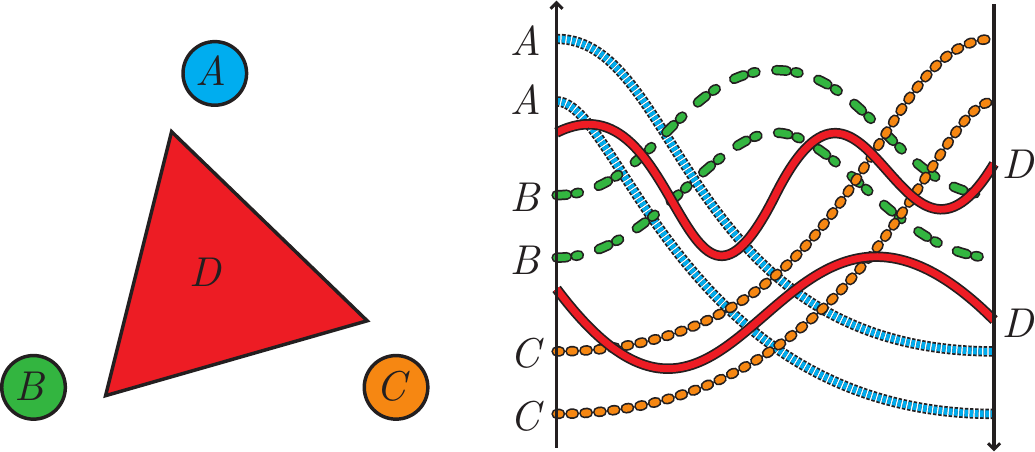}}
	\caption[A configuration of~$4$ disjoint convex bodies and its dual double pseudoline arrangement]{A configuration of four disjoint convex bodies and its dual double pseudoline arrangement.}
	\label{fig:dplactrexm1}
\end{figure}

\begin{figure}[b]
	\capstart
	\centerline{\includegraphics[scale=1]{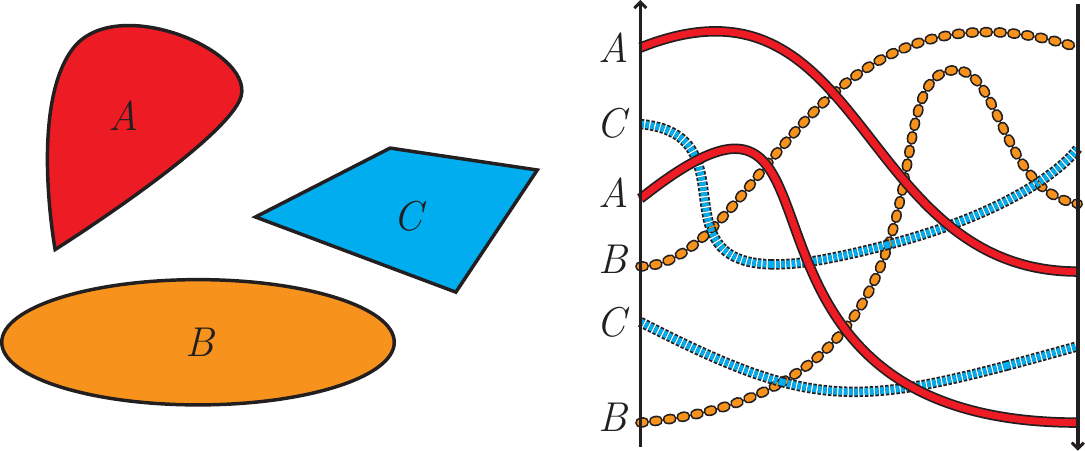}}
	\caption[A configuration of~$3$ disjoint convex bodies and its dual double pseudoline arrangement]{A configuration of three disjoint convex bodies and its dual double pseudoline arrangement.}
	\label{fig:dpla}
\end{figure}

The canonical example of a double pseudoline is the set~$C^*$ of tangents to a convex body~$C$ of the plane. Observe also that the $p$th level of a pseudoline arrangement is a double pseudoline. If~$C$ is a convex body of the plane, then the M\"obius strip~$M_{C^*}$ corresponds to lines that pierce~$C$, while~$C_{C^*}$ corresponds to lines disjoint from~$C$. If~$C$ and~$C'$ are two disjoint convex bodies, the two corresponding double pseudolines~$C^*$ and~$C'^*$ cross four times (see \fref{fig:dplactrexm1} and~\fref{fig:dpla}). Each of these four crossings corresponds to one of the four \defn{bitangents} (or \defn{common tangents}) between~$C$~and~$C'$.

\begin{definition}[\cite{hp-adp-08}]
A \defn{double pseudoline arrangement} is a finite set of double pseudolines such that any two of them have exactly four intersection points, cross transversally at these points, and induce a cell decomposition of the M\"obius strip (\ie the complement of their union is a union of topological disks, together with the external cell).
\end{definition}

Given a set~$Q$ of disjoint convex bodies in the plane (or in any topological plane), its dual $Q^* \eqdef \set{C^*}{C\in Q}$ is an arrangement of double pseudolines (see \fref{fig:dplactrexm1} and \fref{fig:dpla}). Furthermore, as for pseudoline arrangements, any double pseudoline arrangement can be represented by (\ie is the dual of) a set of disjoint convex bodies in a topological plane~\cite{hp-adp-08}.

\begin{figure}[b]
	\capstart
	\centerline{\includegraphics[scale=1]{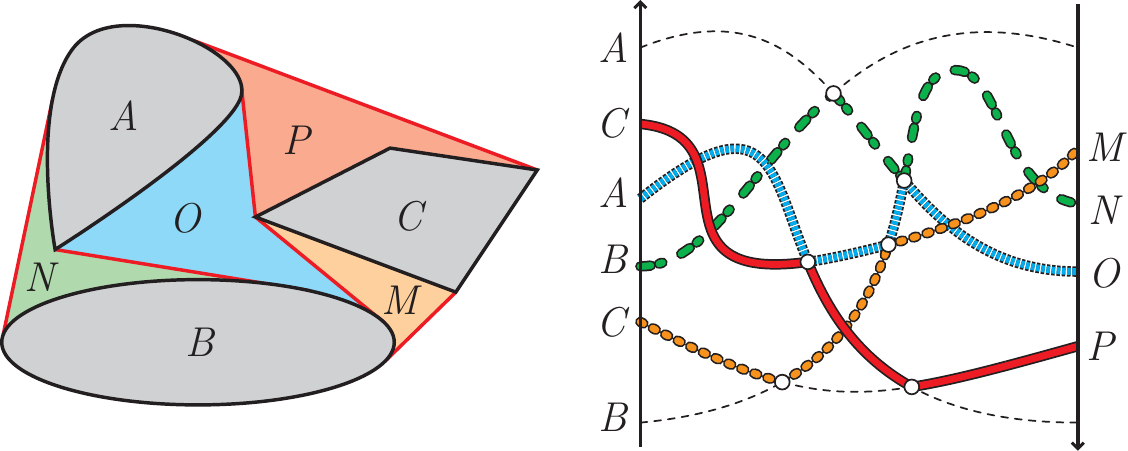}}
	\caption[A pseudotriangulation of a set of disjoint convex bodies]{A pseudotriangulation of the set of disjoint convex bodies of \fref{fig:dpla}.}
	\label{fig:pseudotriangulation}
\end{figure}

In this paper, we only consider \defn{simple} arrangements of double pseudolines. Defining the \defn{support}, the \defn{levels}, and the \defn{kernels} of double pseudoline arrangements as for pseudoline arrangements, we can extend \mpt{}s to double pseudoline arrangements (see \fref{fig:pseudotriangulation}):

\begin{definition}
A \pt{k} of a double pseudoline arrangement~$L$ is a pseudoline arrangement supported by the \kkernel{k} of~$L$.
\end{definition}

All the properties related to flips developed in Section~\ref{sec:enumeration} apply in this context. In the end of this section, we only revisit the properties of the primal of a \mpt{} of a double pseudoline arrangement.

\subsubsection{Elementary properties}

Let~$Q$ be a set disjoint convex bodies in general position in the plane and~$Q^*$ be its dual arrangement. Let~$\Lambda$ be a \pt{k} of~$Q^*$, $V(\Lambda)$~denote all crossing points of~$Q^*$ that are not crossing points of~$\Lambda$, and~$E$ denote the corresponding set of bitangents of~$Q$. As in Subsection~\ref{subsec:mpt:pointedcrossing}, we discuss the properties of the primal configuration~$E$:

\begin{lemma}
The set~$E$ has~$4|Q|k-|Q|-2k^2-k$ edges.
\end{lemma}

\begin{proof}
The number of edges of~$E$ is the number of crossing points of~$Q^*$ minus the number of crossing points of~$\Lambda$, \ie
\begin{align*}
|E| & =4{|Q^*| \choose 2}-{|\Lambda| \choose 2}=4{|Q| \choose 2}-{2|Q|-2k \choose 2} \\
& =4|Q|k-|Q|-2k^2-k.\qedhere
\end{align*}
\end{proof}

We now discuss pointedness. For any convex body~$C$ of~$Q$, we arbitrarily choose a point~$p_C$ in the interior of~$C$, and we consider the set~$X_C$ of all segments between~$p_C$ and a sharp boundary point of~$C$. We denote by~$X \eqdef \bigcup_{C\in Q} X_C$ the set of all these segments.

\begin{lemma}\label{lem:alternationfree}
The set $E\cup X$ cannot contain a \kalter{k}.
\end{lemma}

\begin{figure}
	\capstart
	\centerline{\includegraphics[scale=1]{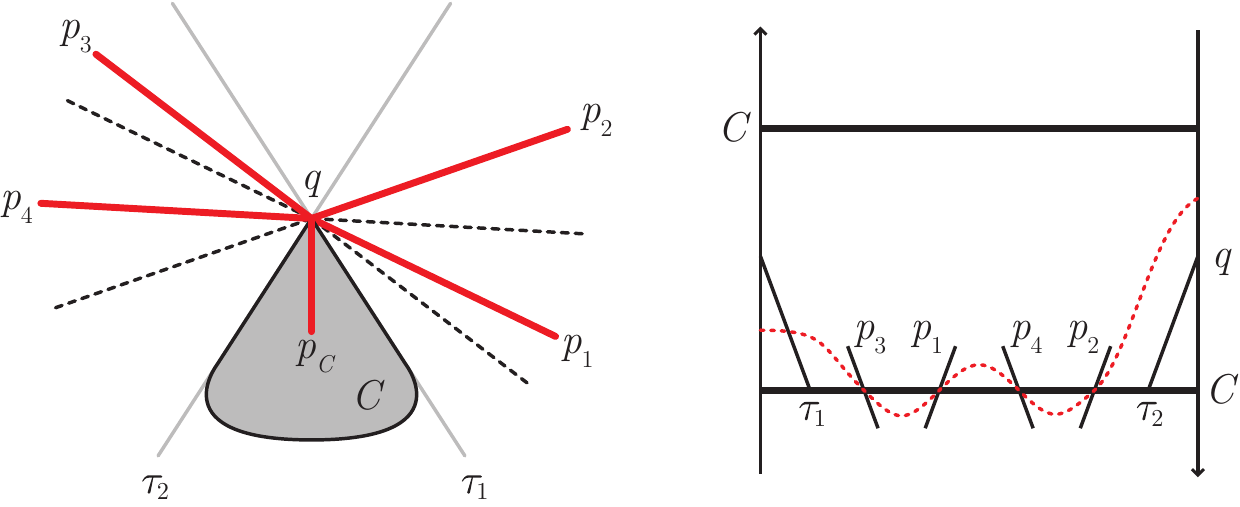}}
	\caption[A \kalter{k} at a sharp point]{A \kalter{k} at a sharp point.}
	\label{fig:sharp}
\end{figure}

\begin{proof}
Let~$C$ be a convex body of~$Q$, let $q$~be a sharp point of~$C$ and let $F \eqdef \set{[p_i,q]}{i\in[2k]}$ be a set of edges incident to~$q$ such that $\{[p_C,q]\}\cup F$ is a \kalter{k}. We prove that $F$ is not contained in~$E$. Indeed, the dual pseudolines~$\set{p_i^*}{i\in[2k]}$ intersect alternately the double pseudoline~$C^*$ between the tangents to~$C$ at~$q$ (see \fref{fig:sharp}). This ensures the existence of a witness pseudoline~$\ell$ that separates all the contact points~$p_i^*\wedge C^*$, while crossing~$C^*$ exactly~$2k$ times and the other double pseudolines of~$Q^*$ exactly has~$q^*$ does. (As usual, we obtain it by a perturbation of the pseudoline~$q^*$.) Counting the crossings of~$\ell$ with~$Q^*$ and~$\Lambda$ respectively, we obtain:
\begin{enumerate}[(i)]
\item $\ell$~crosses~$Q^*$ exactly~$2|(Q\ssm\{C\})^*|+2k=2|Q|+2k-2$ times;
\item $\ell$~crosses~$\Lambda$ at least~$|\Lambda|=2|Q|-2k$ times;
\item for each of the points~$p_i\wedge q^*$, replacing the crossing point by a contact point removes two crossings with~$\ell$.\qedhere
\end{enumerate}
\end{proof}

\subsubsection{Stars}

If~$\lambda$ is a pseudoline of~$\Lambda$, we call \defn{star} the envelope~$S(\lambda)$ of the primal lines of the points of~$\lambda$. The star~$S(\lambda)$ contains:
\begin{enumerate}[(i)]
\item all bitangents~$\tau$ between two convex bodies of~$Q$ such that~$\tau$ is a contact point of~$\lambda$; and
\item all convex arcs formed by the tangent points of the lines covered by~$\lambda$ with the convex bodies of~$Q$.
\end{enumerate}
This star is a (non-necessarily simple) closed curve. We again have bounds on the number of \defn{corners} (\ie convex internal angles) of~$S(\lambda)$:

\begin{proposition}\label{prop:cornersconv}
The number of corners of~$S(\lambda)$ is odd and between $2k-1$ and~$4k|Q|-2k-1$.
\end{proposition}

\begin{proof}
In the case of double pseudoline arrangements, corners are even easier to characterize: a bitangent~$\tau$ between two convex~$C$ and~$C'$ of~$Q$ always defines two corners, one at each extremity. These corners are contained in one of the two stars adjacent to~$\tau$. Let~$\lambda$ be a pseudoline with a contact point at~$\tau$. In a neighborhood of~$\tau$, the pseudoline~$\lambda$ can be contained either in~$M_{C^*}$ or in~$C_{C^*}$. In the first case, the star~$S(\lambda)$ contains the corner formed by the bitangent~$\tau$ and the convex~$C$ (or possibly, by the bitangent~$\tau$ and another tangent to~$C$); while in the second case, it does not. (The same observation holds for~$C'$.)

\begin{figure}
	\capstart
	\centerline{\includegraphics[scale=1]{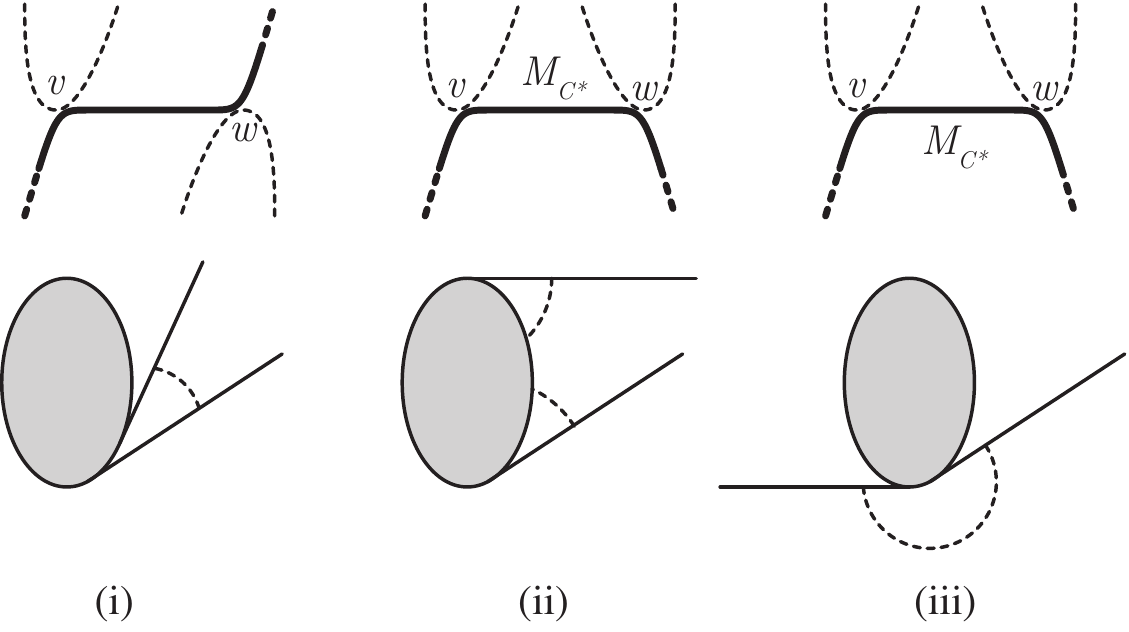}}
	\caption[Corners of a star]{The three possible situations for two consecutive contact points on~$\ell$.}
	\label{fig:corners}
\end{figure}

In other words, if the double pseudoline~$C^*$ supports a pseudoline~$\lambda$ between two contact points~$v$ and~$w$, then one of the three following situations occurs (see \fref{fig:corners}):
\begin{enumerate}[(i)]
\item either~$v$ and~$w$ lie on opposite sides of~$\lambda$; then exactly one of these contact points lies in~$M_{C^*}$, and~$S(\lambda)$ has one corner at~$C$.
\item or~$v$ and~$w$ both lie on~$M_{C^*}$, and~$S(\lambda)$ has two corners at~$C$.
\item or~$v$ and~$w$ both lie on~$C_{C^*}$, and~$S(\lambda)$ has no corners at~$C$.
\end{enumerate}
In particular, the number~$c=c(\lambda)$ of corners of~$\lambda$ is the number of situations~(i) plus twice the number of situations~(ii). This proves that~$c$ is odd and (at least) bigger than the number of opposite consecutive contact points of the pseudoline~$\lambda$.

In order to get a lower bound on this number, we construct (as in the proof of Proposition~\ref{prop:corners}) a witness pseudoline~$\mu$ that crosses~$\lambda$ between each pair of opposite contact points and passes on the opposite side of each contact point. It crosses~$\lambda$ at most~$c$ times and~$\Lambda\ssm\{\lambda\}$ exactly~$|\Lambda|-1$ times. Moreover, if~$\alpha$ is a pseudoline and~$\beta$ is a double pseudoline of~$\cM$, then either~$\alpha$ is contained in~$M_\beta$ and has no crossing with~$\beta$, or~$\alpha$ and~$\beta$ have an even number of crossings. Since~$\mu$ is a pseudoline and can be contained in at most one M\"obius strip~$M_C^*$ (for~$C\in Q$), the number of crossings of~$\mu$ with~$Q^*$ is at least~$2(|Q|-1)$. Thus, we obtain the lower bound $2(|Q|-1)\le2|Q|-2k-1+c$, \ie $c\ge 2k-1$.

From this lower bound, we obtain the upper bound: the total number of corners is at most twice the number of bitangents:
$$2(4k|Q|-|Q|-2k^2-k)\ge\sum_{\mu\in\Lambda} c(\mu)\ge c(\lambda)+(2|Q|-2k-1)(2k-1),$$
and we get~$c\le4|Q|k-2k-1$.
\end{proof}

When~$k=1$, we can even prove that all stars are pseudotriangles. Indeed, since any star has at least~$3$ corners, the upper bound calculus gives~$2(3|Q|-3)\ge c+3(2|Q|-3)$, \ie $c\le 3$.

\begin{figure}
	\capstart
	\centerline{\includegraphics[scale=1]{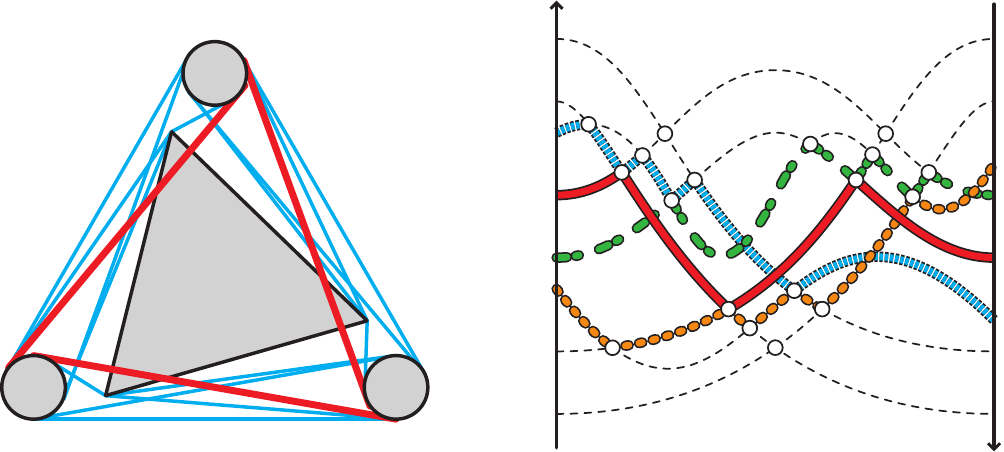}}
	\caption[A \pt{2} of the double pseudoline arrangement of \fref{fig:dplactrexm1}]{A \pt{2} of the double pseudoline arrangement of \fref{fig:dplactrexm1}. Observe that the bolded red star has only~$3$ corners.}
	\label{fig:dplactrexm2}
\end{figure}

For general~$k$, observe that contrary to the case of pseudoline arrangements, a star may have~$2k-1$ corners (see \fref{fig:dplactrexm2}).

\svs
Let us now give an analogue of Proposition~\ref{prop:decomposition}. For any point~$q$ in the plane, we denote by~$\eta^k(q)$ the number of crossings between~$q^*$ and the support of~$Q^*$ minus its first~$k$ levels. Let~$\delta^k(q) \eqdef \eta^k(q)/2-|Q|+k$. For any~$\lambda\in\Lambda(U)$ and any point~$q$ in the plane, we still denote by~$\sigma_\lambda(q)$ the \defn{winding number} of~$S(\lambda)$ around~$q$.

\begin{proposition}\label{prop:decompositionconv}
For any point~$q$ of the plane~$\delta^k(q)=\sum_{\lambda\in\Lambda} \sigma_\lambda(q)$.
\end{proposition}

\begin{proof}
Remember that if~$\tau_\lambda(q)$ denotes the number of intersection points between~$q^*$ and~$\lambda$, then~$\sigma_\lambda(q)=(\tau_\lambda(q)-1)/2$. Thus, we have
$$\eta^k(q)=\sum_{\lambda\in\Lambda} \tau_\lambda(q)=|\Lambda|+2\sum_{\lambda\in\Lambda}\sigma_\lambda(q),$$
\nobreak and we get the result since~$|\Lambda|=2|Q|-2k$.
\end{proof}

When~$k=1$, it is easy to see that~$\delta^1(q)$ is~$1$ if~$q$ is inside the \defn{free space} of the convex hull of~$Q$ (\ie in the convex hull of~$Q$, but not in~$Q$), and~$0$ otherwise. Remember that a \defn{pseudotriangulation} of~$Q$ is a pointed set of bitangents that decomposes the free space of the convex hull of~$Q$ into pseudotriangles~\cite{pv-ptta-96}. Propositions~\ref{prop:cornersconv} and ~\ref{prop:decompositionconv} provide, when~$k=1$, the following analogue of Theorem~\ref{theo:dualitypt}:

\begin{theorem}\label{theo:dualityptconv}
Let~$Q$ be a set of disjoint convex bodies (in general position) and~$Q^*$ denote its dual arrangement. Then:
\begin{enumerate}[(i)]
\item The dual arrangement~$T^* \eqdef \set{\Delta^*}{\Delta\text{ pseudotriangle of }T}$ of a pseudotriangulation~$T$ of~$Q$ is a \pt{1} of~$Q^*$.
\item The primal set of edges $$\quad E \eqdef \set{[p,q]}{p,q\in P,\; p^*\wedge q^*\text{ is not a crossing point of } \Lambda}$$ of a \pt{1}~$\Lambda$ of~$Q^*$ is a pseudotriangulation of~$Q$.
\end{enumerate}
\end{theorem}

Observe that at least two other arguments are possible to prove~(ii):
\begin{enumerate}
\item either comparing the degrees of the flip graphs as in our first proof of Theorem~\ref{theo:dualitypt};
\item or checking that all forbidden configurations of the primal (two crossing bitangents, a non-pointed sharp vertex, a non-free bitangent) may not appear in the dual, as in our second proof of Theorem~\ref{theo:dualitypt}.
\end{enumerate}

\begin{figure}[b]
	\capstart
	\centerline{\includegraphics[scale=1]{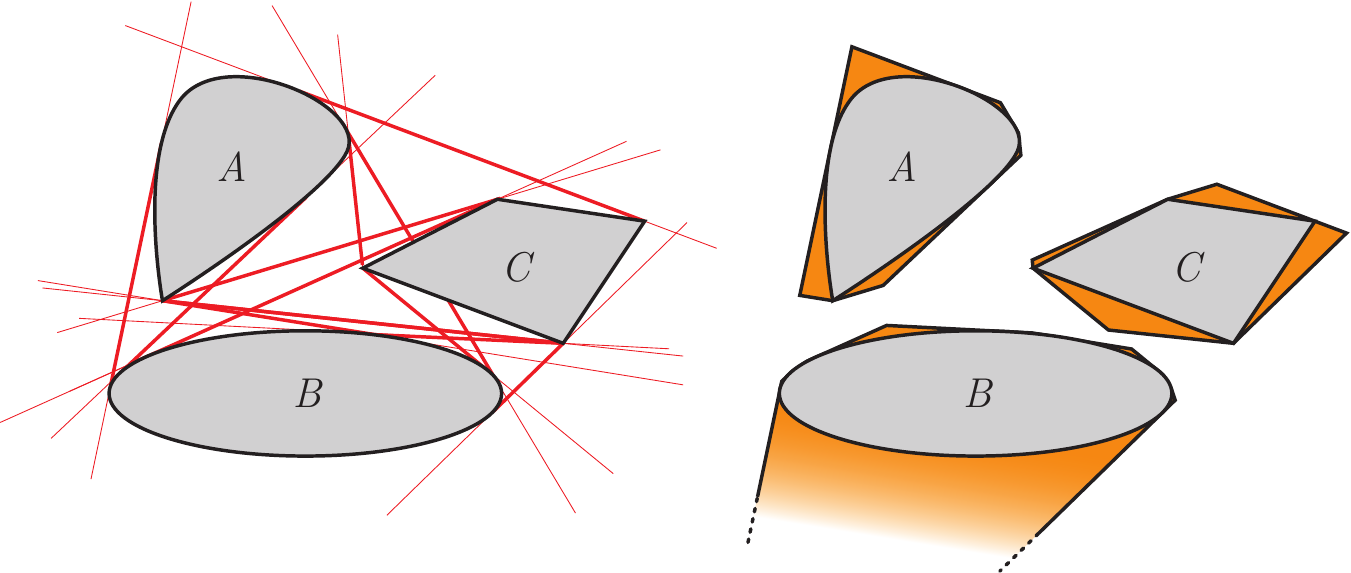}}
	\caption[The set of all bitangents to the arrangement of convex bodies of \fref{fig:dpla} and the corresponding maximal convex bodies]{The set of all bitangents to the arrangement of convex bodies of \fref{fig:dpla} and the corresponding maximal convex bodies.}
	\label{fig:maxconv}
\end{figure}

\svs
Let us finish this discussion about stars by interpreting the number~$\delta^k(q)$ for general~$k$ and for ``almost every'' point~$q$. For any convex body~$C$ of~$Q$, let~$\nabla C$ denote the intersection of all closed half-planes delimited by a bitangent between two convex bodies of~$Q$, and containing~$C$ (see \fref{fig:maxconv}). By definition, the bitangents between two convex bodies~$C$ and~$C'$ of~$Q$ coincide with the bitangent of~$\nabla C$ and~$\nabla C'$. Furthermore, the convex bodies~$\nabla C$~($C\in Q$) are maximal for this property. We denote~$\nabla Q \eqdef \set{\nabla C}{C\in Q}$ the set of maximal convex bodies of~$Q$.

For a point~$q$ outside~$\nabla Q$, the interpretation of~$\delta^k(q)$ is similar to the case of points. We call \defn{level} of a bitangent~$\tau$ the level of the corresponding crossing point in~$Q^*$. Given a point~$q$ outside~$\nabla Q$, the number~$\delta^k(q)$ is the number of bitangents of level~$k$ crossed by any (generic continuous) path from~$q$ to the external face (in the complement of~$\nabla Q$), counted positively when passing from the ``big'' side to the ``small side'', and negatively otherwise.


\section{Open questions}\label{sec:open}

We finish by a short presentation of some open questions that have arisen out of this work. Since the submission of this paper, some of these questions were (at least partially) answered in~\cite{s-npkt-11, ss-mfmpscsp-10, ps-bp-12}. We have decided to keep these questions in the discussion and to refer to the relevant articles in side remarks.

\medskip

\paragraph{\sc Primal of a \mpt{}.}

When $k=1$ or in the case of convex position, primals of \pt{k}s are characterized by simple non-crossing and pointedness conditions. 
For general $k$ and general position, we know that the primal of a \pt{k} is \kalter{k}-free (Lemma~\ref{lem:pointed}), but there exist \pt{k}s containing \kcross{(k+1)}s as well as \kcross{(k+1)}-free \kalter{k}-free sets of edges not contained in \pt{k}s (see \fref{fig:3crossing}).
Thus, we still miss a simple condition to characterize \mpt{}s of points in general position:

\begin{question}\label{ques:primal}
Characterize primals of \mpt{}s.
\end{question}

Another question related to the primal of a \mpt{} would be to determine whether for every point set there exists a \mpt{} that looks as simple as possible.
When $k=1$, we know that every point set in general position has:
\begin{enumerate}[(i)]
\item a pointed pseudotriangulation consisting only of triangles and four-sided pseudotriangles \cite{kkmsst-tdbptp-03};
\item a pointed pseudotriangulation whose maximum degree is at most five \cite{kkmsst-tdbptp-03}.
\end{enumerate}

\begin{question}
Does every point set in general position have a \pt{k} with only ``simple" stars (resp.~only vertices with ``little" degree)?
\end{question}

In the previous question, ``simple" may be interpreted either as ``with exactly $2k+1$ corners" or as ``with at most $2k+t$ edges" (for a minimal $t$). Similarily, ``little" means smaller than a constant (as small as possible).

\medskip

\paragraph{\sc Diameter of the graph of flips.}

The graph of flips on \pt{k}s of a pseudoline arrangement $L$ is connected, and an easy inductive argument shows that its diameter at is most quadratic in $|L|$. For certain particular cases, even better bounds are known:
\begin{enumerate}[(a)]
\item The diameter of the graph of flips on pointed pseudotriangulations of a set of $n$ points is at most~$O(n\ln n)$ \cite{b-tpt-04}.
\item The graph of flips on the \ktri{k}s of the \gon{n} has diameter bounded by~$2k(n-2k-1)$ \cite{n-gdfcp-00,ps-mtcsp-09}.
\item For triangulations of the \gon{n}, the diameter is exactly~$2n-10$ \cite{stt-rdthg-88}.
\end{enumerate}

\begin{question}
What is the (asymptotic) diameter of the flip graph?
\end{question}

\begin{figure}
	\centerline{\includegraphics[scale=1]{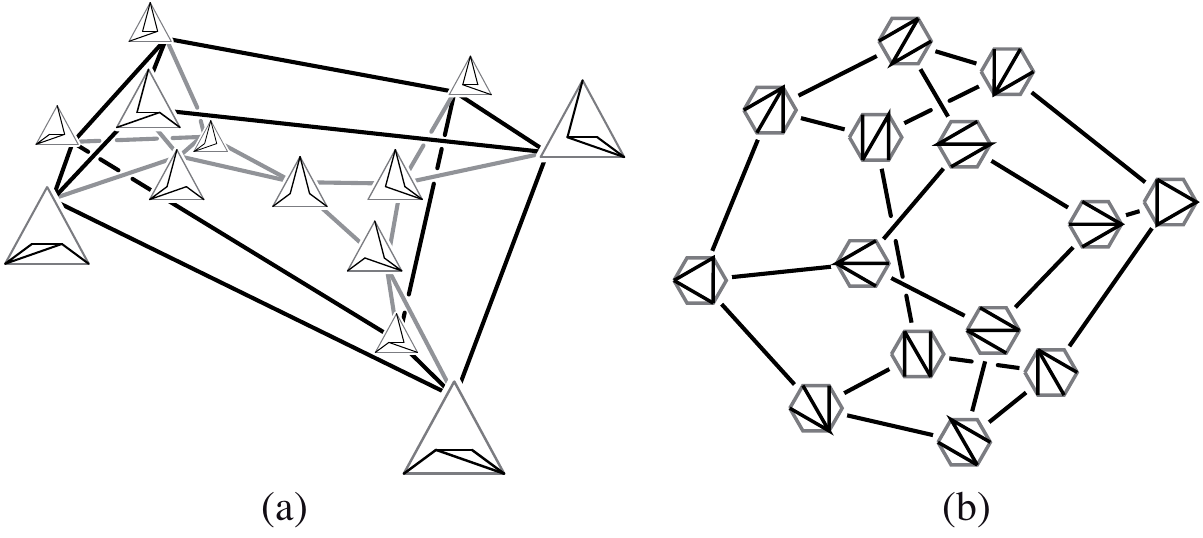}}
	\caption[The polytope of pointed pseudotriangulations and the associahedron]{(a) The polytope of pointed pseudotriangulations of a set of 5 points. (b) The $6$-dimensional associahedron.}
	\label{fig:associahedron}
\end{figure}

\medskip

\paragraph{\sc Polytopality.}

Let $\cS$ be the support of a pseudoline arrangement. Let $\Delta(\cS)$ denote the simplicial complex of subsets of contact points of pseudoline arrangements supported by~$\cS$. Our results ensure that $\Delta(\cS)$ is an abstract polytope whose ridge graph is the graph of flips (see the discussion in~\cite[Subsection 2.2]{bkps-ceppgfa-06}). When $\cS$ is the first kernel of the dual pseudoline arrangement of a point set of the Euclidean plane, it turns out that this abstract polytope can be realized effectively as a polytope of $\mathbb{R}^d$ (where $d$ is the number of flippable edges), which is the polar of the \emph{polytope of pointed pseudotriangulations} of \cite{rss-empppt-03}. An example of this polytope is presented in \fref{fig:associahedron}(a). When the points are in convex position, this polytope is the \emph{associahedron} (see \fref{fig:associahedron}(b)). This leads to the following question:

\begin{question}\label{qu:polytope}
Is $\Delta(\cS)$ the boundary complex of a polytope?
\end{question}

This question specializes for pseudotriangulations and multitriangulations to the following interesting questions:
\begin{enumerate}[(a)]
\item Is the graph of flips on pseudotriangulations of a point set polytopal? Since~\cite{rss-empppt-03} answers positively for Euclidean point sets, this question only remains open for non-stretchable arrangements. We have represented in \fref{fig:nonPappus} a pseudotriangulation of the smallest non-stretchable simple pseudoline arrangement (the non-Pappus pseudoline arrangement). Since this arrangement is symmetric under the dihedral group~$D_6$, it could be worked out with methods similar to those in~\cite{bp-srsc3cf-09}.

\begin{figure}
	\centerline{\includegraphics[scale=.8]{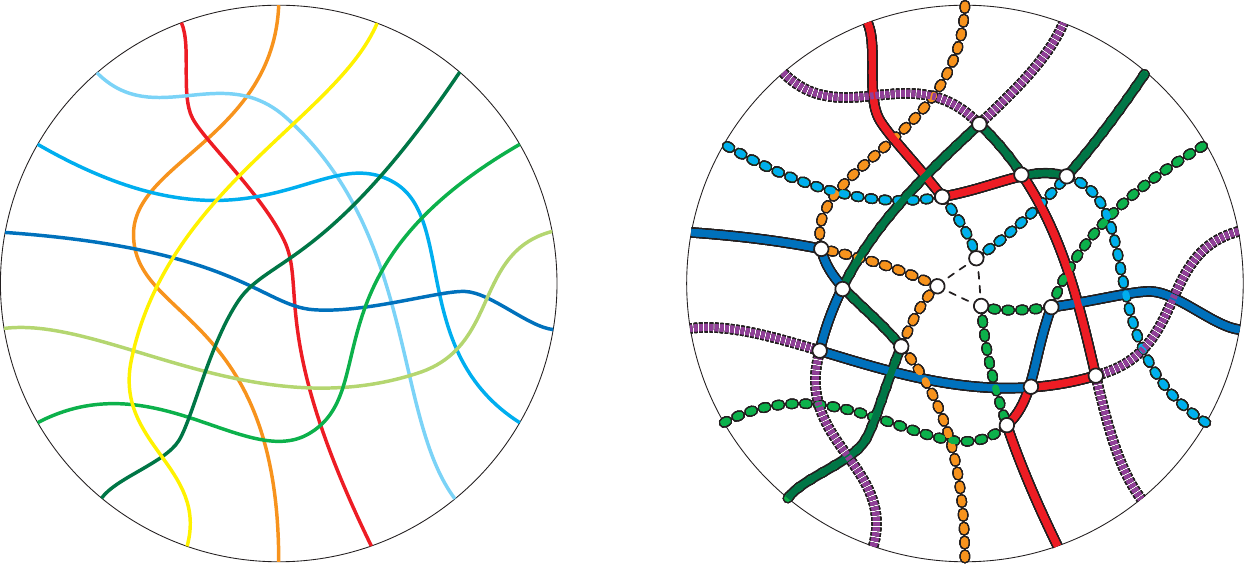}}
	\caption[The non-Pappus configuration]{A non-stretchable arrangement of $9$ pseudolines (left) represented on the projective plane to show its symmetries. Remove the center point to obtain the arrangement in the M\"obius strip. A pseudotriangulation of it (right).}
	\label{fig:nonPappus}
\end{figure}

\item Is the graph of flips on \ktri{k}s of the \gon{n} polytopal? Jonsson proved that the simplicial complex of \kcross{(k+1)}-free sets of \krel{k} diagonals of the \gon{n} is a combinatorial sphere~\cite{j-gt-03}. However, except for little cases, the question of the polytopality of this complex remains open. We refer to \cite{ps-mtcsp-09,bp-srsc3cf-09} and \cite[Section 4.3]{p} for a detailed discussion on this question.
\end{enumerate}

\begin{remark}
Since the submission of this paper, our knowledge on this question has progressed. In~\cite{s-npkt-11}, Stump made the connection between the multitriangulations and the type~$A$ subword complexes defined in~\cite{km-sccg-04}. These simplicial complexes are precisely the simplicial complexes~$\Delta(\cS)$ defined above. It implies in particular that
\begin{itemize}
\item these simplicial complexes are either topological spheres or topological balls~\cite{km-sccg-04}, and
\item Question~\ref{qu:polytope} is the special type~$A$ case of Question~6.4 in~\cite{km-sccg-04}.
\end{itemize}
We close this discussion by mentioning the related construction in~\cite{ps-bp-12}. It associates to each sorting network~$\cN$ its so-called \defn{brick polytope} whose vertices correspond to certain pseudoline arrangements with contact points supported by~$\cN$. For certain well-chosen networks, the brick polytope specializes to Hohlweg and Lange's realizations of the associahedron~\cite{hl-rac-07}. This construction answers Question~\ref{qu:polytope} positively for a certain class of supports~$\cS$. It was moreover extended recently to spherical subword complexes of any finite type in~\cite{ps-bpssc-11}.
\end{remark}

\medskip

\paragraph{\sc Number of \mpt{}s.}

In his paper~\cite{j-gtdfssp-05}, Jonsson proved that the number of \ktri{k}s of the \gon{n} is equal to the determinant $\det(C_{n-i-j})_{1\le i,j\le k}$ (where $C_m=\frac{1}{m+1}{2m\choose m}$ denotes the $m$th Catalan number). This determinant also counts non-crossing $k$-tuples of Dyck paths of semi-length $n-2k$ (see \cite{ps-mtcsp-09} and \cite[Section 4.1]{p} for a more detailed discussion). It raises the following question:

\begin{question}
Find an explicit bijection between Dyck multipaths and multitriangulations.
\end{question}

\begin{figure}[h]
	\capstart
	\centerline{\includegraphics[scale=1]{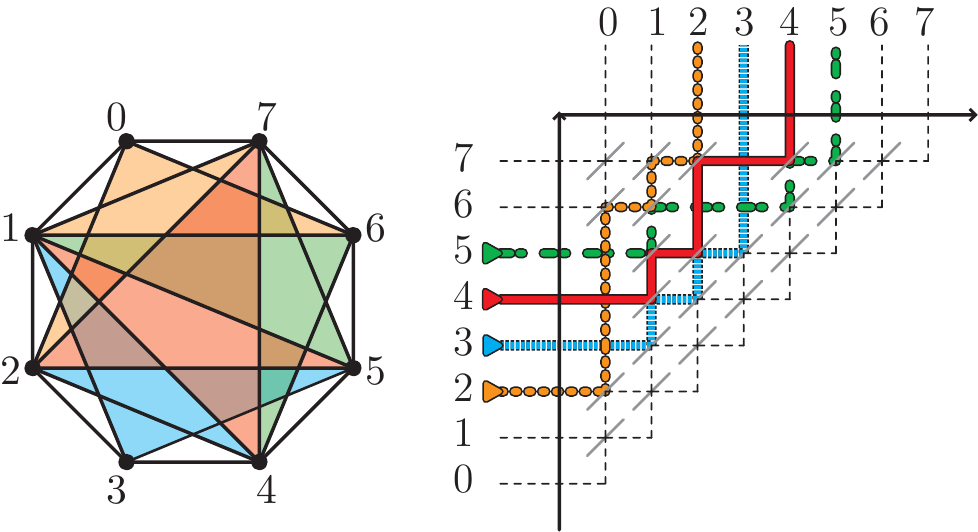}}
	\caption[Beam arrangement of a \ktri{k}]{The beam arrangement of a \ktri{2} of the octagon. Each pseudoline is a lattice path.}
	\label{fig:mirrorslasers}
\end{figure}

Our interpretation of multitriangulations in terms of pseudoline arrangements naturally associates to a \ktri{k} of the \gon{n} a set of $n-2k$ lattice paths as follows. Let~$T$ be a \ktri{k} of the \gon{n}. For each edge~$(u,v)$ of~$T$ (with~$u<v$), we place a \defn{mirror} at the grid point of coordinates~$(u,v)$. This mirror is a double faced mirror parallel to the diagonal~$x=y$ so that it reflects a ray coming from~$(-\infty,v)$ to a ray going to~$(u,+\infty)$, and a ray coming from~$(u,-\infty)$ to a ray going to~$(+\infty,v)$. Furthermore, for~$1\le i\le n-2k$, we place a \defn{light beam} at~$(-\infty,k-1+i)$ pointing horizontally. We obtain~$n-2k$ beams which reflect on the  mirrors of~$T$ (see \fref{fig:mirrorslasers}).

The results of this paper imply the following properties:
\begin{enumerate}[(i)]
\item All beams are $x$- and $y$-monotone lattice paths.
\item The $i$th beam comes from the direction $(-\infty,k-1+i)$ and goes to the direction $(k-1+i,+\infty)$.
\item Each beam reflects exactly~$2k+1$ times, and thus, has~$k$ vertical segments (plus one vertical half-line) and~$k$ horizontal segments (plus one horizontal half-line).
\item The beams form a pseudoline arrangement: any two of them cross exactly once.
\end{enumerate}

The $i$th beam~$B_i$ ``corresponds'' via duality to the \kstar{k}~$S_i$ of~$T$ whose $(k+1)$th vertex is the vertex~$k-1+i$ (in other words, the \kstar{k} bisected by the line passing through the vertex~$k-1+i$ and through the midpoint of~$[0,n-1]$). Indeed:
\begin{enumerate}[(i)]
\item the beam~$B_i$ is (by duality) the set of all bisectors of~$S_i$;
\item the mirrors which reflect~$B_i$ are the edges of~$S_i$;~and
\item the intersection of two beams~$B_i$ and~$B_j$ is the common bisector of~$S_i$ and~$S_j$.
\end{enumerate}
Observe that instead of $k$~Dyck paths of semi-length~$n-2k$, the beam arrangement of a \ktri{k} has~$n-2k$ beams which all have~$k$ horizontal steps. 

\begin{remark}
In their recent paper~\cite{ss-mfmpscsp-10}, Serrano and Stump also observe the correspondence between multitriangulations and beam arrangements (called ``reduced pipe dream'' in their paper). Starting from this correspondence, they provide an explicit bijection between multitriangulations and Dyck multipaths. We refer to~\cite{s-npkt-11,ss-mfmpscsp-10,r-mfmprclg-10} and the references therein.
\end{remark}

\medskip
Another interesting question concerning the number of multipseudo\-triangulations would be to determine what point sets give the maximum and minimum number of \mpt{}s. For example, when $k=1$, every point set in general position has at least as many pointed pseudotriangulations as the convex polygon with the same number of points~\cite{aaks-cmpt-04}.

\begin{question}
What point sets have the maximum and minimum number of \mpt{}s?
\end{question}

\medskip

\paragraph{\sc Decomposition of a \pt{k}.}

We have seen in Section~\ref{sec:iterated} that certain \mpt{}s can be decomposed into iterated \mpt{}s. Remember however that there exists irreductible \mpt{}s (Example~\ref{exm:irred}).

\begin{question}\label{q:decomp}
Characterize (completely) decomposable \mpt{}s.
\end{question}

The same question can be asked in a more algorithmical flavor:

\begin{question}\label{q:decompalgo}
How can we decide algorithmically whether a \pt{k} is decomposable?
\end{question}

Obviously, a brute-force algorithm is not considered as a good solution. A good way to test efficiency of the answers to Questions~\ref{q:decomp} and~\ref{q:decompalgo} would be to prove/disprove that Example~\ref{exm:irred} is the minimal irreducible \ktri{2} (or, in other words, that any \ktri{2} of an \gon{n}, with $n\le 14$, contains a triangulation).

\smallskip

Remember also that when a \mpt{} is decomposable, the graph of partial flips is not necessarily connected. In particular, finding all decompositions of a \mpt{} cannot be achieved just by searching in the graph of partial flips.

\begin{question}\label{q:alldecompalgo}
How can we enumerate all the decompositions of a \mpt{}?
\end{question}

\medskip

\paragraph{\sc Computing a \pt{k}.}

An initial pseudotriangulation of a set of $n$ points can be computed in $O(n\ln n)$ time, using only the predicate of the chirotope. A similar result would be interesting for \pt{k}s:

\begin{question}
Compute an initial \pt{k} of a given (double) pseudoline arrangement, using only its chirotope.
\end{question}


\section*{Acknowledgments}

We thank Francisco Santos for fruitful discussions on the subject (especially for pointing out to us the counter-example of \fref{fig:15gonctrexm}). During his visit at the \'Ecole Normale Sup\'erieure, Micha Sharir contributed a lot to the improvement of the paper. We thank Jakob Jonsson for pointed out an error in a previous formulation of Lemma~\ref{lem:alternationfree}. Finally, we are grateful to an anonymous reviewer for valuable feedback on the presentation of a preliminary version of this paper.


\bibliographystyle{alpha}
\bibliography{biblio.bib}

\end{document}